\documentclass[11pt]{article}
\usepackage{color}
\usepackage{srcltx}
\usepackage{graphicx}
\usepackage{amssymb,amsfonts,amsthm}
\usepackage{amsmath,amsthm,amssymb}
\usepackage{amsfonts}
\newcommand{\bbb}{{\mathcal B}}

\newcommand{\p}{{\mathfrak p}}
\newcommand{\z}{{\mathfrak z}}
\newcommand{\q}{{\mathfrak q}}
\newcommand{\g}{{\mathfrak g}}

\newcommand{\x}{{\mathfrak x}}
\newcommand{\uu}{{\mathfrak u}}
\newcommand{\w}{{\mathfrak w}}
\newcommand{\vv}{{\mathfrak v}}
\newcommand{\y}{{\mathfrak y}}
\newcommand{\dist}{\mathrm{dist}}

\newcommand{\Area}{\mathrm{Area}}
\newcommand{\irr}{\mathrm{curv}}

\newcommand{\ctt}{\{p,q\}}

\newcommand{\R}{{\mathbb R}}

\newtheorem{theorem}{Theorem}[section]
\newtheorem{lemma}[theorem]{Lemma}

\theoremstyle{definition}

\newtheorem{rk}[theorem]{Remark}

\newcounter{ppp}

\newcommand{\la}{\langle}
\newcommand{\ra}{\rangle}

\newcommand{\area}{\mathrm{area}}

\newcommand{\Z}{{\mathbb Z}}

\newcommand{\ccc}{{\mathcal C}}

\newcommand{\iv}{^{-1}}

\begin{document}

\renewcommand{\theequation}{\thesection.\arabic{equation}}

\title{On flat submaps of maps of non-positive curvature}
\author{A.Yu. Olshanskii, M.V. Sapir\thanks{Both authors were supported in part by the NSF grant DMS 1418506.The first author is also supported by RFFI grant 15-01-05823}}
\date{}
\maketitle

\begin{abstract} We prove that for every $r>0$ if a non-positively curved $(p,q)$-map $M$ contains no flat submaps of radius $r$, then the area of $M$ does not exceed $Crn$ for some constant $C$. This strengthens a theorem of Ivanov and Schupp. We show that an infinite $(p,q)$-map which tessellates the plane is quasi-isometric to the Euclidean plane if and only if the map contains only finitely many non-flat vertices and faces. We also generalize Ivanov and Schupp's result to a much larger class of maps, namely to maps with angle functions.
\end{abstract}

\tableofcontents

\section{Introduction}\label{s:i}

Recall that a \emph{map} is a finite, connected and simply-connected
2-complex embedded in the Euclidean plane. So its $1$-skeleton is a finite, connected plane graph.
The cells of dimensions $0$, $1$ and $2$ are called {\it vertices, edges} and {\it faces}, respectively.
Every edge $e$ has an orientation; so it starts at the vertex $e_-$ and ends at $e_+$, and $(e^{-1})_-=e_+$,
$(e^{-1})_+=e_-$ for the inverse edge $e^{-1}$, which has the same support as $e$.
The {\it degree} $d(o)$ of a vertex $o$ is the number of oriented edges $e$ with $e_-=o$. In particular, every loop $e$ (an edge which connects a vertex $o$ with itself) together with $e^{-1}$ contributes 2 to the degree of o.

If a closed path $\q=e_1\dots e_k$ is the boundary of a face $\Pi$, then the {\it degree} $d(\Pi)$
of $\Pi$ is the length $|\q| = k$. In particular, if both $e$ and $e^{-1}$ occur in the boundary path of a face, they contribute 2 to the degree of that face. Similarly one defines the {\it perimeter} $|\partial M|$ of a map $M$
as the length of a closed boundary path of $M$.

A \emph{submap} $N$ of a map $M$ is the subcomplex bounded by a closed  curve which can be made simple by an arbitrary small transformation. So either $N$ is a map or it can be turned into a map after such small transformation.

In group theory, maps appear most often as van Kampen diagrams. Many algebraic and geometric results about groups (say, the small cancelation theory, and construction of various groups with ``extreme properties" such as Tarski monsters \cite{LS, Olbook}) are obtained by establishing combinatorial properties of corresponding maps and their submaps. A typical example of such a statement: The area (i.e., the number of faces) of every reduced van Kampen diagram over a finite group presentation is at most linear in terms of its perimeter if and only if  the group is hyperbolic. In other words, hyperbolic groups are precisely the finitely presented groups with linear Dehn functions. One of recurrent features of van Kampen diagrams is existence of ``special" submaps in every van Kampen diagram of large area. For example, in the proof of the upper bound of the Dehn function of a group constructed in \cite{SBR} using an $S$-machine, it is proved that if the area of a reduced diagram is large enough, then up to a homotopy which does not change the area very much, the area is ``concentrated" is a few special subdiagrams called ``discs" (these are the subdiagrams simulating the work of the $S$-machine).

A remarkable result of this kind was proved by Ivanov and Schupp in \cite{IS}.
Recall that an edge $e$ of a map $M$ is called \emph{exterior} if it belongs to a boundary path of $M$. A face  $\Pi$ of $M$ is called {\it exterior} if its boundary $\partial\Pi$ has a common edge with $\partial(M)$.
An {\it exterior vertex} is one of the vertices of the boundary path. Non-exterior faces, vertices and edges are called \emph{interior}.


A map $M$ is called a $(p,q)$-map if every interior face $\Pi$ in $M$ has degree at least $p$
and
the degree of every interior vertex is at least $q$.
Note that if a group presentation ${\mathcal P}=\la X\ R\ra$ satisfies the small cancelation condition $C(p)-T(q)$, then every reduced van Kampen diagram over ${\mathcal P}$ is a $(p,q)$-map if we ignore all interior vertices of degree 2 (as in \cite{LS, IS}).
It is well known (see \cite{LS}) that if $\frac 1p+\frac 1q$ is smaller than $\frac 12$ (i.e., the curvature of the presentation is negative), then the group is hyperbolic and its Dehn function is linear. The case when $\frac 1p+\frac 1q>\frac 12$ is not interesting. Indeed, by a result of Gol'berg \cite{Gol, Olbook}, every group can be given by a presentation satisfying $C(5)$ and $T(3)$ and by a presentation satisfying $C(3)$ and $T(5)$ and hence by a presentation satisfying $C(p)$ and $T(q)$ for every $p\ge 3,q\ge 3$ with $\frac 1p +\frac 1q>\frac 12$ (see \cite{Gol, Olbook}). If $\frac 1p + \frac 1q=\frac 12$ (i.e., the ``curvature" is non-positive), and so $(p,q)$ is either $(3,6)$, $(4,4)$ or $(6,3)$), then the group has at most quadratic Dehn function \cite[Theorem V.6.2]{LS}.

A submap of a $(p,q)$-map is called \emph{flat} if each of its faces is {\it flat}, i.e. has degree $p$, and each interior vertex is {\it flat}, i.e. has degree $q$. The \emph{radius} of a map is the maximal distance from a  vertex to the boundary of the map.

Ivanov and Schupp \cite{IS} proved that if a $(p,q)$-map  $M$, $\frac 1p+\frac 1q=\frac 12$, has no flat submaps of radius $r$ (they call flat submaps \emph{regular}), then the area of the map is linear in terms of its perimeter with the multiplicative constant depending on $r$. More precisely, Ivanov and Schupp proved the following

\begin{theorem} [Ivanov, Schupp, \cite{IS}]\label{t:is} Let $M$ be a finite $(p,q)$-map with perimeter $n$ such that the maximal distance from a vertex in $M$ to a boundary vertex or to a non-flat vertex or face is $r$. Then the area of $M$ does not exceed $L(r)n$, where $L(r)$ is some exponential function in $r$.
\end{theorem}


Theorem \ref{t:is} implies that if a group presentation ${\mathcal P}=\la X\ R\ra$ satisfies conditions $C(p)$ and $T(q)$, $\frac 1p+\frac 1q=\frac 12$, and the radius of  every flat van Kampen diagram over ${\mathcal P}$ does not exceed certain constant, then the Dehn function of the group given by ${\mathcal P}$ is linear, hence the group is hyperbolic. Using this, Ivanov and Schupp proved hyperbolicity of many 1-related groups.

In this paper, we strengthen Theorem \ref{t:is} in two ways. First, we replace the exponential upper bound for $L(r)$ by a linear upper bound. Second, we extend Theorem \ref{t:is} to a much larger class of maps called ``maps with angle functions".

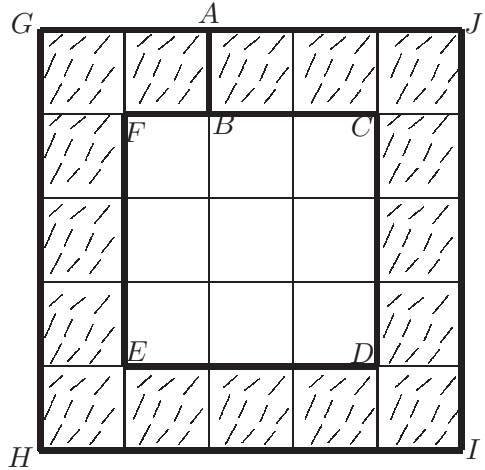
\begin{figure}[ht]
\begin{center}
\unitlength .8mm 
\linethickness{0.4pt}
\ifx\plotpoint\undefined\newsavebox{\plotpoint}\fi 
\begin{picture}(99,77.5)(0,0)
\put(27,63.5){\line(0,-1){59.25}}
\put(86,74.75){\line(-1,0){59.25}}
\put(41,63.5){\line(0,-1){59.25}}
\put(86,60.75){\line(-1,0){59.25}}
\put(55,63.5){\line(0,-1){59.25}}
\put(86,46.75){\line(-1,0){59.25}}
\put(69,63.5){\line(0,-1){59.25}}
\put(86,32.75){\line(-1,0){59.25}}
\put(83,63.5){\line(0,-1){59.25}}
\put(86,18.75){\line(-1,0){59.25}}
\put(97,63.5){\line(0,-1){59.25}}
\put(86,4.75){\line(-1,0){59.25}}
\put(27,74.5){\line(0,-1){13.5}}
\put(41,61){\line(0,1){13.75}}
\put(55,61){\line(0,1){13.75}}
\put(69,61){\line(0,1){13.75}}
\put(83,60.75){\line(0,1){14}}
\put(83,74.75){\line(1,0){14}}
\put(83,60.75){\line(1,0){14.25}}
\put(83,46.75){\line(1,0){14}}
\put(83,32.75){\line(1,0){14}}
\put(83,18.75){\line(1,0){14}}
\put(83,4.75){\line(1,0){14}}
\put(97,60.75){\line(0,1){14}}
\linethickness{2 pt}
\put(55,74.75){\line(1,0){42}}
\put(97,74.75){\line(0,-1){69.75}}
\put(96.75,4.75){\line(-1,0){69.75}}
\put(27,4.75){\line(0,1){70}}
\put(27,74.75){\line(1,0){28.25}}
\put(55,74.75){\line(0,-1){14}}
\put(55,60.75){\line(1,0){28}}
\put(83,60.75){\line(0,-1){42}}
\put(83,18.75){\line(-1,0){42}}
\put(41,18.75){\line(0,1){42}}
\put(41,60.75){\line(1,0){28.25}}
\linethickness{.1 pt}
\multiput(30.75,74)(-.0375,-.03333333){60}{\line(-1,0){.0375}}
\multiput(30.75,60.75)(-.0375,-.03333333){60}{\line(-1,0){.0375}}
\multiput(87,60.75)(-.0375,-.03333333){60}{\line(-1,0){.0375}}
\multiput(30.75,45.75)(-.0375,-.03333333){60}{\line(-1,0){.0375}}
\multiput(87,45.75)(-.0375,-.03333333){60}{\line(-1,0){.0375}}
\multiput(30.75,30.75)(-.0375,-.03333333){60}{\line(-1,0){.0375}}
\multiput(87,30.75)(-.0375,-.03333333){60}{\line(-1,0){.0375}}
\multiput(30.75,17.25)(-.0375,-.03333333){60}{\line(-1,0){.0375}}
\multiput(45.5,74)(-.0375,-.03333333){60}{\line(-1,0){.0375}}
\multiput(45.5,17.25)(-.0375,-.03333333){60}{\line(-1,0){.0375}}
\multiput(59.75,74)(-.0375,-.03333333){60}{\line(-1,0){.0375}}
\multiput(59.75,17.25)(-.0375,-.03333333){60}{\line(-1,0){.0375}}
\multiput(74,74)(-.0375,-.03333333){60}{\line(-1,0){.0375}}
\multiput(74,17.25)(-.0375,-.03333333){60}{\line(-1,0){.0375}}
\multiput(88.25,74)(-.0375,-.03333333){60}{\line(-1,0){.0375}}
\multiput(88.25,17.25)(-.0375,-.03333333){60}{\line(-1,0){.0375}}
\multiput(34.5,73.5)(-.04166667,-.03333333){60}{\line(-1,0){.04166667}}
\multiput(34.5,60.25)(-.04166667,-.03333333){60}{\line(-1,0){.04166667}}
\multiput(90.75,60.25)(-.04166667,-.03333333){60}{\line(-1,0){.04166667}}
\multiput(34.5,45.25)(-.04166667,-.03333333){60}{\line(-1,0){.04166667}}
\multiput(90.75,45.25)(-.04166667,-.03333333){60}{\line(-1,0){.04166667}}
\multiput(34.5,30.25)(-.04166667,-.03333333){60}{\line(-1,0){.04166667}}
\multiput(90.75,30.25)(-.04166667,-.03333333){60}{\line(-1,0){.04166667}}
\multiput(34.5,16.75)(-.04166667,-.03333333){60}{\line(-1,0){.04166667}}
\multiput(49.25,73.5)(-.04166667,-.03333333){60}{\line(-1,0){.04166667}}
\multiput(49.25,16.75)(-.04166667,-.03333333){60}{\line(-1,0){.04166667}}
\multiput(63.5,73.5)(-.04166667,-.03333333){60}{\line(-1,0){.04166667}}
\multiput(63.5,16.75)(-.04166667,-.03333333){60}{\line(-1,0){.04166667}}
\multiput(77.75,73.5)(-.04166667,-.03333333){60}{\line(-1,0){.04166667}}
\multiput(77.75,16.75)(-.04166667,-.03333333){60}{\line(-1,0){.04166667}}
\multiput(92,73.5)(-.04166667,-.03333333){60}{\line(-1,0){.04166667}}
\multiput(92,16.75)(-.04166667,-.03333333){60}{\line(-1,0){.04166667}}
\multiput(38.5,73.75)(-.03333333,-.04){75}{\line(0,-1){.04}}
\multiput(38.5,60.5)(-.03333333,-.04){75}{\line(0,-1){.04}}
\multiput(94.75,60.5)(-.03333333,-.04){75}{\line(0,-1){.04}}
\multiput(38.5,45.5)(-.03333333,-.04){75}{\line(0,-1){.04}}
\multiput(94.75,45.5)(-.03333333,-.04){75}{\line(0,-1){.04}}
\multiput(38.5,30.5)(-.03333333,-.04){75}{\line(0,-1){.04}}
\multiput(94.75,30.5)(-.03333333,-.04){75}{\line(0,-1){.04}}
\multiput(38.5,17)(-.03333333,-.04){75}{\line(0,-1){.04}}
\multiput(53.25,73.75)(-.03333333,-.04){75}{\line(0,-1){.04}}
\multiput(53.25,17)(-.03333333,-.04){75}{\line(0,-1){.04}}
\multiput(67.5,73.75)(-.03333333,-.04){75}{\line(0,-1){.04}}
\multiput(67.5,17)(-.03333333,-.04){75}{\line(0,-1){.04}}
\multiput(81.75,73.75)(-.03333333,-.04){75}{\line(0,-1){.04}}
\multiput(81.75,17)(-.03333333,-.04){75}{\line(0,-1){.04}}
\multiput(96,73.75)(-.03333333,-.04){75}{\line(0,-1){.04}}
\multiput(96,17)(-.03333333,-.04){75}{\line(0,-1){.04}}
\multiput(29.5,70.75)(-.03365385,-.07692308){52}{\line(0,-1){.07692308}}
\multiput(29.5,57.5)(-.03365385,-.07692308){52}{\line(0,-1){.07692308}}
\multiput(85.75,57.5)(-.03365385,-.07692308){52}{\line(0,-1){.07692308}}
\multiput(29.5,42.5)(-.03365385,-.07692308){52}{\line(0,-1){.07692308}}
\multiput(85.75,42.5)(-.03365385,-.07692308){52}{\line(0,-1){.07692308}}
\multiput(29.5,27.5)(-.03365385,-.07692308){52}{\line(0,-1){.07692308}}
\multiput(85.75,27.5)(-.03365385,-.07692308){52}{\line(0,-1){.07692308}}
\multiput(29.5,14)(-.03365385,-.07692308){52}{\line(0,-1){.07692308}}
\multiput(44.25,70.75)(-.03365385,-.07692308){52}{\line(0,-1){.07692308}}
\multiput(44.25,14)(-.03365385,-.07692308){52}{\line(0,-1){.07692308}}
\multiput(58.5,70.75)(-.03365385,-.07692308){52}{\line(0,-1){.07692308}}
\multiput(58.5,14)(-.03365385,-.07692308){52}{\line(0,-1){.07692308}}
\multiput(72.75,70.75)(-.03365385,-.07692308){52}{\line(0,-1){.07692308}}
\multiput(72.75,14)(-.03365385,-.07692308){52}{\line(0,-1){.07692308}}
\multiput(87,70.75)(-.03365385,-.07692308){52}{\line(0,-1){.07692308}}
\multiput(87,14)(-.03365385,-.07692308){52}{\line(0,-1){.07692308}}
\multiput(32.25,69.75)(-.03289474,-.07236842){38}{\line(0,-1){.07236842}}
\multiput(32.25,56.5)(-.03289474,-.07236842){38}{\line(0,-1){.07236842}}
\multiput(88.5,56.5)(-.03289474,-.07236842){38}{\line(0,-1){.07236842}}
\multiput(32.25,41.5)(-.03289474,-.07236842){38}{\line(0,-1){.07236842}}
\multiput(88.5,41.5)(-.03289474,-.07236842){38}{\line(0,-1){.07236842}}
\multiput(32.25,26.5)(-.03289474,-.07236842){38}{\line(0,-1){.07236842}}
\multiput(88.5,26.5)(-.03289474,-.07236842){38}{\line(0,-1){.07236842}}
\multiput(32.25,13)(-.03289474,-.07236842){38}{\line(0,-1){.07236842}}
\multiput(47,69.75)(-.03289474,-.07236842){38}{\line(0,-1){.07236842}}
\multiput(47,13)(-.03289474,-.07236842){38}{\line(0,-1){.07236842}}
\multiput(61.25,69.75)(-.03289474,-.07236842){38}{\line(0,-1){.07236842}}
\multiput(61.25,13)(-.03289474,-.07236842){38}{\line(0,-1){.07236842}}
\multiput(75.5,69.75)(-.03289474,-.07236842){38}{\line(0,-1){.07236842}}
\multiput(75.5,13)(-.03289474,-.07236842){38}{\line(0,-1){.07236842}}
\multiput(89.75,69.75)(-.03289474,-.07236842){38}{\line(0,-1){.07236842}}
\multiput(89.75,13)(-.03289474,-.07236842){38}{\line(0,-1){.07236842}}
\multiput(35.75,69)(-.03289474,-.07236842){38}{\line(0,-1){.07236842}}
\multiput(35.75,55.75)(-.03289474,-.07236842){38}{\line(0,-1){.07236842}}
\multiput(92,55.75)(-.03289474,-.07236842){38}{\line(0,-1){.07236842}}
\multiput(35.75,40.75)(-.03289474,-.07236842){38}{\line(0,-1){.07236842}}
\multiput(92,40.75)(-.03289474,-.07236842){38}{\line(0,-1){.07236842}}
\multiput(35.75,25.75)(-.03289474,-.07236842){38}{\line(0,-1){.07236842}}
\multiput(92,25.75)(-.03289474,-.07236842){38}{\line(0,-1){.07236842}}
\multiput(35.75,12.25)(-.03289474,-.07236842){38}{\line(0,-1){.07236842}}
\multiput(50.5,69)(-.03289474,-.07236842){38}{\line(0,-1){.07236842}}
\multiput(50.5,12.25)(-.03289474,-.07236842){38}{\line(0,-1){.07236842}}
\multiput(64.75,69)(-.03289474,-.07236842){38}{\line(0,-1){.07236842}}
\multiput(64.75,12.25)(-.03289474,-.07236842){38}{\line(0,-1){.07236842}}
\multiput(79,69)(-.03289474,-.07236842){38}{\line(0,-1){.07236842}}
\multiput(79,12.25)(-.03289474,-.07236842){38}{\line(0,-1){.07236842}}
\multiput(93.25,69)(-.03289474,-.07236842){38}{\line(0,-1){.07236842}}
\multiput(93.25,12.25)(-.03289474,-.07236842){38}{\line(0,-1){.07236842}}
\multiput(30.5,65.75)(-.03365385,-.0625){52}{\line(0,-1){.0625}}
\multiput(30.5,52.5)(-.03365385,-.0625){52}{\line(0,-1){.0625}}
\multiput(86.75,52.5)(-.03365385,-.0625){52}{\line(0,-1){.0625}}
\multiput(30.5,37.5)(-.03365385,-.0625){52}{\line(0,-1){.0625}}
\multiput(86.75,37.5)(-.03365385,-.0625){52}{\line(0,-1){.0625}}
\multiput(30.5,22.5)(-.03365385,-.0625){52}{\line(0,-1){.0625}}
\multiput(86.75,22.5)(-.03365385,-.0625){52}{\line(0,-1){.0625}}
\multiput(30.5,9)(-.03365385,-.0625){52}{\line(0,-1){.0625}}
\multiput(45.25,65.75)(-.03365385,-.0625){52}{\line(0,-1){.0625}}
\multiput(45.25,9)(-.03365385,-.0625){52}{\line(0,-1){.0625}}
\multiput(59.5,65.75)(-.03365385,-.0625){52}{\line(0,-1){.0625}}
\multiput(59.5,9)(-.03365385,-.0625){52}{\line(0,-1){.0625}}
\multiput(73.75,65.75)(-.03365385,-.0625){52}{\line(0,-1){.0625}}
\multiput(73.75,9)(-.03365385,-.0625){52}{\line(0,-1){.0625}}
\multiput(88,65.75)(-.03365385,-.0625){52}{\line(0,-1){.0625}}
\multiput(88,9)(-.03365385,-.0625){52}{\line(0,-1){.0625}}
\multiput(33.5,65)(-.03333333,-.0375){60}{\line(0,-1){.0375}}
\multiput(33.5,51.75)(-.03333333,-.0375){60}{\line(0,-1){.0375}}
\multiput(89.75,51.75)(-.03333333,-.0375){60}{\line(0,-1){.0375}}
\multiput(33.5,36.75)(-.03333333,-.0375){60}{\line(0,-1){.0375}}
\multiput(89.75,36.75)(-.03333333,-.0375){60}{\line(0,-1){.0375}}
\multiput(33.5,21.75)(-.03333333,-.0375){60}{\line(0,-1){.0375}}
\multiput(89.75,21.75)(-.03333333,-.0375){60}{\line(0,-1){.0375}}
\multiput(33.5,8.25)(-.03333333,-.0375){60}{\line(0,-1){.0375}}
\multiput(48.25,65)(-.03333333,-.0375){60}{\line(0,-1){.0375}}
\multiput(48.25,8.25)(-.03333333,-.0375){60}{\line(0,-1){.0375}}
\multiput(62.5,65)(-.03333333,-.0375){60}{\line(0,-1){.0375}}
\multiput(62.5,8.25)(-.03333333,-.0375){60}{\line(0,-1){.0375}}
\multiput(76.75,65)(-.03333333,-.0375){60}{\line(0,-1){.0375}}
\multiput(76.75,8.25)(-.03333333,-.0375){60}{\line(0,-1){.0375}}
\multiput(91,65)(-.03333333,-.0375){60}{\line(0,-1){.0375}}
\multiput(91,8.25)(-.03333333,-.0375){60}{\line(0,-1){.0375}}
\multiput(39.25,68.5)(-.03333333,-.04166667){60}{\line(0,-1){.04166667}}
\multiput(39.25,55.25)(-.03333333,-.04166667){60}{\line(0,-1){.04166667}}
\multiput(95.5,55.25)(-.03333333,-.04166667){60}{\line(0,-1){.04166667}}
\multiput(39.25,40.25)(-.03333333,-.04166667){60}{\line(0,-1){.04166667}}
\multiput(95.5,40.25)(-.03333333,-.04166667){60}{\line(0,-1){.04166667}}
\multiput(39.25,25.25)(-.03333333,-.04166667){60}{\line(0,-1){.04166667}}
\multiput(95.5,25.25)(-.03333333,-.04166667){60}{\line(0,-1){.04166667}}
\multiput(39.25,11.75)(-.03333333,-.04166667){60}{\line(0,-1){.04166667}}
\multiput(54,68.5)(-.03333333,-.04166667){60}{\line(0,-1){.04166667}}
\multiput(54,11.75)(-.03333333,-.04166667){60}{\line(0,-1){.04166667}}
\multiput(68.25,68.5)(-.03333333,-.04166667){60}{\line(0,-1){.04166667}}
\multiput(68.25,11.75)(-.03333333,-.04166667){60}{\line(0,-1){.04166667}}
\multiput(82.5,68.5)(-.03333333,-.04166667){60}{\line(0,-1){.04166667}}
\multiput(82.5,11.75)(-.03333333,-.04166667){60}{\line(0,-1){.04166667}}
\multiput(96.75,68.5)(-.03333333,-.04166667){60}{\line(0,-1){.04166667}}
\multiput(96.75,11.75)(-.03333333,-.04166667){60}{\line(0,-1){.04166667}}
\multiput(37,64.75)(-.03365385,-.04326923){52}{\line(0,-1){.04326923}}
\multiput(37,51.5)(-.03365385,-.04326923){52}{\line(0,-1){.04326923}}
\multiput(93.25,51.5)(-.03365385,-.04326923){52}{\line(0,-1){.04326923}}
\multiput(37,36.5)(-.03365385,-.04326923){52}{\line(0,-1){.04326923}}
\multiput(93.25,36.5)(-.03365385,-.04326923){52}{\line(0,-1){.04326923}}
\multiput(37,21.5)(-.03365385,-.04326923){52}{\line(0,-1){.04326923}}
\multiput(93.25,21.5)(-.03365385,-.04326923){52}{\line(0,-1){.04326923}}
\multiput(37,8)(-.03365385,-.04326923){52}{\line(0,-1){.04326923}}
\multiput(51.75,64.75)(-.03365385,-.04326923){52}{\line(0,-1){.04326923}}
\multiput(51.75,8)(-.03365385,-.04326923){52}{\line(0,-1){.04326923}}
\multiput(66,64.75)(-.03365385,-.04326923){52}{\line(0,-1){.04326923}}
\multiput(66,8)(-.03365385,-.04326923){52}{\line(0,-1){.04326923}}
\multiput(80.25,64.75)(-.03365385,-.04326923){52}{\line(0,-1){.04326923}}
\multiput(80.25,8)(-.03365385,-.04326923){52}{\line(0,-1){.04326923}}
\multiput(94.5,64.75)(-.03365385,-.04326923){52}{\line(0,-1){.04326923}}
\multiput(94.5,8)(-.03365385,-.04326923){52}{\line(0,-1){.04326923}}
\put(55,77.5){\makebox(0,0)[cc]{$A$}}
\put(57.5,59){\makebox(0,0)[cc]{$B$}}
\put(80.5,58.75){\makebox(0,0)[cc]{$C$}}
\put(80.5,21){\makebox(0,0)[cc]{$D$}}
\put(43,21.25){\makebox(0,0)[cc]{$E$}}
\put(42.75,57.75){\makebox(0,0)[cc]{$F$}}
\put(99,75.75){\makebox(0,0)[cc]{$J$}}
\put(99,5){\makebox(0,0)[cc]{$I$}}
\put(23.75,3.75){\makebox(0,0)[cc]{$H$}}
\put(24,75.75){\makebox(0,0)[cc]{$G$}}
\end{picture}

\caption{The thick path $ABCDEFBAGHIJA$ is simple up to an arbitrary small deformation. It bounds a non-simple submap.}
\label{f:0}
\end{center}

\end{figure}

Let us call a submap of a $(p,q)$-map \emph{simple} if it is bounded by a simple closed curve.
Note that the closure of the union of faces from a non-simple submap may not be simply connected (Figure \ref{f:0}) while the closure of the union of faces from a simple submap is always simply connected.

Let $p=4,3$ or $6$. Let $S^p$ be the usual tessellation of the plane by $p$-gons.  Then for every $n\ge 0$ the {\emph standard} map $S^p_n$ is constructed as follows. By definition $S^p_0$ is a vertex $o$ in $T$. If the submap $S^p_{n}$ of $T$ is constructed, then $S^p_{n+1}$ is the (closure of the) union of all faces having a common vertex with $S^p_n$. Then $S^p_n$ is a simple $(p,q)$-map. For example, then $S^4_n$ is the $2n\times 2n$-square tessellated by unit squares, $S^3_n$ is a regular hexagon with side length $n$ tessellated by equilateral triangles with side length 1, $S^6_n$ can be viewed as the weak dual\footnote{Recall that if $M$ is a map, then the weak dual map $\bar M$ is obtained by putting a vertex in every (bounded) face, and for every edge shared by two faces of $M$, connect the two vertices from these faces by an edge crossing that edge. Thus the vertices of $\bar M$ correspond to faces of $M$, edges of $\bar M$ correspond to interior edges of $M$, faces of $\bar M$ correspond to interior vertices of $M$.} to the $(3,6)$-map constructed just as  $S^3_{n+1}$, only
starting with a triangle face instead of a vertex.

\begin{rk} \label{r:01}
Every simple $(p,q)$-map, $\frac 1p+\frac 1q=\frac 12$, of radius $r$, contains a simple submap $M'$ which is isomorphic to $S^p_n$ for $n=O(r)$.  The submap $M'$ can be obtained in a similar manner as $S^p_n$.
Pick a vertex $o$ in $M$ at distance $r$ from the boundary of $M$. This is the submap $M_0$. If $M_i$ is already constructed, then $M_{i+1}$ is obtained from $M_i$ by adding all faces having a vertex in common with $M_i$. The process continues until one of the vertices in $M_i$ is exterior. In that case we set $M'=M_i$. Of course it should be explained why $M_{i+1}$ is indeed a simple standard submap provided $M_i$ is already a simple standard submap. It is not as obvious as it seems. The explanation uses Lemma \ref{l:0} below, it is given in Remark \ref{st}.
\end{rk}

Our main result is the following

\begin{theorem}[See Theorem  \ref{main} below] If a $(p,q)$-map $M$ does not contain flat simple submaps of radius $\ge r$, then the area of $M$ is at most $crn$ for some constant $c$.
\end{theorem}


Note that the statement of Theorem \ref{main} is non-trivial even for the van Kampen diagrams over the standard presentation of $\Z^2$ although there is a significantly easier proof in this case.
Theorem \ref{main} is applicable to van Kampen diagrams over any $C(p)-T(q)$-presentations with $\frac 1p+\frac 1q=\frac 12$, say, the standard presentations of 2-dimensional Right Angled Artin groups or the fundamental groups of alternating knots.

Theorem \ref{main} is proved in Section \ref{s:0}. The plan of the proof is the following. First for every map $M$ and every two (real) numbers $p,q$ with $\frac 1p+\frac 1q=\frac 12$ we define curvature of faces and vertices of $M$ as numbers proportional to the excessive degrees, and show that the sum of all curvatures is equal to $p$. Then we assume that $p,q$ are positive integers (so $(p,q)$ is $(3,6)$, $(4,4)$ or $(6,3)$) and
note that by a simple transformation of the map, we can assume that all faces of $M$ have degrees between $p$ and $2p-1$, the perimeter of the map after this transformation increases by a factor of $\le p-1$ and the set of (simple) flat submaps does not change. The key ``contraction" Lemma \ref{l:1} says that  the perimeter of the interior $M^0$ which is the union of all faces of $M$ having no boundary vertices of $M$, is ``substantially smaller" than the perimeter of $M$. From this, we deduce, first, that the area of $M$ is $O(Rn)$, where $R$ is the radius of $M$. Second, we deduce that one can cut $M$ along paths of linear in $n$ total length so that in the resulting map $\tilde M$ all non-flat vertices and faces are on the boundary. Then the radius $\tilde R$ of $\tilde M$ is less than $r+p$.
 Hence $\Area(M)=\Area(\tilde M) =  O(\tilde R \tilde  n) = O(rn)$.

In Section \ref{s:3} we consider infinite maps on the plane, i.e. tessellations of $\R^2$. An infinite map is called \emph{proper} if its support, i.e. the union of all faces, edges and vertices is the whole plane $\R^2$ and every disc in $\R^2$ intersects only a finite number of faces, edges and vertices of the map. Our main result is the following:

\begin{theorem}[See Theorem \ref{t:qi}]\label{t:qi1} Let $M$ be a proper $(p,q)$-map with $\frac 1p+\frac 1q=\frac 12$. Then the $1$-skeleton of $M$ with its path metric is quasi-isometric to the Euclidean plane if an only if $M$ has only finitely many non-flat vertices and faces.
\end{theorem}

Our proof of Theorem \ref{t:qi1} proceeds as follows. Suppose that a proper map $M$ has finite number non-flat vertices and faces. Then we modify it in a finite sequence of steps. At each step we reduce the number of non-flat faces and vertices while keeping the map quasi-isometric to $M$. As a result we get a proper map $M'$ with at most one non-flat vertex and no non-flat faces. Such a map is naturally subdivided by infinite paths emanating from the non-flat vertex into a finite number of convex infinite submaps, each of which is quasi-isometric to a quadrant of the Euclidean plane. Combining the corresponding quasi-isometries, we get a quasi-isometry between $M'$ and $\R^2$. The main tool in the proof is the notion of infinite \emph{corridor}, that is an infinite sequence of faces in $M$ where  each consecutive faces share an edge. This gives the ``if" part of the  theorem.

For the ``only if" part, we prove that if a proper map $M$ has infinitely many non-flat vertices or faces, then for every constant $c$ it contains an infinite $c$-separated set $S$ of vertices which has super-quadratic growth function (that is the function that for every $n$ gives the number of vertices from $S$ at distance $\le n$ from a given vertex is super-quadratic). This cannot happen if $M$ was quasi-isometric to $\R^2$. The key tool in proving this part of the theorem is the ``contraction" Lemma \ref{l:1} from the proof of Theorem \ref{main}. We construct a sequence of submaps $N(r)$ such that the boundaries of $N(r)$ contain large $c$-separated sets of vertices. In order to prove this property of $N(r)$ we use winding numbers of piece-wise geodesic paths in $M$ passing through vertices of $\partial(N(r))$ around a vertex which is deep inside $N(r)$.

Bruce Kleiner and Michah Sageev explained to us that  the ``only if" part of Theorem \ref{t:qi1} can be deduced  from Theorem 4.1 of their paper \cite{BKS} (joint with Mladen Bestvina).
If we view every face of a proper $(p,q)$-map $M$ as a regular Euclidean $n$-gone, then the  map $M$ turns into a CAT(0) 2-complex $M'$ which is quasi-isometric to the original map. Then Part 1 of \cite[Theorem 4.1]{BKS} implies that $M$ has a locally finite second homology class whose support $S$ is locally isometric to the Euclidean plane outside some finite ball. It remains to notice that the only such homology class is (up to a scalar multiple) the fundamental class of $M$. Hence $M$ is locally flat outside a finite ball. Bruce Kleiner also explained how to deduce the ``if" part of Theorem \ref{t:qi1} using Riemannian geometry. First we ``smooth out" the CAT(0) 2-complex $M$ which is locally flat outside a finite ball to obtain (using the Cartan-Hadamard theorem) a 2-dimensional Riemannian manifold $M'$ with the same property and which is quasi-isometric to the map $M$. Then we use the Rauch comparison theorem to establish a bi-Lipschitz equivalence between $M'$ and the Euclidean plane.

Note that our proof of Theorem \ref{t:qi1} is completely self-contained, short  and uses only basic graph theory.

In Section \ref{s:4} we consider the class of maps with angle functions. Let $o$ be a vertex on the boundary of a face $\Pi$ of a map $M$. A \emph{corner} of $\Pi$ at $o$ is the pair of two consecutive oriented edges $e$ and $f$ of $\partial(\Pi)$ with $e_+=f_-=o$. ($f^{-1}$ and $e^{-1}$ define the same corner.) An angle function assigns a non-negative number (\emph{angle}) to each corner of each cell. Then the \emph{curvature} of an interior vertex $o$ is $2\pi$ minus the sum of all angles at $o$. The curvature of an exterior vertex is defined in a similar way.  The curvature of a face of perimeter $d$ is the sum of angles of corners of this face minus the sum of angles of an Euclinean $d$-gone (that is $\pi(d-2)$).
A map with an angle function is called \emph{flat} if all its faces and interior vertices are of curvature 0. A map with an angle function is called a $(\delta, b)$-map, $\delta>0, b> 0$ if the curvature of every non-flat vertex and face does not exceed $-\delta$ and the degree of every vertex and face does not exceed $b$.

The class of $(\delta,b)$-maps is very large. By F\'ary's theorem (see \cite{K, F}), every finite planar graph $M$ without double edges and loops can be drawn on the Euclidean plane using only
straight line segments for edges. The proof from \cite{F} shows that if $M$ is a plane map, then one can assume that the graph with straight edges is isomorphic to $M$ as a 2-complex.
For a map with straight edges, we can assign to each corner its Euclidean angle, making the map flat.

Note that many authors considered van Kampen diagrams as maps with angle functions. Some of the earliest implementations of this idea are in the papers \cite{Gersten} by Steve Gersten, \cite{Pride} by Steve Pride  and  \cite{How} by Jim Howie.

It is easy to see that a (finitely presented) group $G$ has a (finite) presentation ${\mathcal P}=\la X\mid R\ra$ such that every van Kampen diagram over ${\mathcal P}$ can be assigned an angle function which makes the diagram a flat map if and only if one can find a finite generating set of the group which does not contain involutions. Such a finite generating set exists if and only if the group is not
an extension of an Abelian group $A$ by the automorphism of order $2$ which takes every element of $A$ to its inverse.

We will show in Remark \ref{last} that every $(p,q)$-map can be transferred into a $(\delta,b)$-map with angle function without decreasing the area, or increasing the perimeter or the set of flat vertices and faces. Thus the following theorem is a generalization of Ivanov and Schupp's Theorem \ref{t:is} to a much wider class of maps.

 \begin{theorem}\label{t:dense} Suppose that $M$ is a
$(\delta, b)$-map of  perimeter $n$ and the distance of every vertex of $M$ to a boundary vertex or to a non-flat vertex or face of $M$ is at most $r$. Then $\Area(M)\le L n$,
 where $L$ is exponential in $r$.
 \end{theorem}

The key part of the (very short) proof of Theorem \ref{t:dense} (see Section \ref{s:4}) is Lemma \ref{A} which shows that in every non-positively curved map with an angle function, the sum of curvatures of all faces and interior vertices exceeds $\pi(2-n)$ where $n$ is the perimeter of the map.


\section{Large flat submaps of $(p,q)$-maps}\label{s:0}


Although every edge has an orientation, when
counting the numbers of edges (or faces) in a map, we take usually any pair $(e, e^{-1})$ as one non-oriented edge
(e.g., $E$ is the number of non-oriented edges, when we apply Euler's formula).
The boundary path $\p$ of a map or a face is considered up to cyclic permutations and taking
inverse paths $\p^{-1}$.


The number of faces
in $M$ is called the {\it area} of $M$, denoted $\Area(M)$.

Here is a precise formulation of our main result.

\begin{theorem}\label{main} Let $p, q$ be positive integers with $\frac 1p+\frac 1q=\frac 12$,
$C=\frac32(p-1)(q+1)$.
Then the  area
of $M$ does not exceed $C(r+p)n$,
provided $M$ contains no simple
flat submaps of radius greater than $r$.
\end{theorem}

\subsection{A lemma about curvatures}

Given a pair $(p,q)$ of arbitrary (possibly negative) real numbers with $\frac{1}{p}+\frac{1}{q}=\frac {1}{2}$, the {\em curvature} of a face $\Pi$ in a  map $M$ is defined as $\irr_{p,q}(\Pi)=p-d(\Pi)$.  Let $o$ be a vertex in $M$. Then let $\mu(o)$ be the number of times the boundary path $\x$ goes through the vertex $o$.
For example,  the multiplicity $\mu(o)$ is $1$ if the closed path $\x$ passes through $o$ only once, and $\mu(o)=0$ if $o$ is an interior vertex.
The {\em curvature} $\irr_{p,q}(o)$ of a vertex $o$ is defined as $\frac pq (q-d(o))-\mu(o)$.  Let $I_v=I_v(M)$ be the sum of curvatures of all vertices of $M$, and let $I_f=I_f(M)$ be the sum of curvatures of all faces of $M$.


The following lemma follows from Theorem 3.1 from \cite[Capter V]{LS}, we include the proof for completeness and because it is significantly easier than in \cite{LS}.

\begin{lemma} \label{l:0} For an arbitrary map $M$ and arbitrary real $p,q$ with $\frac 1p +\frac 1q=\frac 12$, we have $I_v+I_f=p$.
 \end{lemma}

\proof The statement is obvious for a map consisting of one vertex. So assume that $M$ has more than one vertex, hence it has no vertices of degree $0$. Let us assign weight 1 to every non-oriented edge of the map $M$. Then the sum of all weights is the number $E$ of edges in $M$. Now let us make each edge give $\frac 1q$ of its weight to each of its vertices (if it has only one vertex, the edge gives it $\frac 2q$) and $\frac 1p$ to each of the (at most two) faces containing that edge. Thus the sum of weights of all vertices of $M$ is equal to $\sum_o \frac 1q d(o)$, where $o$ runs over all vertices of $M$. The sum of weights of all faces is $\sum_\Pi \frac 1p d(\Pi)$, where the sum runs over all faces of $M$.

By the assumption $\frac{1}{p}+\frac{1}{q}=\frac {1}{2}$, every edge separating two faces becomes completely
weightless (it gives $\frac 1p$ to each of the faces and $\frac 1q$ to each of its vertices). For the same reason, an edge $e$ of the boundary path $\x=\partial(M)$ becomes of weight $\frac{1}{p}$ if it lies on the boundary of a face, or $\frac{2}{p}$ otherwise. In the later case, the non-oriented edge $e$  occurs in
the path $\x$ twice (with different orientations).
Therefore after the redistribution of weights, the sum of weighs of all edges in $M$ is equal to $\frac 1p n$ where $n$ is the perimeter of the map.

Thus, the total weight is equal to $$E=\sum_o \frac 1q d(o) + \sum_\Pi \frac 1p d(\Pi)+\frac 1p n.$$ Since $E-V-F=-1$ by Euler's formula, where $V$ and $F$ are numbers of vertices and faces in $M$ respectively, we have:

\begin{equation}\label{1} -1=\sum_o \frac 1q d(o) + \sum_\Pi \frac 1p d(\Pi)+\frac 1p n - V - F=\sum_o \left(\frac 1q d(o)-1\right) + \sum_\Pi \left(\frac 1p d(\Pi)-1\right)+\frac 1p n.\end{equation}

Notice also that $n=\sum_v \mu(o)$ where $v$ runs over all vertices of $M$ (indeed, $\mu(o)=0$ for all interior vertices and $\mu(o)$ is the number of times the boundary path passes through $o$).
Therefore we can rewrite (\ref{1}) as follows:

$$-1=\sum_o \left(\frac 1q d(o)-1+\frac 1p \mu(o)\right) + \sum_\Pi \left(\frac 1p d(\Pi)-1\right). $$ Since the first of these sums is $-\frac 1p I_v$ and the second sum is $-\frac 1p I_f$, we deduce that
$I_v+I_f=p$.

\endproof

\begin{rk}\label{r:0} A result similar to Lemma \ref{l:0} is true for maps on arbitrary surfaces $S$ with boundary.  It is easy to see that in that case the right hand side is equal to $p\chi(S)$ where $\chi(S)$ is the Euler characteristic of the smallest subsurface of $S$ containing the map.
\end{rk}

\begin{rk} \label{st} Let us use Lemma \ref{l:0} to complete the proof from Remark \ref{r:01} about standard submaps of simple flat maps. Here we consider the case $p=q=4$ only, leaving the other two cases to the reader as an exercise.
Let, as in Remark \ref{r:01}, $M$ be a simple flat map of radius $r$. We set $M_0$ to be a vertex $o$ at distance $r$ from the boundary of $M$, and assume that for $i\ge 0$, the submap $M_i$ is constructed and this submap is simple and isomorphic to the standard map $S^4_i$ which is the $n\times n$-square (where $n=2i$)  tessellated by unit squares. Counting the difference between the degrees of the vertices from $\partial M_i$ in $M$ and in $M_i$, we obtain
exactly $4n+4$ oriented edges $e_1, ..., e_{4n+4}$ with $(e_j)_-\in \partial(M_i)$ and $(e_j)_+\not\in M_i$.
We claim that no two of the edges $e_j, e_k$, $j\ne k$, are mutually inverse and no two of the vertices $(e_j)_+, (e_k)_+$ coincide.

Indeed, suppose that  $e_j=e_k\iv$. Then this edge and a subpath of $\partial M_i$ bound a submap $N$
having no faces from $M_i$.
\begin{figure}[ht]
\begin{center}
\unitlength .7mm 
\linethickness{0.4pt}
\ifx\plotpoint\undefined\newsavebox{\plotpoint}\fi 
\begin{picture}(98,108.25)(0,0)
\put(5.25,19){\line(0,-1){12.5}}
\put(56.25,19){\line(0,-1){12.5}}
\put(5.25,31.5){\line(0,-1){12.5}}
\put(56.25,31.5){\line(0,-1){12.5}}
\put(5.25,44){\line(0,-1){12.5}}
\put(56.25,44){\line(0,-1){12.5}}
\put(5.25,56.5){\line(0,-1){12.5}}
\put(56.25,56.5){\line(0,-1){12.5}}
\put(18,19){\line(0,-1){12.5}}
\put(18,44){\line(0,-1){12.5}}
\put(18,56.5){\line(0,-1){12.5}}
\put(30.75,19){\line(0,-1){12.5}}
\put(30.75,44){\line(0,-1){12.5}}
\put(30.75,56.5){\line(0,-1){12.5}}
\put(43.5,19){\line(0,-1){12.5}}
\put(43.5,56.5){\line(0,-1){12.5}}
\put(17.75,6.5){\line(-1,0){12.5}}
\put(17.75,19){\line(-1,0){12.5}}
\put(17.75,31.5){\line(-1,0){12.5}}
\put(17.75,44){\line(-1,0){12.5}}
\put(30.5,6.5){\line(-1,0){12.5}}
\put(30.5,19){\line(-1,0){12.5}}
\put(30.5,31.5){\line(-1,0){12.5}}
\put(44.75,31.5){\line(-1,0){12.5}}
\put(43.25,6.5){\line(-1,0){12.5}}
\put(43.25,31.5){\line(-1,0){12.5}}
\put(43.25,44){\line(-1,0){12.5}}
\put(56,6.5){\line(-1,0){12.5}}
\put(56,19){\line(-1,0){12.5}}
\put(56,31.5){\line(-1,0){12.5}}
\put(56,44){\line(-1,0){12.5}}
\put(5.5,56.5){\line(1,0){12.5}}
\put(18.25,56.5){\line(1,0){12.5}}
\put(31,56.5){\line(1,0){12.5}}
\put(43.75,56.5){\line(1,0){12.5}}
\put(18,43.75){\line(0,-1){24.75}}
\put(18,19){\line(1,0){25.5}}
\put(43.5,19){\line(0,1){25}}
\put(43.5,44){\line(-1,0){25.5}}
\put(43.5,19){\line(1,0){12.75}}
\put(56.25,19){\line(0,-1){12.5}}
\put(56.25,6.5){\line(-1,0){12.75}}
\put(30.75,31.5){\line(0,-1){12.5}}
\qbezier(56.25,19.25)(98,108.25)(30.75,56.25)
\put(44.25,64.25){\line(0,1){0}}
\put(61.5,74){\line(0,1){0}}
\put(55,64.5){\line(0,1){0}}
\put(64.75,64){\line(0,1){0}}
\put(65.25,58){\line(0,1){0}}
\put(65,52.5){\line(0,1){0}}
\put(65,45.5){\line(0,1){0}}
\put(65.25,40.25){\line(0,1){0}}
\put(65.25,35.25){\line(0,1){0}}
\put(44.25,64.25){\line(0,1){0}}
\put(61.5,74){\line(0,1){0}}
\put(55,64.5){\line(0,1){0}}
\put(64.75,64){\line(0,1){0}}
\put(65.25,58){\line(0,1){0}}
\put(65,52.5){\line(0,1){0}}
\put(65,45.5){\line(0,1){0}}
\put(65.25,40.25){\line(0,1){0}}
\put(65.25,35.25){\line(0,1){0}}
\put(44.25,64.25){\line(0,1){0}}
\put(61.5,74){\line(0,1){0}}
\put(55,64.5){\line(0,1){0}}
\put(64.75,64){\line(0,1){0}}
\put(65.25,58){\line(0,1){0}}
\put(65,52.5){\line(0,1){0}}
\put(65,45.5){\line(0,1){0}}
\put(65.25,40.25){\line(0,1){0}}
\put(65.25,35.25){\line(0,1){0}}
\linethickness{.1 pt}
\multiput(40.5,61.25)(.04605263,.03289474){38}{\line(1,0){.04605263}}
\multiput(57.75,71)(.04605263,.03289474){38}{\line(1,0){.04605263}}
\multiput(51.25,61.5)(.04605263,.03289474){38}{\line(1,0){.04605263}}
\multiput(61,61)(.04605263,.03289474){38}{\line(1,0){.04605263}}
\multiput(61.5,55)(.04605263,.03289474){38}{\line(1,0){.04605263}}
\multiput(61.25,49.5)(.04605263,.03289474){38}{\line(1,0){.04605263}}
\multiput(61.25,42.5)(.04605263,.03289474){38}{\line(1,0){.04605263}}
\multiput(61.5,37.25)(.04605263,.03289474){38}{\line(1,0){.04605263}}
\multiput(59.5,32.25)(.04605263,.03289474){38}{\line(1,0){.04605263}}
\multiput(46,63.75)(.05182927,.03353659){82}{\line(1,0){.05182927}}
\multiput(66.5,63.5)(.05182927,.03353659){82}{\line(1,0){.05182927}}
\multiput(67,57.5)(.05182927,.03353659){82}{\line(1,0){.05182927}}
\multiput(42.75,60)(.06402439,.03353659){82}{\line(1,0){.06402439}}
\multiput(60,69.75)(.06402439,.03353659){82}{\line(1,0){.06402439}}
\multiput(63.25,59.75)(.06402439,.03353659){82}{\line(1,0){.06402439}}
\multiput(63.75,53.75)(.06402439,.03353659){82}{\line(1,0){.06402439}}
\multiput(63.5,48.25)(.06402439,.03353659){82}{\line(1,0){.06402439}}
\multiput(36.25,58.25)(.08653846,.03365385){52}{\line(1,0){.08653846}}
\multiput(53.5,68)(.08653846,.03365385){52}{\line(1,0){.08653846}}
\multiput(47,58.5)(.08653846,.03365385){52}{\line(1,0){.08653846}}
\multiput(56.75,58)(.08653846,.03365385){52}{\line(1,0){.08653846}}
\multiput(57.25,52)(.08653846,.03365385){52}{\line(1,0){.08653846}}
\multiput(57,46.5)(.08653846,.03365385){52}{\line(1,0){.08653846}}
\multiput(57,39.5)(.08653846,.03365385){52}{\line(1,0){.08653846}}
\multiput(57.25,34.25)(.08653846,.03365385){52}{\line(1,0){.08653846}}
\multiput(57.25,29.25)(.08653846,.03365385){52}{\line(1,0){.08653846}}
\multiput(64,67)(.05333333,.03333333){75}{\line(1,0){.05333333}}
\multiput(57.5,44)(.05833333,.03333333){60}{\line(1,0){.05833333}}
\multiput(63.75,45.75)(.06666667,.03333333){45}{\line(1,0){.06666667}}
\multiput(58,37.25)(.05416667,.03333333){60}{\line(1,0){.05416667}}
\multiput(57.5,25.5)(.03358209,.05223881){67}{\line(0,1){.05223881}}
\multiput(50,64)(.05487805,.03353659){82}{\line(1,0){.05487805}}
\linethickness{0.4pt}
\put(58,63.75){\makebox(0,0)[cc]{$N$}}
\put(25.75,35.75){\makebox(0,0)[cc]{$M_i$}}
\end{picture}

\caption{The case when $e_i=e_k\iv$.}
\label{f:15}
\end{center}
\end{figure}
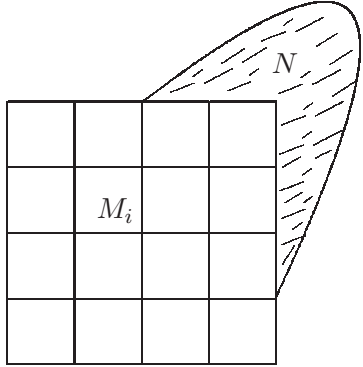
All exterior vertices of $N$, except for $(e_j)_{\pm}$, have degree at least $3$ in $N$ because the degrees of these vertices in $M_i$ are at most $3$ while there degrees in $M$ are $4$.
Therefore the $(4,4)$-curvature of every face and vertex in the flat map $N$, except for $(e_j)_{\pm}$, are non-positive, while the curvature of  each of $(e_j)_{\pm}$ is at most $1$. Hence the sum $I_f+I_v$ of all these curvatures for $N$ is at most $2<4$, which contradicts Lemma \ref{l:0}.

In case when $(e_j)_+=(e_k)_+$, we consider the submap $N$, without faces from $M_i$,  bounded by  $e_j$, $e_k^{-1}$ and by a subpath of $\partial M_i$.  It has at most $3$ vertices of positive
curvature (equal to $1$), namely, $(e_j)_-$, $(e_j)_+=(e_k)_+$ and $(e_k)_-$. This again contradicts
Lemma \ref{l:0} since $3<4$.

Now if we assume that edges $e_1,...,e_{4n+4}$  are enumerated clockwise, we see that for each $j$, $e_j, e_{j+1}$, $j=1,...,4n+4$ (where $j+1$ is $1$ if $j=4n+4$) must belong to the same face $\Pi_j$ which shares a vertex with $\partial(M_i)$ and every face sharing a vertex with $\partial(M_i)$ is one of the $\Pi_j$. Each pair of consecutive faces $\Pi_j, \Pi_{j+1}$ share exactly one edge $e_{j+1}$. Hence $M_{i+1}$ is isomorphic to the $(n+2)\times (n+2)$-square tessellated by unit squares.
The boundary of $M_{i+1}$ is simple since otherwise a part of this boundary bounds a flat submap $N$ with
at most one exterior vertex of positive curvature contrary to Lemma \ref{l:0} for $N$ since $1<4$.
\end{rk}

\subsection{Weakly exterior faces and the interior of a ($p,q$)-map.} \label{s:2}

In this section, $p$ and $q$ are postive integers satisfying
$\frac 1p +\frac 1q=\frac 12$.

A face (edge) in a map $M$ is called \emph{strongly interior} if it does not share a vertex with $\partial(M)$, otherwise it is called \emph{weakly exterior}.

\begin{lemma}\label{l:7} Let $M$ be a $(p,q)$-map, $o_1,...,o_m$ be its exterior vertices (counted counterclockwise). Then the number of weakly exterior faces does not exceed $\sum d(o_i)-2m$.
\end{lemma}

\proof Induction on the number of faces in $M$. If $M$ has a cut vertex, then  $M$ is a union of two submaps $M_1$ and $M_2$ intersecting by a vertex, and it is easy to see that the statement for $M$ follows from the statements for $M_1$ and $M_2$. Thus we can assume that the boundary path of $M$ is simple. Then every exterior vertex $o_j$ belongs to $d(o_j)-1$ weakly exterior faces, and two vertices $o_j, o_{j+1}$ (addition modulo $m$) belong to one face. So the sum $\sum (d(o_j)-1) = \sum d(o_j) -m$ overcounts weakly exterior faces by at least $m$. The statement of the lemma then follows.
\endproof

By the  \emph{interior} $M^0$ of $M$ we mean the union of all strongly interior faces of $M$, their vertices and edges. Note that $M^0$ may be empty. It may also be not connected in which case, it coincides with the
union of its maximal simple submaps $M^0_1, M^0_2,...$
It follows that every edge of $\partial M^0_i$ belongs
to a face of $M^0_i$ and to a weakly exterior face of $M$.

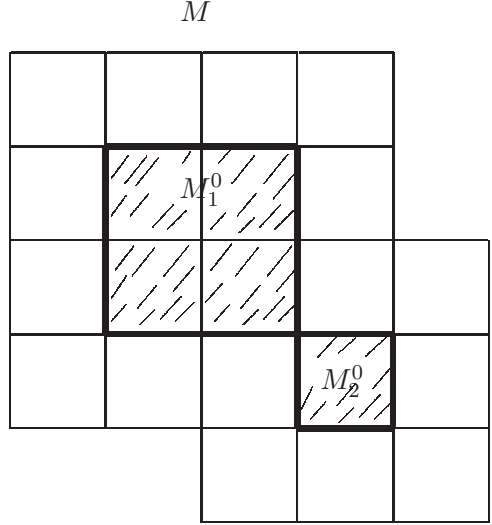
\begin{figure}[ht]
\begin{center}
\unitlength 1mm 
\linethickness{0.4pt}
\ifx\plotpoint\undefined\newsavebox{\plotpoint}\fi 
\begin{picture}(73.25,86)(0,0)
\put(9.5,43){\line(0,-1){12.5}}
\put(60.5,43){\line(0,-1){12.5}}
\put(60.5,30.5){\line(0,-1){12.5}}
\put(73.25,43){\line(0,-1){12.5}}
\put(73.25,30.5){\line(0,-1){12.5}}
\put(9.5,55.5){\line(0,-1){12.5}}
\put(60.5,55.5){\line(0,-1){12.5}}
\put(73.25,55.5){\line(0,-1){12.5}}
\put(9.5,68){\line(0,-1){12.5}}
\put(60.5,68){\line(0,-1){12.5}}
\put(9.5,80.5){\line(0,-1){12.5}}
\put(60.5,80.5){\line(0,-1){12.5}}
\put(22.25,43){\line(0,-1){12.5}}
\put(22.25,68){\line(0,-1){12.5}}
\put(22.25,80.5){\line(0,-1){12.5}}
\put(35,43){\line(0,-1){12.5}}
\put(35,30.5){\line(0,-1){12.5}}
\put(35,68){\line(0,-1){12.5}}
\put(35,80.5){\line(0,-1){12.5}}
\put(47.75,43){\line(0,-1){12.5}}
\put(47.75,30.5){\line(0,-1){12.5}}
\put(47.75,55.5){\line(0,-1){12.5}}
\put(47.75,80.5){\line(0,-1){12.5}}
\put(22,30.5){\line(-1,0){12.5}}
\put(22,43){\line(-1,0){12.5}}
\put(22,55.5){\line(-1,0){12.5}}
\put(22,68){\line(-1,0){12.5}}
\put(34.75,30.5){\line(-1,0){12.5}}
\put(34.75,43){\line(-1,0){12.5}}
\put(34.75,55.5){\line(-1,0){12.5}}
\put(49,55.5){\line(-1,0){12.5}}
\put(47.5,30.5){\line(-1,0){12.5}}
\put(47.5,18){\line(-1,0){12.5}}
\put(47.5,55.5){\line(-1,0){12.5}}
\put(47.5,68){\line(-1,0){12.5}}
\put(60.25,30.5){\line(-1,0){12.5}}
\put(60.25,18){\line(-1,0){12.5}}
\put(73,30.5){\line(-1,0){12.5}}
\put(73,18){\line(-1,0){12.5}}
\put(60.25,43){\line(-1,0){12.5}}
\put(73,43){\line(-1,0){12.5}}
\put(60.25,55.5){\line(-1,0){12.5}}
\put(73,55.5){\line(-1,0){12.5}}
\put(60.25,68){\line(-1,0){12.5}}
\put(9.75,80.5){\line(1,0){12.5}}
\put(22.5,80.5){\line(1,0){12.5}}
\put(35.25,80.5){\line(1,0){12.5}}
\put(48,80.5){\line(1,0){12.5}}
\linethickness{2 pt}
\put(22.25,67.75){\line(0,-1){24.75}}
\put(22.25,43){\line(1,0){25.5}}
\put(47.75,43){\line(0,1){25}}
\put(47.75,68){\line(-1,0){25.5}}
\put(47.75,43){\line(1,0){12.75}}
\put(60.5,43){\line(0,-1){12.5}}
\put(60.5,30.5){\line(-1,0){12.75}}
\put(47.75,30.5){\line(0,1){12.75}}
\linethickness{.4 pt}
\put(34.25,86){\makebox(0,0)[cc]{$M$}}
\put(53.75,36.75){\makebox(0,0)[cc]{$M_2^0$}}
\put(35,62.25){\makebox(0,0)[cc]{$M_1^0$}}
\put(35,55.5){\line(0,-1){12.5}}
\linethickness{.1 pt}
\multiput(27.75,66.75)(-.03370787,-.04213483){89}{\line(0,-1){.04213483}}
\multiput(42,66.75)(-.03370787,-.04213483){89}{\line(0,-1){.04213483}}
\multiput(28.75,54.5)(-.03370787,-.04213483){89}{\line(0,-1){.04213483}}
\multiput(42.5,54.5)(-.03370787,-.04213483){89}{\line(0,-1){.04213483}}
\multiput(28,61.75)(-.03333333,-.04){75}{\line(0,-1){.04}}
\multiput(42.25,61.75)(-.03333333,-.04){75}{\line(0,-1){.04}}
\multiput(29,49.5)(-.03333333,-.04){75}{\line(0,-1){.04}}
\multiput(42.75,49.5)(-.03333333,-.04){75}{\line(0,-1){.04}}
\multiput(55.5,36.5)(-.03333333,-.04){75}{\line(0,-1){.04}}
\multiput(31.25,60.25)(-.03353659,-.03658537){82}{\line(0,-1){.03658537}}
\multiput(45.5,60.25)(-.03353659,-.03658537){82}{\line(0,-1){.03658537}}
\multiput(32.25,48)(-.03353659,-.03658537){82}{\line(0,-1){.03658537}}
\multiput(46,48)(-.03353659,-.03658537){82}{\line(0,-1){.03658537}}
\multiput(58.75,35)(-.03353659,-.03658537){82}{\line(0,-1){.03658537}}
\multiput(33,59.5)(-.03333333,-.03333333){75}{\line(0,-1){.03333333}}
\multiput(47.25,59.5)(-.03333333,-.03333333){75}{\line(0,-1){.03333333}}
\multiput(34,47.25)(-.03333333,-.03333333){75}{\line(0,-1){.03333333}}
\multiput(47.75,47.25)(-.03333333,-.03333333){75}{\line(0,-1){.03333333}}
\multiput(60.5,34.25)(-.03333333,-.03333333){75}{\line(0,-1){.03333333}}
\multiput(46.5,67)(-.033505155,-.041237113){97}{\line(0,-1){.041237113}}
\multiput(33.25,54.75)(-.033505155,-.041237113){97}{\line(0,-1){.041237113}}
\multiput(47,54.75)(-.033505155,-.041237113){97}{\line(0,-1){.041237113}}
\multiput(47.25,64.25)(-.03353659,-.03963415){82}{\line(0,-1){.03963415}}
\multiput(34,52)(-.03353659,-.03963415){82}{\line(0,-1){.03963415}}
\multiput(47.75,52)(-.03353659,-.03963415){82}{\line(0,-1){.03963415}}
\multiput(60.5,39)(-.03353659,-.03963415){82}{\line(0,-1){.03963415}}
\multiput(38.25,59.25)(-.03333333,-.0375){60}{\line(0,-1){.0375}}
\multiput(25,47)(-.03333333,-.0375){60}{\line(0,-1){.0375}}
\multiput(38.75,47)(-.03333333,-.0375){60}{\line(0,-1){.0375}}
\multiput(51.5,34)(-.03333333,-.0375){60}{\line(0,-1){.0375}}
\multiput(42,58.5)(-.03333333,-.03333333){60}{\line(0,-1){.03333333}}
\multiput(28.75,46.25)(-.03333333,-.03333333){60}{\line(0,-1){.03333333}}
\multiput(42.5,46.25)(-.03333333,-.03333333){60}{\line(0,-1){.03333333}}
\multiput(55.25,33.25)(-.03333333,-.03333333){60}{\line(0,-1){.03333333}}
\multiput(53,42.75)(-.03333333,-.03666667){75}{\line(0,-1){.03666667}}
\multiput(49.75,36)(-.03333333,-.06666667){45}{\line(0,-1){.06666667}}
\multiput(39,55)(-.03353659,-.04268293){82}{\line(0,-1){.04268293}}
\put(36.25,51.5){\line(0,1){0}}
\multiput(37.75,49.75)(-.03358209,-.03731343){67}{\line(0,-1){.03731343}}
\multiput(26,54.75)(-.03333333,-.04333333){75}{\line(0,-1){.04333333}}
\multiput(25,50.75)(-.03333333,-.05){60}{\line(0,-1){.05}}
\multiput(25,61.5)(-.03333333,-.04583333){60}{\line(0,-1){.04583333}}
\multiput(25.25,66.75)(-.03358209,-.04477612){67}{\line(0,-1){.04477612}}
\multiput(29.25,66.5)(-.03370787,-.04213483){89}{\line(0,-1){.04213483}}
\multiput(59.75,42.75)(-.03658537,-.03353659){82}{\line(-1,0){.03658537}}
\multiput(55.5,42.5)(-.0375,-.03333333){60}{\line(-1,0){.0375}}
\multiput(39,67.75)(-.03365385,-.03846154){52}{\line(0,-1){.03846154}}
\multiput(33.5,67.25)(-.0333333,-.05){30}{\line(0,-1){.05}}
\end{picture}

\caption{The interior of a map and its components}
\label{f:1}
\end{center}
\end{figure}

Hence the intersection of
different submaps is either empty or consists of one vertex. Let us call these submaps
the {\it components} of $M^0$. Thus the boundary paths $\y_1^0, \y_2^0, \dots$ of the components $M_1^0, M^0_2,...$
are simple (see Fig. 1). Below we denote by $\y$ the union of these boundaries and set $|\y|=\sum_i |\y_i^0|$.

In the next lemma, we induct on the \emph{type} of a map $M$. By definition, the \emph{type} $\tau=\tau(\Pi)$
of a face $\Pi$ in $M$ is the number of interior (non-oriented) edges in $\partial\Pi$. If $m=\Area(M)$ and $(\tau_1,\tau_2,...,\tau_m)$ is the $m$-tuple of the types of all the faces of $M$ with
$\tau_1\ge\tau_2\ge...$, then $\tau(M)$ is the infinite string $(\tau_1,\dots,\tau_m,-1,-1,\dots)$.
We order types lexicographically: $\tau(M)\ge \tau(M')=(\tau'_1,\tau'_2, \dots)$ if $\tau_1>\tau'_1$ or
$\tau_1=\tau'_1$ but $\tau_2>\tau'_2$, and so on. For example, a map with only one face
has the type $(0,-1,-1,\dots)$. The set of types is obviously well ordered.

\begin{lemma} \label{l:1} Assume that a $(p,q)$  map $M$ has at least one face.
Also assume that the degree of every face of $M$ is at least $p$.
Let $\x$ be the boundary path of $M$ and $\y$ is defined as above. Then
\begin{equation}\label{e:1}
|\x|-|\y|\ge -I_f-2I^i_v+p
\end{equation}
where $I_f$ is the sum of the $(p,q)$-curvatures of all faces and $I^i_v$ is the sum of $(p,q)$-curvatures of all interior vertices of $M$.
\end{lemma}

\proof Let us denote $-I_f-2I_v^i$ by $J=J(M)$. Thus we need to prove that $|\x|-|\y|\ge J+p$.

{\bf Step 1.} The statement of the lemma is true if $M$ has only one face $\Pi$,
because we have $|\x|=d(\Pi) \ge p$, $|\y|=0$, $I_f= p-d(\Pi)$ and $I^i_v=0$.
Since the smallest type of a $(p,q)$ map that has faces is the type $(0,-1,-1,...)$ of a map consisting of one face (recall that we assume that $(p,q)$-maps do not have vertices of degree $1$), this gives the base of induction and we can assume that

\begin{quote}
($U_1$) The area of $M$ is greater than $1$.
\end{quote}

{\bf Step 2.} Suppose that $M$ can be cut into two maps $M_1$  and $M_2$
with smaller number of faces by a path of length at most $1$. Defining parameters $\x(j), \y(j), I_f(j), I^i_v(j)$ of the map $M_j$, $j=1,2$, in the natural way ($\x(j)$ is the boundary path of $M_j$, etc.), we have:
\begin{itemize}
\item $|\x(1)|+|\x(2)|\le |\x|+2$,
\item $|\y|=|\y(1)|+|\y(2)|$,
\item $I_f=I_f(1)+I_f(2)$ and $I^i_v=I^i_v(1)+I^i_v(2)$ since
no interior vertex of $M$ became exterior after the cutting.
\end{itemize}
Since $p\ge 2$, the statement
of the lemma follows from inequalities  $|\x(j)|-|\y(j)|\ge -I_f(j)-2I^i_v(j)+p$
(where $j=1,2$), which hold, since $\tau(M_j)<\tau(M)$ because $\Area(M_j)<\Area(M)$.
Thus, we may assume further that

\begin{quote}
($U_2$) $M$
has no cutting paths of length $\le 1$, hence, in particular, the boundary path $\x$ is simple.
\end{quote}

{\bf Step 3.} Assume there is an exterior vertex $o$ of degree $d>3$. Let $e_1,e_2,...,e_d$ be all edges ending in $o$, so that $e_1$ and $e_2$ are on the boundary path of a face $\Pi_1$, $e_2$ and $e_3$ are on the boundary path of a face $\Pi_2$, etc. Suppose $(e_1)_-=o'$ .  Then let us split the edge $e_2$ into two edges by a new vertex $o''$ in the middle of $e_2$, and replace $e_1$ with a new edge $e_1'$ going from $o'$ to $o''$. Note that this transformation does not change the type of $M$. Indeed, since $e_1$ is an exterior edge, the only faces that are changed by this transformation are $\Pi_1, \Pi_2$, but the number of interior edges on $\partial(\Pi_j), j=1,2,$ does not change (one of the two edges which are parts of the exterior edge $e_2$ is exterior in the new map and one is interior, and the new edge $e_1'$ is exterior as was $e_1$), hence $\tau(\Pi_j), j=1,2,$ does not change, and the type of the map does not change as well.  Also this transformation does not change the set of interior vertices of the map and their degrees. The degree of $\Pi_2$ increases by $1$ (because of the new vertex $o'$). Thus both $-I_f$ and $|\x|$ increases by $1$ and $I^i_v$ does not change. Hence the value of $|\x| -J$ does not change. The degree of the vertex $o$ decreases by 1 and the degree of the new vertex $o'$ is $3$. Therefore by doing this transformation, we will eventually get a map with the degrees of all exterior vertices at most $3$ with the same path $\y$ and the difference $|\x|-J$ as for $M$.

\begin{figure}[ht]
\begin{center}
\unitlength .7mm 
\linethickness{0.4pt}
\ifx\plotpoint\undefined\newsavebox{\plotpoint}\fi 
\begin{picture}(342.75,46.125)(0,0)
\thicklines
\put(18.25,39.5){\line(1,0){36.5}}
\multiput(54.75,39.5)(4.3125,-.03125){8}{\line(1,0){4.3125}}
\multiput(184.75,39.25)(4.3125,-.03125){8}{\line(1,0){4.3125}}
\put(52.75,39.25){\line(-6,-5){24}}
\put(185.75,39.5){\line(-6,-5){24}}
\qbezier(18.5,39.5)(11.25,39.5)(10,32.5)
\qbezier(151.5,39.75)(144.25,39.75)(143,32.75)
\qbezier(10,32.5)(9.625,28.625)(11.75,26.25)
\qbezier(143,32.75)(142.625,28.875)(144.75,26.5)
\qbezier(11.75,26.25)(13.5,24.75)(21.25,25.25)
\qbezier(144.75,26.5)(146.5,25)(154.25,25.5)
\qbezier(21.25,25.25)(27.25,25)(25.25,15.75)
\qbezier(154.25,25.5)(160.25,25.25)(158.25,16)
\qbezier(25.25,15.75)(26.25,15.375)(29.25,19.5)
\qbezier(158.25,16)(159.25,15.625)(162.25,19.75)
\qbezier(28.75,19)(30,10.75)(35.25,7.5)
\qbezier(161.75,19.25)(163,11)(168.25,7.75)
\qbezier(35.25,7.5)(45,14.625)(49.75,8.25)
\qbezier(168.25,7.75)(178,14.875)(182.75,8.5)
\qbezier(49.75,8.25)(55.625,1.5)(60,8.75)
\qbezier(182.75,8.5)(188.625,1.75)(193,9)
\multiput(60,8.75)(-.03372093,.141860465){215}{\line(0,1){.141860465}}
\multiput(193,9)(-.03372093,.141860465){215}{\line(0,1){.141860465}}
\put(29.5,31.75){\makebox(0,0)[cc]{$\Pi_1$}}
\put(162.5,32){\makebox(0,0)[cc]{$\Pi_1$}}
\put(45.5,21.75){\makebox(0,0)[cc]{$\Pi_2$}}
\put(178.5,22){\makebox(0,0)[cc]{$\Pi_2$}}
\put(35.25,43.25){\makebox(0,0)[cc]{$e_1$}}
\put(70.25,42.5){\makebox(0,0)[cc]{$e_d$}}
\put(111,24.5){\makebox(0,0)[cc]{$\Rightarrow$}}
\put(203.25,42.75){\makebox(0,0)[cc]{$e_d$}}
\put(37.75,31){\makebox(0,0)[cc]{$e_2$}}
\put(52.5,41.25){\makebox(0,0)[cc]{$o$}}
\put(185.25,41.5){\makebox(0,0)[cc]{$o$}}
\put(52.5,39.25){\circle*{2.5}}
\put(185.5,38.75){\circle*{2.5}}
\put(173.25,29.25){\circle*{2.5}}
\qbezier(151.5,39.75)(166.75,46.125)(173,30)
\put(174.6,33.25){\makebox(0,0)[cc]{$o'$}}
\put(165,43.75){\makebox(0,0)[cc]{$e_1'$}}
\linethickness{0.1pt}
\multiput(15.5,37.25)(-.03358209,-.04104478){67}{\line(0,-1){.04104478}}
\multiput(10.5,23.25)(-.03358209,-.04104478){67}{\line(0,-1){.04104478}}
\multiput(144.25,23.25)(-.03358209,-.04104478){67}{\line(0,-1){.04104478}}
\multiput(20.25,16.25)(-.03358209,-.04104478){67}{\line(0,-1){.04104478}}
\multiput(154,16.25)(-.03358209,-.04104478){67}{\line(0,-1){.04104478}}
\multiput(28.75,7)(-.03358209,-.04104478){67}{\line(0,-1){.04104478}}
\multiput(162.5,7)(-.03358209,-.04104478){67}{\line(0,-1){.04104478}}
\multiput(37,7)(-.03358209,-.04104478){67}{\line(0,-1){.04104478}}
\multiput(170.75,7)(-.03358209,-.04104478){67}{\line(0,-1){.04104478}}
\multiput(51.75,6)(-.03358209,-.04104478){67}{\line(0,-1){.04104478}}
\multiput(185.5,6)(-.03358209,-.04104478){67}{\line(0,-1){.04104478}}
\multiput(65,7.5)(-.03358209,-.04104478){67}{\line(0,-1){.04104478}}
\multiput(198.75,7.5)(-.03358209,-.04104478){67}{\line(0,-1){.04104478}}
\multiput(150,36.5)(-.03358209,-.04104478){67}{\line(0,-1){.04104478}}
\multiput(14,33)(-.03358209,-.04104478){67}{\line(0,-1){.04104478}}
\multiput(9,19)(-.03358209,-.04104478){67}{\line(0,-1){.04104478}}
\multiput(142.75,19)(-.03358209,-.04104478){67}{\line(0,-1){.04104478}}
\multiput(18.75,12)(-.03358209,-.04104478){67}{\line(0,-1){.04104478}}
\multiput(152.5,12)(-.03358209,-.04104478){67}{\line(0,-1){.04104478}}
\multiput(27.25,2.75)(-.03358209,-.04104478){67}{\line(0,-1){.04104478}}
\multiput(161,2.75)(-.03358209,-.04104478){67}{\line(0,-1){.04104478}}
\multiput(50.25,1.75)(-.03358209,-.04104478){67}{\line(0,-1){.04104478}}
\multiput(184,1.75)(-.03358209,-.04104478){67}{\line(0,-1){.04104478}}
\multiput(63.5,3.25)(-.03358209,-.04104478){67}{\line(0,-1){.04104478}}
\multiput(197.25,3.25)(-.03358209,-.04104478){67}{\line(0,-1){.04104478}}
\multiput(148.5,32.25)(-.03358209,-.04104478){67}{\line(0,-1){.04104478}}
\multiput(26.25,27.25)(-.03358209,-.04104478){67}{\line(0,-1){.04104478}}
\multiput(160.75,26.5)(-.03358209,-.04104478){67}{\line(0,-1){.04104478}}
\multiput(31.5,18.75)(-.03358209,-.04104478){67}{\line(0,-1){.04104478}}
\multiput(166,18)(-.03358209,-.04104478){67}{\line(0,-1){.04104478}}
\multiput(43.5,17.25)(-.03358209,-.04104478){67}{\line(0,-1){.04104478}}
\multiput(178,16.5)(-.03358209,-.04104478){67}{\line(0,-1){.04104478}}
\multiput(44,31)(-.03358209,-.04104478){67}{\line(0,-1){.04104478}}
\multiput(178.5,30.25)(-.03358209,-.04104478){67}{\line(0,-1){.04104478}}
\multiput(34.75,38)(-.03358209,-.04104478){67}{\line(0,-1){.04104478}}
\multiput(169.25,37.25)(-.03358209,-.04104478){67}{\line(0,-1){.04104478}}
\multiput(56.5,35.75)(-.03358209,-.04104478){67}{\line(0,-1){.04104478}}
\multiput(191,35)(-.03358209,-.04104478){67}{\line(0,-1){.04104478}}
\multiput(71,35.75)(-.03358209,-.04104478){67}{\line(0,-1){.04104478}}
\multiput(205.5,35)(-.03358209,-.04104478){67}{\line(0,-1){.04104478}}
\multiput(78.75,29.25)(-.03358209,-.04104478){67}{\line(0,-1){.04104478}}
\multiput(213.25,28.5)(-.03358209,-.04104478){67}{\line(0,-1){.04104478}}
\multiput(57.25,23.75)(-.03358209,-.04104478){67}{\line(0,-1){.04104478}}
\multiput(191.75,23)(-.03358209,-.04104478){67}{\line(0,-1){.04104478}}
\multiput(69.75,22.5)(-.03358209,-.04104478){67}{\line(0,-1){.04104478}}
\multiput(204.25,21.75)(-.03358209,-.04104478){67}{\line(0,-1){.04104478}}
\multiput(59.25,18.25)(-.03358209,-.04104478){67}{\line(0,-1){.04104478}}
\multiput(193.75,17.5)(-.03358209,-.04104478){67}{\line(0,-1){.04104478}}
\multiput(70.75,17.5)(-.03358209,-.04104478){67}{\line(0,-1){.04104478}}
\multiput(204.75,16.5)(-.03358209,-.04104478){67}{\line(0,-1){.04104478}}
\multiput(62.5,11.5)(-.03358209,-.04104478){67}{\line(0,-1){.04104478}}
\multiput(197,10.75)(-.03358209,-.04104478){67}{\line(0,-1){.04104478}}
\multiput(75.5,10.75)(-.03358209,-.04104478){67}{\line(0,-1){.04104478}}
\multiput(210,10)(-.03358209,-.04104478){67}{\line(0,-1){.04104478}}
\multiput(19.75,37.5)(-.03365385,-.05288462){52}{\line(0,-1){.05288462}}
\multiput(14.75,23.5)(-.03365385,-.05288462){52}{\line(0,-1){.05288462}}
\multiput(148.5,23.5)(-.03365385,-.05288462){52}{\line(0,-1){.05288462}}
\multiput(24.5,16.5)(-.03365385,-.05288462){52}{\line(0,-1){.05288462}}
\multiput(158.25,16.5)(-.03365385,-.05288462){52}{\line(0,-1){.05288462}}
\multiput(33,7.25)(-.03365385,-.05288462){52}{\line(0,-1){.05288462}}
\multiput(166.75,7.25)(-.03365385,-.05288462){52}{\line(0,-1){.05288462}}
\multiput(41.25,7.25)(-.03365385,-.05288462){52}{\line(0,-1){.05288462}}
\multiput(175,7.25)(-.03365385,-.05288462){52}{\line(0,-1){.05288462}}
\multiput(69.25,7.75)(-.03365385,-.05288462){52}{\line(0,-1){.05288462}}
\multiput(203,7.75)(-.03365385,-.05288462){52}{\line(0,-1){.05288462}}
\multiput(154.25,36.75)(-.03365385,-.05288462){52}{\line(0,-1){.05288462}}
\multiput(18.25,33.25)(-.03365385,-.05288462){52}{\line(0,-1){.05288462}}
\multiput(13.25,19.25)(-.03365385,-.05288462){52}{\line(0,-1){.05288462}}
\multiput(147,19.25)(-.03365385,-.05288462){52}{\line(0,-1){.05288462}}
\multiput(23,12.25)(-.03365385,-.05288462){52}{\line(0,-1){.05288462}}
\multiput(156.75,12.25)(-.03365385,-.05288462){52}{\line(0,-1){.05288462}}
\multiput(31.5,3)(-.03365385,-.05288462){52}{\line(0,-1){.05288462}}
\multiput(165.25,3)(-.03365385,-.05288462){52}{\line(0,-1){.05288462}}
\multiput(39.75,3)(-.03365385,-.05288462){52}{\line(0,-1){.05288462}}
\multiput(173.5,3)(-.03365385,-.05288462){52}{\line(0,-1){.05288462}}
\multiput(54.5,2)(-.03365385,-.05288462){52}{\line(0,-1){.05288462}}
\multiput(188.25,2)(-.03365385,-.05288462){52}{\line(0,-1){.05288462}}
\multiput(67.75,3.5)(-.03365385,-.05288462){52}{\line(0,-1){.05288462}}
\multiput(201.5,3.5)(-.03365385,-.05288462){52}{\line(0,-1){.05288462}}
\multiput(152.75,32.5)(-.03365385,-.05288462){52}{\line(0,-1){.05288462}}
\multiput(30.5,27.5)(-.03365385,-.05288462){52}{\line(0,-1){.05288462}}
\multiput(165,26.75)(-.03365385,-.05288462){52}{\line(0,-1){.05288462}}
\multiput(35.75,19)(-.03365385,-.05288462){52}{\line(0,-1){.05288462}}
\multiput(170.25,18.25)(-.03365385,-.05288462){52}{\line(0,-1){.05288462}}
\multiput(47.75,17.5)(-.03365385,-.05288462){52}{\line(0,-1){.05288462}}
\multiput(182.25,16.75)(-.03365385,-.05288462){52}{\line(0,-1){.05288462}}
\multiput(52,13)(-.03365385,-.05288462){52}{\line(0,-1){.05288462}}
\multiput(186.5,12.25)(-.03365385,-.05288462){52}{\line(0,-1){.05288462}}
\multiput(48.25,31.25)(-.03365385,-.05288462){52}{\line(0,-1){.05288462}}
\multiput(182.75,30.5)(-.03365385,-.05288462){52}{\line(0,-1){.05288462}}
\multiput(39,38.25)(-.03365385,-.05288462){52}{\line(0,-1){.05288462}}
\multiput(60.75,36)(-.03365385,-.05288462){52}{\line(0,-1){.05288462}}
\multiput(195.25,35.25)(-.03365385,-.05288462){52}{\line(0,-1){.05288462}}
\multiput(75.25,36)(-.03365385,-.05288462){52}{\line(0,-1){.05288462}}
\multiput(209.75,35.25)(-.03365385,-.05288462){52}{\line(0,-1){.05288462}}
\multiput(83,29.5)(-.03365385,-.05288462){52}{\line(0,-1){.05288462}}
\multiput(217.5,28.75)(-.03365385,-.05288462){52}{\line(0,-1){.05288462}}
\multiput(83.25,22.5)(-.03365385,-.05288462){52}{\line(0,-1){.05288462}}
\multiput(217.75,21.75)(-.03365385,-.05288462){52}{\line(0,-1){.05288462}}
\multiput(61.25,30)(-.03365385,-.05288462){52}{\line(0,-1){.05288462}}
\multiput(195.75,29.25)(-.03365385,-.05288462){52}{\line(0,-1){.05288462}}
\multiput(71.25,30.75)(-.03365385,-.05288462){52}{\line(0,-1){.05288462}}
\multiput(205.75,30)(-.03365385,-.05288462){52}{\line(0,-1){.05288462}}
\multiput(61.5,24)(-.03365385,-.05288462){52}{\line(0,-1){.05288462}}
\multiput(196,23.25)(-.03365385,-.05288462){52}{\line(0,-1){.05288462}}
\multiput(74,22.75)(-.03365385,-.05288462){52}{\line(0,-1){.05288462}}
\multiput(208.5,22)(-.03365385,-.05288462){52}{\line(0,-1){.05288462}}
\multiput(63.5,18.5)(-.03365385,-.05288462){52}{\line(0,-1){.05288462}}
\multiput(198,17.75)(-.03365385,-.05288462){52}{\line(0,-1){.05288462}}
\multiput(75,17.75)(-.03365385,-.05288462){52}{\line(0,-1){.05288462}}
\multiput(209,16.75)(-.03365385,-.05288462){52}{\line(0,-1){.05288462}}
\multiput(66.75,11.75)(-.03365385,-.05288462){52}{\line(0,-1){.05288462}}
\multiput(201.25,11)(-.03365385,-.05288462){52}{\line(0,-1){.05288462}}
\multiput(79.75,11)(-.03365385,-.05288462){52}{\line(0,-1){.05288462}}
\multiput(214.25,10.25)(-.03365385,-.05288462){52}{\line(0,-1){.05288462}}
\multiput(23.75,37.25)(-.03333333,-.06666667){45}{\line(0,-1){.06666667}}
\multiput(18.75,23.25)(-.03333333,-.06666667){45}{\line(0,-1){.06666667}}
\multiput(152.5,23.25)(-.03333333,-.06666667){45}{\line(0,-1){.06666667}}
\multiput(28.5,16.25)(-.03333333,-.06666667){45}{\line(0,-1){.06666667}}
\multiput(162.25,16.25)(-.03333333,-.06666667){45}{\line(0,-1){.06666667}}
\multiput(45.25,7)(-.03333333,-.06666667){45}{\line(0,-1){.06666667}}
\multiput(179,7)(-.03333333,-.06666667){45}{\line(0,-1){.06666667}}
\multiput(60,6)(-.03333333,-.06666667){45}{\line(0,-1){.06666667}}
\multiput(193.75,6)(-.03333333,-.06666667){45}{\line(0,-1){.06666667}}
\multiput(73.25,7.5)(-.03333333,-.06666667){45}{\line(0,-1){.06666667}}
\multiput(207,7.5)(-.03333333,-.06666667){45}{\line(0,-1){.06666667}}
\multiput(158.25,36.5)(-.03333333,-.06666667){45}{\line(0,-1){.06666667}}
\multiput(22.25,33)(-.03333333,-.06666667){45}{\line(0,-1){.06666667}}
\multiput(17.25,19)(-.03333333,-.06666667){45}{\line(0,-1){.06666667}}
\multiput(151,19)(-.03333333,-.06666667){45}{\line(0,-1){.06666667}}
\multiput(27,12)(-.03333333,-.06666667){45}{\line(0,-1){.06666667}}
\multiput(160.75,12)(-.03333333,-.06666667){45}{\line(0,-1){.06666667}}
\multiput(35.5,2.75)(-.03333333,-.06666667){45}{\line(0,-1){.06666667}}
\multiput(169.25,2.75)(-.03333333,-.06666667){45}{\line(0,-1){.06666667}}
\multiput(43.75,2.75)(-.03333333,-.06666667){45}{\line(0,-1){.06666667}}
\multiput(177.5,2.75)(-.03333333,-.06666667){45}{\line(0,-1){.06666667}}
\multiput(58.5,1.75)(-.03333333,-.06666667){45}{\line(0,-1){.06666667}}
\multiput(192.25,1.75)(-.03333333,-.06666667){45}{\line(0,-1){.06666667}}
\multiput(71.75,3.25)(-.03333333,-.06666667){45}{\line(0,-1){.06666667}}
\multiput(205.5,3.25)(-.03333333,-.06666667){45}{\line(0,-1){.06666667}}
\multiput(156.75,32.25)(-.03333333,-.06666667){45}{\line(0,-1){.06666667}}
\multiput(34.5,27.25)(-.03333333,-.06666667){45}{\line(0,-1){.06666667}}
\multiput(169,26.5)(-.03333333,-.06666667){45}{\line(0,-1){.06666667}}
\multiput(39.75,18.75)(-.03333333,-.06666667){45}{\line(0,-1){.06666667}}
\multiput(174.25,18)(-.03333333,-.06666667){45}{\line(0,-1){.06666667}}
\multiput(51.75,17.25)(-.03333333,-.06666667){45}{\line(0,-1){.06666667}}
\multiput(186.25,16.5)(-.03333333,-.06666667){45}{\line(0,-1){.06666667}}
\multiput(56,12.75)(-.03333333,-.06666667){45}{\line(0,-1){.06666667}}
\multiput(190.5,12)(-.03333333,-.06666667){45}{\line(0,-1){.06666667}}
\multiput(52.25,31)(-.03333333,-.06666667){45}{\line(0,-1){.06666667}}
\multiput(186.75,30.25)(-.03333333,-.06666667){45}{\line(0,-1){.06666667}}
\multiput(43,38)(-.03333333,-.06666667){45}{\line(0,-1){.06666667}}
\multiput(64.75,35.75)(-.03333333,-.06666667){45}{\line(0,-1){.06666667}}
\multiput(199.25,35)(-.03333333,-.06666667){45}{\line(0,-1){.06666667}}
\multiput(79.25,35.75)(-.03333333,-.06666667){45}{\line(0,-1){.06666667}}
\multiput(213.75,35)(-.03333333,-.06666667){45}{\line(0,-1){.06666667}}
\multiput(87,29.25)(-.03333333,-.06666667){45}{\line(0,-1){.06666667}}
\multiput(221.5,28.5)(-.03333333,-.06666667){45}{\line(0,-1){.06666667}}
\multiput(87.25,22.25)(-.03333333,-.06666667){45}{\line(0,-1){.06666667}}
\multiput(221.75,21.5)(-.03333333,-.06666667){45}{\line(0,-1){.06666667}}
\multiput(65.25,29.75)(-.03333333,-.06666667){45}{\line(0,-1){.06666667}}
\multiput(199.75,29)(-.03333333,-.06666667){45}{\line(0,-1){.06666667}}
\multiput(75.25,30.5)(-.03333333,-.06666667){45}{\line(0,-1){.06666667}}
\multiput(209.75,29.75)(-.03333333,-.06666667){45}{\line(0,-1){.06666667}}
\multiput(65.5,23.75)(-.03333333,-.06666667){45}{\line(0,-1){.06666667}}
\multiput(200,23)(-.03333333,-.06666667){45}{\line(0,-1){.06666667}}
\multiput(78,22.5)(-.03333333,-.06666667){45}{\line(0,-1){.06666667}}
\multiput(212.5,21.75)(-.03333333,-.06666667){45}{\line(0,-1){.06666667}}
\multiput(67.5,18.25)(-.03333333,-.06666667){45}{\line(0,-1){.06666667}}
\multiput(201.5,17.25)(-.03333333,-.06666667){45}{\line(0,-1){.06666667}}
\multiput(79,17.5)(-.03333333,-.06666667){45}{\line(0,-1){.06666667}}
\multiput(213,16.5)(-.03333333,-.06666667){45}{\line(0,-1){.06666667}}
\multiput(70.75,11.5)(-.03333333,-.06666667){45}{\line(0,-1){.06666667}}
\multiput(205.25,10.75)(-.03333333,-.06666667){45}{\line(0,-1){.06666667}}
\multiput(83.75,10.75)(-.03333333,-.06666667){45}{\line(0,-1){.06666667}}
\multiput(218.25,10)(-.03333333,-.06666667){45}{\line(0,-1){.06666667}}
\put(65,17.75){\makebox(0,0)[cc]{$.$}}
\put(199,17.25){\makebox(0,0)[cc]{$.$}}
\put(70.25,25){\makebox(0,0)[cc]{$.$}}
\put(204.25,24.5){\makebox(0,0)[cc]{$.$}}
\put(74.25,32.75){\makebox(0,0)[cc]{$.$}}
\put(208.25,32.25){\makebox(0,0)[cc]{$.$}}

\end{picture}

\caption{Step 3.}
\label{f:2}

\end{center}

\end{figure}

So we continue the proof under the additional assumption

\begin{quote}
($U_3$) The degree of every exterior vertex is at most $3$.
\end{quote}
\begin{rk}\label{r:1}
Note that ($U_3$) implies that every weakly exterior face of $M$ is exterior, i.e., every face which shares a vertex with $\partial(M)$ also shares an edge with $\partial(M)$.
\end{rk}

{\bf Step 4.} Assume there is a vertex $o$ of degree $2$ on an exterior face $\Pi$. Let us join two edges incident with $o$  into one edge and remove the vertex $o$. This does not change $\tau(M)$ because only exterior edges and vertices are affected.  Then the boundary $\y$ of the interior of the map does not change, the contribution of $\Pi$ to $J$ decreases by $1$, contributions of all other faces and vertices remain the same, and $|\x|$ also decreases by 1, hence $|\x|-|\y|-J$ will not change. Hence

\begin{quote} ($U_4$) We can remove vertices of degree $2$ of $\x$ (joining pairs of edges that share these vertices), provided the property $d(\Pi)\ge p$ is preserved, and we can split edges of $\x$ by new vertices of degree $2$  without changing $|\x|-|\y|-J$.
\end{quote}

{\bf Step 5.} Suppose an exterior face $\Pi$ of $M$ has boundary path of the form
$\uu\w$, where $\uu$ is a maximal subpath of the boundary path $\partial\Pi$ contained as a subpath in the
boundary path $\x$ of $M$. Note that $|\w|\ge 2$, because we excluded cutting paths of
length $\le 1$ by ($U_2$) and $M$ has more than one face by ($U_1$). Also note that  $|\uu|>0$ by Remark \ref{r:1}.

Suppose that $\w$ has an exterior vertex $o$ which is not equal to the end vertices of $\w$.
Then we add a vertex $o'$ of degree 2 on $\uu$ (using ($U_4$) ) and connect $o$ and $o'$ by a new edge $g$ cutting
up $\Pi$ into two faces of degrees $d_1$ and $d_2$, where $d_1+d_2-3 = d=d(\Pi)$.
For the new map $M'$ (with parameters $\x', \y'$, etc.) we have $|\x'|=|\x|+1$, $\y'=\y$. Instead of the face $\Pi$ of curvature $p-d$,
we have two faces with curvatures $p-d_1$ and $p-d_2$. Hence $I_f -I'_f = 3-p$.
Since the degrees of interior vertices were preserved, we have $J'-J=3-p$ and
$(|\x'|-J')-(|\x|-J)=p-2$. According to ($U_4$), the same difference $p-2$ has a map
with additional vertices of degree $2$ on $\uu$. So we may assume that $d_1, d_2\ge p$.

Cutting along $g$, we obtain new maps $M_1$ and $M_2$ with $\tau(M_j)<\tau(M)$ ($j=1,2$),
since $\Pi$ is subdivided into two faces with $\tau(\Pi_j)<\tau(\Pi), j=1,2$, where $\tau(\Pi_j)$ is computed in $M_j$. For the parameters $\x_j, \y_j, J_j$, etc. of the maps $M_j$, $j=1,2$, we have
$|\x_1|+|\x_2|=|\x|+3$, $|\y_1|+|\y_2|=|\y|$ and $J_1+J_2=J' = J+3-p$. So, by induction on  the type, we obtain $$|\x|-|\y|-J-p= (|\x_1|-|\y_1|-J_1-p)+(|\x_2|-|\y_2|-J_2-p) -3+3\ge 0+0-3+3=0,$$
as desired. Thus, we may assume that

\begin{quote} ($U_5$) For every exterior face $\Pi$ as above, the path $\w$ has no exterior vertices except its end vertices.
\end{quote}

\begin{figure}[ht]
\begin{center}
\unitlength .6 mm 
\linethickness{0.4pt}
\ifx\plotpoint\undefined\newsavebox{\plotpoint}\fi 
\begin{picture}(240.5,37.5)(0,0)
\thicklines
\put(5.5,33.5){\line(1,0){106}}
\put(134.5,33.5){\line(1,0){106}}
\put(5.5,5.75){\line(1,0){106}}
\put(134.5,5.75){\line(1,0){106}}
\multiput(27,33.75)(.03365385,-.28365385){52}{\line(0,-1){.28365385}}
\multiput(156,33.75)(.03365385,-.28365385){52}{\line(0,-1){.28365385}}
\multiput(28.75,19)(.0595854922,-.0336787565){386}{\line(1,0){.0595854922}}
\multiput(157.75,19)(.0595854922,-.0336787565){386}{\line(1,0){.0595854922}}
\multiput(51.75,6)(.0810397554,.0336391437){327}{\line(1,0){.0810397554}}
\multiput(180.75,6)(.0810397554,.0336391437){327}{\line(1,0){.0810397554}}
\multiput(78.25,17)(.033730159,.128968254){126}{\line(0,1){.128968254}}
\multiput(207.25,17)(.033730159,.128968254){126}{\line(0,1){.128968254}}
\multiput(23.25,20.25)(.2391304,-.0326087){23}{\line(1,0){.2391304}}
\multiput(152.25,20.25)(.2391304,-.0326087){23}{\line(1,0){.2391304}}
\multiput(28.75,19.5)(-.033505155,-.043814433){97}{\line(0,-1){.043814433}}
\multiput(157.75,19.5)(-.033505155,-.043814433){97}{\line(0,-1){.043814433}}
\multiput(78.25,17.25)(.12083333,.03333333){60}{\line(1,0){.12083333}}
\multiput(207.25,17.25)(.12083333,.03333333){60}{\line(1,0){.12083333}}
\multiput(78.25,17)(.040865385,-.033653846){104}{\line(1,0){.040865385}}
\multiput(207.25,17)(.040865385,-.033653846){104}{\line(1,0){.040865385}}
\put(123,20.25){\makebox(0,0)[cc]{$\Rightarrow$}}
\put(180.75,6){\line(0,1){27.5}}
\put(180.75,33.25){\circle*{1.581}}
\put(183,21.25){\makebox(0,0)[cc]{$g$}}
\put(168.5,23.25){\makebox(0,0)[cc]{$\Pi_1$}}
\put(194,23.25){\makebox(0,0)[cc]{$\Pi_2$}}
\put(53.25,22.5){\makebox(0,0)[cc]{$\Pi$}}
\put(52.25,3){\makebox(0,0)[cc]{$o$}}
\put(180.75,3){\makebox(0,0)[cc]{$o$}}
\put(180.75,36.75){\makebox(0,0)[cc]{$o'$}}
\put(58.25,37.5){\makebox(0,0)[cc]{$u$}}
\put(194,36.5){\makebox(0,0)[cc]{$u$}}
\put(65.75,15){\makebox(0,0)[cc]{$w$}}
\linethickness{.1 pt}
\multiput(13,31)(-.03358209,-.04850746){67}{\line(0,-1){.04850746}}
\multiput(142.25,30.5)(-.03358209,-.04850746){67}{\line(0,-1){.04850746}}
\multiput(11.25,25.75)(-.03358209,-.04850746){67}{\line(0,-1){.04850746}}
\multiput(140.5,25.25)(-.03358209,-.04850746){67}{\line(0,-1){.04850746}}
\multiput(8,18.75)(-.03358209,-.04850746){67}{\line(0,-1){.04850746}}
\multiput(137.25,18.25)(-.03358209,-.04850746){67}{\line(0,-1){.04850746}}
\multiput(10.75,15.75)(-.03358209,-.04850746){67}{\line(0,-1){.04850746}}
\multiput(140,15.25)(-.03358209,-.04850746){67}{\line(0,-1){.04850746}}
\multiput(21.25,12.75)(-.03358209,-.04850746){67}{\line(0,-1){.04850746}}
\multiput(150.5,12.25)(-.03358209,-.04850746){67}{\line(0,-1){.04850746}}
\multiput(28.75,12.5)(-.03358209,-.04850746){67}{\line(0,-1){.04850746}}
\multiput(158,12)(-.03358209,-.04850746){67}{\line(0,-1){.04850746}}
\multiput(32.5,32.5)(-.03358209,-.04850746){67}{\line(0,-1){.04850746}}
\multiput(161.75,32)(-.03358209,-.04850746){67}{\line(0,-1){.04850746}}
\multiput(31.5,26.75)(-.03358209,-.04850746){67}{\line(0,-1){.04850746}}
\multiput(160.75,26.25)(-.03358209,-.04850746){67}{\line(0,-1){.04850746}}
\multiput(37.5,21.25)(-.03358209,-.04850746){67}{\line(0,-1){.04850746}}
\multiput(166.75,20.75)(-.03358209,-.04850746){67}{\line(0,-1){.04850746}}
\multiput(42,31)(-.03358209,-.04850746){67}{\line(0,-1){.04850746}}
\multiput(171.25,30.5)(-.03358209,-.04850746){67}{\line(0,-1){.04850746}}
\multiput(57.5,31.25)(-.03358209,-.04850746){67}{\line(0,-1){.04850746}}
\multiput(186.75,30.75)(-.03358209,-.04850746){67}{\line(0,-1){.04850746}}
\multiput(65,28)(-.03358209,-.04850746){67}{\line(0,-1){.04850746}}
\multiput(194.25,27.5)(-.03358209,-.04850746){67}{\line(0,-1){.04850746}}
\multiput(54,22.25)(-.03358209,-.04850746){67}{\line(0,-1){.04850746}}
\multiput(183.25,21.75)(-.03358209,-.04850746){67}{\line(0,-1){.04850746}}
\multiput(50.75,14.75)(-.03358209,-.04850746){67}{\line(0,-1){.04850746}}
\multiput(180,14.25)(-.03358209,-.04850746){67}{\line(0,-1){.04850746}}
\multiput(86.5,32)(-.03358209,-.04850746){67}{\line(0,-1){.04850746}}
\multiput(215.75,31.5)(-.03358209,-.04850746){67}{\line(0,-1){.04850746}}
\multiput(86.75,26.5)(-.03358209,-.04850746){67}{\line(0,-1){.04850746}}
\multiput(216,26)(-.03358209,-.04850746){67}{\line(0,-1){.04850746}}
\multiput(86.75,19)(-.03358209,-.04850746){67}{\line(0,-1){.04850746}}
\multiput(216,18.5)(-.03358209,-.04850746){67}{\line(0,-1){.04850746}}
\multiput(84.5,12.5)(-.03358209,-.04850746){67}{\line(0,-1){.04850746}}
\multiput(213.75,12)(-.03358209,-.04850746){67}{\line(0,-1){.04850746}}
\multiput(70.75,11.5)(-.03358209,-.04850746){67}{\line(0,-1){.04850746}}
\multiput(200,11)(-.03358209,-.04850746){67}{\line(0,-1){.04850746}}
\multiput(103.25,32.5)(-.03358209,-.04850746){67}{\line(0,-1){.04850746}}
\multiput(232.5,32)(-.03358209,-.04850746){67}{\line(0,-1){.04850746}}
\multiput(98.25,27)(-.03358209,-.04850746){67}{\line(0,-1){.04850746}}
\multiput(227.5,26.5)(-.03358209,-.04850746){67}{\line(0,-1){.04850746}}
\multiput(96.75,19.5)(-.03358209,-.04850746){67}{\line(0,-1){.04850746}}
\multiput(226,19)(-.03358209,-.04850746){67}{\line(0,-1){.04850746}}
\multiput(98.25,12.75)(-.03358209,-.04850746){67}{\line(0,-1){.04850746}}
\multiput(227.5,12.25)(-.03358209,-.04850746){67}{\line(0,-1){.04850746}}
\multiput(16,30.25)(-.03358209,-.04477612){67}{\line(0,-1){.04477612}}
\multiput(145.25,29.75)(-.03358209,-.04477612){67}{\line(0,-1){.04477612}}
\multiput(14.25,25)(-.03358209,-.04477612){67}{\line(0,-1){.04477612}}
\multiput(143.5,24.5)(-.03358209,-.04477612){67}{\line(0,-1){.04477612}}
\multiput(11,18)(-.03358209,-.04477612){67}{\line(0,-1){.04477612}}
\multiput(140.25,17.5)(-.03358209,-.04477612){67}{\line(0,-1){.04477612}}
\multiput(13.75,15)(-.03358209,-.04477612){67}{\line(0,-1){.04477612}}
\multiput(143,14.5)(-.03358209,-.04477612){67}{\line(0,-1){.04477612}}
\multiput(24.25,12)(-.03358209,-.04477612){67}{\line(0,-1){.04477612}}
\multiput(153.5,11.5)(-.03358209,-.04477612){67}{\line(0,-1){.04477612}}
\multiput(31.75,11.75)(-.03358209,-.04477612){67}{\line(0,-1){.04477612}}
\multiput(161,11.25)(-.03358209,-.04477612){67}{\line(0,-1){.04477612}}
\multiput(35.5,31.75)(-.03358209,-.04477612){67}{\line(0,-1){.04477612}}
\multiput(164.75,31.25)(-.03358209,-.04477612){67}{\line(0,-1){.04477612}}
\multiput(34.5,26)(-.03358209,-.04477612){67}{\line(0,-1){.04477612}}
\multiput(163.75,25.5)(-.03358209,-.04477612){67}{\line(0,-1){.04477612}}
\multiput(40.5,20.5)(-.03358209,-.04477612){67}{\line(0,-1){.04477612}}
\multiput(169.75,20)(-.03358209,-.04477612){67}{\line(0,-1){.04477612}}
\multiput(45,30.25)(-.03358209,-.04477612){67}{\line(0,-1){.04477612}}
\multiput(174.25,29.75)(-.03358209,-.04477612){67}{\line(0,-1){.04477612}}
\multiput(60.5,30.5)(-.03358209,-.04477612){67}{\line(0,-1){.04477612}}
\multiput(189.75,30)(-.03358209,-.04477612){67}{\line(0,-1){.04477612}}
\multiput(68,27.25)(-.03358209,-.04477612){67}{\line(0,-1){.04477612}}
\multiput(197.25,26.75)(-.03358209,-.04477612){67}{\line(0,-1){.04477612}}
\multiput(57,21.5)(-.03358209,-.04477612){67}{\line(0,-1){.04477612}}
\multiput(186.25,21)(-.03358209,-.04477612){67}{\line(0,-1){.04477612}}
\multiput(53.75,14)(-.03358209,-.04477612){67}{\line(0,-1){.04477612}}
\multiput(183,13.5)(-.03358209,-.04477612){67}{\line(0,-1){.04477612}}
\multiput(89.5,31.25)(-.03358209,-.04477612){67}{\line(0,-1){.04477612}}
\multiput(218.75,30.75)(-.03358209,-.04477612){67}{\line(0,-1){.04477612}}
\multiput(89.75,25.75)(-.03358209,-.04477612){67}{\line(0,-1){.04477612}}
\multiput(219,25.25)(-.03358209,-.04477612){67}{\line(0,-1){.04477612}}
\multiput(89.75,18.25)(-.03358209,-.04477612){67}{\line(0,-1){.04477612}}
\multiput(219,17.75)(-.03358209,-.04477612){67}{\line(0,-1){.04477612}}
\multiput(87.5,11.75)(-.03358209,-.04477612){67}{\line(0,-1){.04477612}}
\multiput(216.75,11.25)(-.03358209,-.04477612){67}{\line(0,-1){.04477612}}
\multiput(73.75,10.75)(-.03358209,-.04477612){67}{\line(0,-1){.04477612}}
\multiput(203,10.25)(-.03358209,-.04477612){67}{\line(0,-1){.04477612}}
\multiput(106.25,31.75)(-.03358209,-.04477612){67}{\line(0,-1){.04477612}}
\multiput(235.5,31.25)(-.03358209,-.04477612){67}{\line(0,-1){.04477612}}
\multiput(101.25,26.25)(-.03358209,-.04477612){67}{\line(0,-1){.04477612}}
\multiput(230.5,25.75)(-.03358209,-.04477612){67}{\line(0,-1){.04477612}}
\multiput(99.75,18.75)(-.03358209,-.04477612){67}{\line(0,-1){.04477612}}
\multiput(229,18.25)(-.03358209,-.04477612){67}{\line(0,-1){.04477612}}
\multiput(101.25,12)(-.03358209,-.04477612){67}{\line(0,-1){.04477612}}
\multiput(230.5,11.5)(-.03358209,-.04477612){67}{\line(0,-1){.04477612}}
\multiput(20.25,30.25)(-.03365385,-.04807692){52}{\line(0,-1){.04807692}}
\multiput(149.5,29.75)(-.03365385,-.04807692){52}{\line(0,-1){.04807692}}
\multiput(18.5,25)(-.03365385,-.04807692){52}{\line(0,-1){.04807692}}
\multiput(147.75,24.5)(-.03365385,-.04807692){52}{\line(0,-1){.04807692}}
\multiput(15.25,18)(-.03365385,-.04807692){52}{\line(0,-1){.04807692}}
\multiput(144.5,17.5)(-.03365385,-.04807692){52}{\line(0,-1){.04807692}}
\multiput(18,15)(-.03365385,-.04807692){52}{\line(0,-1){.04807692}}
\multiput(147.25,14.5)(-.03365385,-.04807692){52}{\line(0,-1){.04807692}}
\multiput(28.5,12)(-.03365385,-.04807692){52}{\line(0,-1){.04807692}}
\multiput(157.75,11.5)(-.03365385,-.04807692){52}{\line(0,-1){.04807692}}
\multiput(36,11.75)(-.03365385,-.04807692){52}{\line(0,-1){.04807692}}
\multiput(165.25,11.25)(-.03365385,-.04807692){52}{\line(0,-1){.04807692}}
\multiput(39.75,31.75)(-.03365385,-.04807692){52}{\line(0,-1){.04807692}}
\multiput(169,31.25)(-.03365385,-.04807692){52}{\line(0,-1){.04807692}}
\multiput(38.75,26)(-.03365385,-.04807692){52}{\line(0,-1){.04807692}}
\multiput(168,25.5)(-.03365385,-.04807692){52}{\line(0,-1){.04807692}}
\multiput(44.75,20.5)(-.03365385,-.04807692){52}{\line(0,-1){.04807692}}
\multiput(174,20)(-.03365385,-.04807692){52}{\line(0,-1){.04807692}}
\multiput(49.25,30.25)(-.03365385,-.04807692){52}{\line(0,-1){.04807692}}
\multiput(178.5,29.75)(-.03365385,-.04807692){52}{\line(0,-1){.04807692}}
\multiput(64.75,30.5)(-.03365385,-.04807692){52}{\line(0,-1){.04807692}}
\multiput(194,30)(-.03365385,-.04807692){52}{\line(0,-1){.04807692}}
\multiput(72.25,27.25)(-.03365385,-.04807692){52}{\line(0,-1){.04807692}}
\multiput(201.5,26.75)(-.03365385,-.04807692){52}{\line(0,-1){.04807692}}
\multiput(61.25,21.5)(-.03365385,-.04807692){52}{\line(0,-1){.04807692}}
\multiput(190.5,21)(-.03365385,-.04807692){52}{\line(0,-1){.04807692}}
\multiput(58,14)(-.03365385,-.04807692){52}{\line(0,-1){.04807692}}
\multiput(187.25,13.5)(-.03365385,-.04807692){52}{\line(0,-1){.04807692}}
\multiput(93.75,31.25)(-.03365385,-.04807692){52}{\line(0,-1){.04807692}}
\multiput(223,30.75)(-.03365385,-.04807692){52}{\line(0,-1){.04807692}}
\multiput(94,25.75)(-.03365385,-.04807692){52}{\line(0,-1){.04807692}}
\multiput(223.25,25.25)(-.03365385,-.04807692){52}{\line(0,-1){.04807692}}
\multiput(94,18.25)(-.03365385,-.04807692){52}{\line(0,-1){.04807692}}
\multiput(223.25,17.75)(-.03365385,-.04807692){52}{\line(0,-1){.04807692}}
\multiput(91.75,11.75)(-.03365385,-.04807692){52}{\line(0,-1){.04807692}}
\multiput(221,11.25)(-.03365385,-.04807692){52}{\line(0,-1){.04807692}}
\multiput(78,10.75)(-.03365385,-.04807692){52}{\line(0,-1){.04807692}}
\multiput(207.25,10.25)(-.03365385,-.04807692){52}{\line(0,-1){.04807692}}
\multiput(110.5,31.75)(-.03365385,-.04807692){52}{\line(0,-1){.04807692}}
\multiput(239.75,31.25)(-.03365385,-.04807692){52}{\line(0,-1){.04807692}}
\multiput(105.5,26.25)(-.03365385,-.04807692){52}{\line(0,-1){.04807692}}
\multiput(234.75,25.75)(-.03365385,-.04807692){52}{\line(0,-1){.04807692}}
\multiput(104,18.75)(-.03365385,-.04807692){52}{\line(0,-1){.04807692}}
\multiput(233.25,18.25)(-.03365385,-.04807692){52}{\line(0,-1){.04807692}}
\multiput(105.5,12)(-.03365385,-.04807692){52}{\line(0,-1){.04807692}}
\multiput(234.75,11.5)(-.03365385,-.04807692){52}{\line(0,-1){.04807692}}
\end{picture}

\caption{Step 5.}
\label{f:3}

\end{center}

\end{figure}

{\bf Step 6.} Properties ($U_1$)-($U_5$) imply the following property

\begin{quote} ($U_6$) For every weakly exterior face $\Pi$ of $M$, we have $\partial(\Pi)=e\uu f\vv$ where $\uu=\uu(\Pi)$ is the subpath of the boundary path of $M$,  $|\uu(\Pi)|>0$, and $e,f$ are edges with exactly one exterior vertex while $\vv=\vv(\Pi)$ has no exterior vertices.
\end{quote}

Using the notation of Property ($U_6$),
suppose now that $|\vv(\Pi)|>1$, i.e., $\vv(\Pi)=\vv'\vv''$ with $|\vv'|,|\vv''| > 0 $. Let $o$ be the last vertex of $\vv'$. Then we add a new vertex $o'$ of degree 2 on an edge of $\uu$ (subdividing that edge into two edges) and add a new edge $t$ connecting $o$ and $o'$. As a result, the face $\Pi$ of degree $d=d(\Pi)$ is subdivided into two faces: a face $\Pi'$ of degree $d'$ and a face $\Pi''$ of degree $d''$, where $d'+d''=d+3$. Let $M'$ be the new map with parameters $\x', \y', I_f', J'$, etc. We have $\tau(M')<\tau(M)$ since
$\tau(\Pi'), \tau(\Pi'')<\tau(\Pi)$ by ($U_6$).

By ($U_4$), we can add new vertices of degree 2 to $\partial(\Pi'), \partial(\Pi'')$ so we can assume that $d'=d(\Pi')\ge p, d''=d(\Pi'')\ge p$. The contributions of $\Pi'$ and $\Pi''$ to $I'_f$ are $p-d'$ and $p-d''$, respectively,
while the contribution of $\Pi$ to $I_f$ was $p-d$. So $I_f-I'_f = p-(d'+d''-3) + (d'-p)+(d''-p) = 3-p$.
The contribution of the interior vertex $o$ to $2I^i_v$ is greater than its contribution to $2(I^i_v)'$ by $\frac{2p}{q}$ since this vertex is incident to the new edge $t$, i.e., $2(I^i_v)'-2I^i_v =2p/q$. Thus $J-J'=p-3-\frac{2p}q = -1$ because $\frac 1p +\frac 1q=\frac 12$. However we also have $|\x|-|\x'|= -1$ since one edge is subdivided by the vertex $o'$, and so $|\x'|-J'=\x-J$.
Thus we can assume,

\begin{quote} ($U$) $M$ satisfies ($U_6$) and for every exterior face $\Pi$, $|\vv(\Pi)|\le 1$.
\end{quote}

\begin{figure}[ht]

\begin{center}
\unitlength .8mm 
\linethickness{0.4pt}
\ifx\plotpoint\undefined\newsavebox{\plotpoint}\fi 
\begin{picture}(149,48.25)(0,0)
\thicklines
\put(29.5,44.25){\line(1,0){48}}
\put(101,44.75){\line(1,0){48}}
\multiput(29.5,44)(.033653846,-.118589744){156}{\line(0,-1){.118589744}}
\multiput(101,44.5)(.033653846,-.118589744){156}{\line(0,-1){.118589744}}
\multiput(34.75,25.5)(.0720973783,-.0337078652){267}{\line(1,0){.0720973783}}
\multiput(106.25,26)(.0720973783,-.0337078652){267}{\line(1,0){.0720973783}}
\multiput(54,16.5)(.079596413,.033632287){223}{\line(1,0){.079596413}}
\multiput(125.5,17)(.079596413,.033632287){223}{\line(1,0){.079596413}}
\multiput(71.75,24)(.033536585,.12347561){164}{\line(0,1){.12347561}}
\multiput(143.25,24.5)(.033536585,.12347561){164}{\line(0,1){.12347561}}
\put(53.75,47.5){\makebox(0,0)[cc]{$u$}}
\put(76.25,33.75){\makebox(0,0)[cc]{$f$}}
\put(147.75,34.25){\makebox(0,0)[cc]{$f$}}
\put(64.25,18){\makebox(0,0)[cc]{$v'$}}
\put(135.75,18.5){\makebox(0,0)[cc]{$v'$}}
\put(41.25,18.75){\makebox(0,0)[cc]{$v''$}}
\put(112.75,19.25){\makebox(0,0)[cc]{$v''$}}
\put(29.25,34){\makebox(0,0)[cc]{$e$}}
\put(100.75,34.5){\makebox(0,0)[cc]{$e$}}
\put(54.5,32.25){\makebox(0,0)[cc]{$\Pi$}}
\put(54,13.75){\makebox(0,0)[cc]{$o$}}
\put(85.75,32.75){\makebox(0,0)[cc]{$\Rightarrow$}}
\put(125.5,44.75){\circle*{1.581}}
\put(125.5,13.75){\makebox(0,0)[cc]{$o$}}
\put(126.25,48.25){\makebox(0,0)[cc]{$o'$}}
\multiput(125.75,16.75)(-.03125,3.46875){8}{\line(0,1){3.46875}}
\put(125.5,44.5){\line(0,1){.25}}
\put(127.75,30.5){\makebox(0,0)[cc]{$t$}}
\put(114.25,33.75){\makebox(0,0)[cc]{$\Pi''$}}
\put(136.25,33.5){\makebox(0,0)[cc]{$\Pi'$}}
\end{picture}

\caption{Step 6.}
\label{f:4}

\end{center}

\end{figure}

{\bf Step 7.} Now we consider the cases $p=3,4,6$ separately. Since $|\vv (\Pi)|\le 1$ by ($U$), we have $|\uu(\Pi)|\ge p-3$, and by $(U_4)$, one may assume that every exterior face has degree $p$ if $p>3$
and has degree 4 if $p=3$.

Since every exterior vertex has degree $2$ or $3$ (by  ($U_3$)), the difference $|\x|-|\y|$ is not smaller than the sum $S$ of $|\uu(\Pi)|-|\vv(\Pi)|$ for all exterior faces $\Pi$.

{\bf Case $p=4, q=4$.} In this case the degree of every exterior face is $4$. Then a vertex of degree $2$ can occur only on the path $\uu(\Pi)$ for some exterior face $\Pi$ with $|\vv(\Pi)|=0$. Let $N$ be the number of vertices of degree $2$ on $\partial(M)$.
 Then the contribution of the face containing that vertex to the sum $S$ is $2$ and $S\ge 2N$. The contribution of faces $\Pi$ with $|\uu(\Pi)|=|\vv(\Pi)|=1$ is $0$. The sum $I_v^b$ of curvatures of exterior vertices is $$\left(\frac 44(4-2)-1\right)N+\left(\frac 44(4-3)-1\right)(|\x|-N)=N.$$ By Lemma \ref{l:0}, $I_v^i+I_v^b+I_f=4$. Thus $N+I_v^i+I_c=4$, hence  $2N+2(I_v^i+I_f)=8$, and  $2N\ge J+8$ because $J=-2I_v^i-I_f$ and $I_f\le 0$ by the assumption of the lemma. Therefore $$|\x|-|\y|\ge J+8\ge J+p.$$

{\bf Case $p=6, q=3$.} Now every exterior face of $M$ has degree $6$.
Let $N_i$,  be the number of weakly exterior faces $\Pi$ with $|\vv(\Pi)|=i$, $i=0,1$. Then $|\uu(\Pi)|$ is
$4$ or $3$, respectively, and so $|\x| \ge 4N_0+3N_1$ and $|\y|\le N_1.$  So $|\x|-|\y| \ge 4N_0+2N_1$.

Every exterior face has either $3$ or $2$ vertices of degree $2$.
It is easy to compute now  that the sum $I_v^b$ of curvatures of the vertices of $\x$ is $$(3N_0+2N_1)\left(\frac 63(3-2)-1\right)+(|\x|-3N_0-2N_1) \left(\frac 63(3-3)-1\right)=-|\x|+6N_0+4N_1 = 2N_0+N_1.$$
Since by Lemma \ref{l:0}, $I_v^b+I_v^i+I_f=6$, we have $2N_0+N_1+I_v^i+I_f=6$. Hence $J+p<12+J
\le 2(-I_v^i-I_f+6) = 4N_0+2N_1\le |\x|-|\y|$, as required.

{\bf Case $p=3, q=6$.} Let $N_1$ be the
 number of exterior faces $\Pi$ of degree 4 with $|\uu(\Pi)|=2$, $|\vv(\Pi)|=0$, and $N_0$ be the number of exterior faces $\Pi$ of degree 4 with $|\uu(\Pi)|=|\vv(\Pi)|=1$.
Then we have $|\x|\ge 2N_1+N_0$, $|\y|\le N_0$,  and
$I^b_v= N_1+\frac 12 (|\x| - N_1)=\frac 12 (|\x|+N_1)$, because the curvature of
an exterior vertex of degree $2$ (respectively, $3$) is $\frac 12 (8-2)-1=1$ (respectively,  $\frac 12$).
Since $I_v^b+I_v^i+I_f=3$ by Lemma \ref{l:0},
we get $ |\x|= 2 I^b_v-N_1 =2(-I_v^i-I_f+3)-N_1$. Note that $-I_f\ge N_0+N_1$ since each quadrangle
contributes $-1$ to the sum $I_f$. Hence $$|\x|=2(-I_v^i-I_f+3)-N_1=J-I_f+6-N_1\ge J +N_0+N_1+6-N_1=  J+N_0+6$$
Therefore
$|\x|-|\y|\ge J+N_0+6-N_0  \ge J+6 > J+p$,
as desired. \endproof

\subsection{Adjustment}\label{ad}

Note that it is enough to prove Theorem \ref{main} for simple maps. Let $M$ be a simple $(p,q)$-map, $\frac 1p +\frac 1q=\frac 12$.
Every face of $M$ of degree $\ge 2p$ can be subdivided by diagonals into faces of degrees $p+1, ..., 2p-1$.
Every vertex $o$ of degree $d\ge 2q$  can be replaced by two (nearby) vertices $o', o''$ connected by an edge such that $d(o')+d(o'')=d(o)+2$ and $d(o')=q+1$. It can be done so that if $o$ is an exterior vertex, then exactly one of the vertices $o', o''$ is exterior in the resulting map. Each such transformation increases the total number of vertices and faces of negative $(p,q)$-curvature but it does not increase the number and the curvatures of the vertices and faces having non-negative curvature. Since the sum $I_f+I_v$ of $(p,q)$-curvatures of all vertices and faces cannot exceed $p$ by Lemma \ref{l:0}, any sequence of subdivisions of vertices and faces as above terminates. Clearly, if $M$ is a $(p,q)$-map, then the new map $M'$ is again a $(p,q)$-map, it has non-smaller area than $M$, the same perimeter,  and the set of (simple) flat submaps in $M'$ is a subset of the set of (simple) flat submaps of $M$ because the subdivisions do not introduce new flat vertices or faces.
The map $M'$ satisfies the following condition.

\begin{itemize}
\item[($B$)] The degree of every face (every vertex) of the simple $(p,q)$-map $M'$ is less than $2p$ (resp. $2q$).
\end{itemize}

If we have a $(p,q)$-map $M$ with condition ($B$), then one can construct a new map $M'$ satisfying condition ($B$), where {\it all} faces (not just interior ones) have degrees at least $p$. Namely, one subsequently cuts out  the exterior faces of degree less than $p$ (and also the edges containing the vertices of degree $1$ if such edges appear). The perimeter of $M'$ is at most $(p-1)n$, where $n$ is the perimeter of $M$,
$\Area(M')\ge \Area(M)-n$, and the maps $M$ and $M'$ have the same flat submaps.

The map $M'$ is a $(p,q)$-map satisfying the following additional condition

\begin{itemize}
\item[($D$)] The degree of every face of $M'$ is $\ge p$.
\end{itemize}

Let us call a $(p,q)$-map with the additional assumptions ($B$) and ($D$) a $\ctt$-map.
It is easy to see that Theorem \ref{main} follows from

\begin{theorem}\label{C4T4} The area
of an arbitrary  simple $\ctt$-map $M$ does not exceed $(\frac{3q}{2}+1)(r+p)n$,
provided $M$ contains no simple
flat submaps of radius greater than $r$.
\end{theorem}

Indeed, let $M'$ be the $\ctt$-map obtaining from a $(p,q)$-map $M$ satisfying the assumption of Theorem \ref{main} and condition ($B$) after removing some exterior faces. By Theorem \ref{C4T4}, we have $\Area(M')\le(\frac{3q}{2}+1)(r+p)(p-1)n$. Hence $\Area(M)\le(\frac{3q}{2}+1)(r+p)(p-1)n+n\le \frac{3}{2}(q+1)(r+p)(p-1)n.$

\subsection{Connecting non-flat vertices and faces with the boundary}


Note that all non-flat faces and all interior non-flat vertices of a \ctt-map $M$ have negative curvatures.

If $M$ contains non-flat faces or vertices, there exists a subgraph $\Gamma$ of the $1$-skeleton of $M$, such that every non-flat face or vertex
of $M$ can be connected with $\partial M$ by a path in $\Gamma$
We will assume that
$\Gamma$ is chosen with minimal number of edges $D(M)$. Then $\Gamma$ will be called  a \emph{connecting} subgraph of $M$.

\begin{rk}\label{r:2} The minimality of $\Gamma$ implies that every vertex of $\Gamma$ can be connected in $\Gamma$ with $\partial M$
by a unique reduced path. It follows that $\Gamma$ is a forest, where every maximal subtree has
exactly one vertex on the boundary $\partial M$.
\end{rk}

We shall use the notation from Section \ref{s:2}. Thus $\x$ is the boundary path of $M$, $\y$ is the union of the boundary paths of components of $M^0$, etc.

\begin{lemma} \label{l:8} We have  $D(M)\le (p-1)|\x|$.
\end{lemma}

\proof
We shall induct on the area $m$ of $M$. The statement is true
if $m=1$ since in this case $D(M)=0$.

Let $D_0=D(M^0)$ and let $\Gamma_0$ be the corresponding connecting subgraph of the interior $M^0$. By the induction hypothesis, $D_0\le (p-1)|\y|$.

Let $A_0$ (resp., $A$) be the number of non-flat faces and interior  vertices in
$M^0$ (in $M$).
By Remark \ref{r:2}, $\Gamma_0$ has at most $A_0$ vertices on $\partial M^0$, and so one needs
at most $A_0$ paths $\z_1, \z_2,\dots,$ to connect them with $\partial M$. Besides there are $A-A_0$ non-flat faces and interior vertices in $M$ which are not counted in $A_0$. One can connect them with $\partial M$ adding at most $A-A_0$ connecting paths $\y_1, \y_2,\dots$ to obtain
a graph $\Gamma'$ connecting all non-flat faces and interior vertices of $M$ with $\partial(M)$.
The lengths $|\z_i|$ and $|\y_j|$ cannot exceed
a half of the maximum of the perimeters of faces, and so $|\z_i|, |\y_j|\le p-1$ by Property ($B$).
Therefore
$D(M)\le D_0+A(p-1)$.

Since every non-flat face (resp. non-flat interior vertex) has curvature at most $-1$ (respectively, at most $-\frac 12$), we have
$A\le -I_f-2I_v^i=J \le |\x|-|\y|-p$ by Lemma \ref{l:1}.
Therefore
$$D(M)\le D_0+J(p-1) \le (p-1)|\y| +J(p-1)\le (p-1)(|\y| +(|\x|-|\y|-p)) < (p-1)|\x|$$

\endproof

\subsection{Cutting the map along its connecting subgraph and the proof of Theorem \ref{C4T4}}

As before $p,q$ are positive integers with $\frac 1p+\frac 1q=\frac 12$.

\begin{lemma}\label{l:6} Let $M$ be  a $\{p,q\}$-map of radius $0$ (i.e., all vertices of $M$ are exterior) and perimeter $n>0$. Then $\Area(M)\le \frac {q(n-2)}{2p}$.
\end{lemma}

\proof Induction on the number of faces in $M$. The statement is easy to check if $M$ contains only one face.  If $M$ contains more than one face, then $M$ has a cut vertex or cut edge.  In each of the two cases, the cut vertex or the cut edge separates $M$ into two submaps $M_1$ and $M_2$ with perimeters $n_1$, $n_2$ such that $n_1+n_2\le n+2$. Therefore $$\Area(M)=\Area(M_1)+\Area(M_2)\le \frac {q(n_1-2)}{2p}+\frac{q(n_2-2)}{2p}\le \frac {q(n-2)}{2p}.$$
\endproof

\begin{lemma}\label{l:5} Let $M$ be a $\{p,q\}$-map. Then the sum $I_v^b$ of curvatures of exterior vertices satisfies $I_v^b\ge p$.
\end{lemma}
\proof Indeed, we have $I_v^b+I_v^i+I_f=p$ by Lemma \ref{l:0}. Since $I_v^i\le 0$ and $I_f\le 0$, we have $I_v^b\ge p$.
\endproof

\begin{lemma}\label{l:4} Let $M$ be a $\{p,q\}$-map with perimeter $n>0$. Then the number $N$ of weakly exterior faces of $M$ is at most $\frac qp n-q$.
\end{lemma}
\proof Let $o_1,...,o_m$ ($m\le n$) be the exterior vertices of $M$.
Then $$I_v^b=\sum_j (\frac pq(q-d(o_j))-\mu(o_j))= -\frac pq\sum_j d(o_j)+pm-n$$ since $\sum \mu(o_j)=n$. Therefore  we have
$$\sum_j d(o_j)= \frac qp\left(-I_v^b+pm-n\right)\le -q+qm-\frac qp n$$ by Lemma \ref{l:5}. By Lemma \ref{l:7},
the number $N$ of weakly exterior faces in  $M$ is at most $\sum_j d(o_j)-2m$
Therefore $$N\le -q+(q-2)m-\frac qp n\le -q+\left(q-2-\frac qp\right)n=\frac qp n-q$$
since $\frac 1p+\frac 1q=\frac 12$.
\endproof

\begin{lemma}\label{l:3} If the radius of a $\{p,q\}$-map $M$ is at most $r-1$ and perimeter is $n$, then the area of $M$  does not exceed $\frac qp rn$.
\end{lemma}

\proof
If $r=1$ then this follows from Lemma \ref{l:6}, so let $r> 1$. Let $M^0$ be the interior of $M$.
Its boundary $y$ is the union of the boundaries of the components $M^0_1, M^0_2,\dots$ having radii $\le r-2$. So one may assume by induction on $r$ that $\Area (M^0)< \frac qp (r-1)|\y|$, which does not exceed $\frac qp (r-1) (n-p)$ by Lemma \ref{l:1}.  It follows
from Lemma \ref{l:4} that $\Area(M)< \frac qp (r-1)(n-p)+\frac qp n-q < \frac qp r n$.
\endproof

{\bf Proof of Theorem \ref{C4T4}.}
Let $\Gamma$ be a connecting graph of $M$. Let $e$ be an non-oriented edge of $\Gamma$ with one vertex
on $\partial M$. Then cutting $M$ along $e$, one obtains a $\ctt$-map $M_1$ with perimeter $|\x|+2$, where
all non-flat faces and vertices are connected with $\partial M_1$ by paths in the graph $\Gamma_1$, where $\Gamma_1$ is obtained from $\Gamma$ by removing the edge  $e$. We can continue cutting this way along the edges of $\Gamma$, until we obtain a $\ctt$-map
$\bar M$ of perimeter $|\x|+2D(M)\le (2p-1)|\x|$ (by Lemma \ref{l:8}), where every vertex and every face of curvature $>0$ is (weakly) exterior. Thus every component of the interior $\bar M^0$ of $\bar M$ is a simple flat map.

\begin{figure}[ht]
\begin{center}
\unitlength .7mm 
\linethickness{0.4pt}
\ifx\plotpoint\undefined\newsavebox{\plotpoint}\fi 
\begin{picture}(163,54.25)(0,0)
\thicklines
\qbezier(6,23.25)(7.375,40.375)(24.25,38)
\qbezier(24.25,38)(35,33.75)(47.75,41.5)
\qbezier(114.25,38)(125,33.75)(137.75,41.5)
\qbezier(47.75,41.5)(68.625,54.25)(73,35)
\qbezier(137.75,41.5)(158.625,54.25)(163,35)
\qbezier(73,35)(74.625,30.375)(72.75,24.25)
\qbezier(6,23.5)(5.625,12.625)(16.75,7.25)
\qbezier(16.75,7.25)(28.875,2)(40.5,4.75)
\qbezier(106.75,7.25)(118.875,2)(130.5,4.75)
\qbezier(40.5,4.75)(51.625,7.875)(59.25,6.5)
\qbezier(130.5,4.75)(141.625,7.875)(149.25,6.5)
\qbezier(59.25,6.5)(67.125,5.25)(72.5,24)
\qbezier(149.25,6.5)(157.125,5.25)(162.5,24)
\put(6,24.25){\line(1,0){10.5}}
\multiput(16.5,24.25)(.040064103,.033653846){156}{\line(1,0){.040064103}}
\multiput(16,24.5)(.03358209,-.041044776){134}{\line(0,-1){.041044776}}
\put(56.75,29.25){\line(1,0){17}}
\put(61,34.75){\line(0,-1){10.5}}
\qbezier(96.25,29.25)(101.25,41.625)(114.25,37.5)
\qbezier(95.5,20)(95.875,12.625)(106.75,7.75)
\multiput(96.5,29.5)(.14583333,-.03333333){60}{\line(1,0){.14583333}}
\multiput(115.75,28.5)(-.042553191,-.033687943){141}{\line(-1,0){.042553191}}
\multiput(109.75,23.75)(.033687943,-.035460993){141}{\line(0,-1){.035460993}}
\put(104.25,21.75){\line(-5,-1){8.75}}
\multiput(162.75,35.25)(-.096153846,-.033653846){104}{\line(-1,0){.096153846}}
\put(152.75,26.25){\line(5,-1){10}}
\put(12.75,27.5){\makebox(0,0)[cc]{$\Gamma$}}
\put(7.5,6.25){\makebox(0,0)[cc]{$M$}}
\put(84.5,26.5){\makebox(0,0)[cc]{$\Rightarrow$}}
\put(161.5,6.75){\makebox(0,0)[cc]{$\bar M$}}

\multiput(104.75,27.75)(.4565217,.0326087){23}{\line(1,0){.4565217}}
\multiput(103.75,22)(.108247423,-.033505155){97}{\line(1,0){.108247423}}
\put(153,32.25){\line(-1,3){2}}
\multiput(151,38.25)(-.03333333,-.14444444){45}{\line(0,-1){.14444444}}
\multiput(149.5,31.75)(-.09444444,-.03333333){45}{\line(-1,0){.09444444}}
\multiput(145.25,30.25)(.06097561,-.03353659){82}{\line(1,0){.06097561}}
\multiput(150.25,27.5)(.0333333,-.1833333){30}{\line(0,-1){.1833333}}

\multiput(151.25,22)(.03365385,.08653846){52}{\line(0,1){.08653846}}
\linethickness{.1pt}
\multiput(20.5,36.5)(-.036057692,-.033653846){104}{\line(-1,0){.036057692}}
\multiput(111.5,35)(-.036057692,-.033653846){104}{\line(-1,0){.036057692}}
\multiput(104.75,35.75)(-.036057692,-.033653846){104}{\line(-1,0){.036057692}}
\multiput(32.75,32.75)(-.036057692,-.033653846){104}{\line(-1,0){.036057692}}
\multiput(123.75,31.25)(-.036057692,-.033653846){104}{\line(-1,0){.036057692}}
\multiput(22.25,26)(-.036057692,-.033653846){104}{\line(-1,0){.036057692}}
\multiput(33.5,23)(-.036057692,-.033653846){104}{\line(-1,0){.036057692}}
\multiput(124.5,21.5)(-.036057692,-.033653846){104}{\line(-1,0){.036057692}}
\multiput(45.75,19.5)(-.036057692,-.033653846){104}{\line(-1,0){.036057692}}
\multiput(136.75,18)(-.036057692,-.033653846){104}{\line(-1,0){.036057692}}
\multiput(46.25,27)(-.036057692,-.033653846){104}{\line(-1,0){.036057692}}
\multiput(137.25,25.5)(-.036057692,-.033653846){104}{\line(-1,0){.036057692}}
\multiput(45.25,36.5)(-.036057692,-.033653846){104}{\line(-1,0){.036057692}}
\multiput(136.25,35)(-.036057692,-.033653846){104}{\line(-1,0){.036057692}}
\multiput(55.25,42)(-.036057692,-.033653846){104}{\line(-1,0){.036057692}}
\multiput(146.25,40.5)(-.036057692,-.033653846){104}{\line(-1,0){.036057692}}
\multiput(155.5,43.25)(-.036057692,-.033653846){104}{\line(-1,0){.036057692}}
\multiput(155.75,24.25)(-.036057692,-.033653846){104}{\line(-1,0){.036057692}}
\multiput(158.25,38.5)(-.036057692,-.033653846){104}{\line(-1,0){.036057692}}
\multiput(138,40.25)(-.036057692,-.033653846){104}{\line(-1,0){.036057692}}
\multiput(64.5,44.75)(-.036057692,-.033653846){104}{\line(-1,0){.036057692}}
\multiput(66.75,20)(-.036057692,-.033653846){104}{\line(-1,0){.036057692}}
\multiput(55.5,23.5)(-.036057692,-.033653846){104}{\line(-1,0){.036057692}}
\multiput(146.5,22)(-.036057692,-.033653846){104}{\line(-1,0){.036057692}}
\multiput(22,15.75)(-.036057692,-.033653846){104}{\line(-1,0){.036057692}}
\multiput(113,14.25)(-.036057692,-.033653846){104}{\line(-1,0){.036057692}}
\multiput(35.5,12.5)(-.036057692,-.033653846){104}{\line(-1,0){.036057692}}
\multiput(126.5,11)(-.036057692,-.033653846){104}{\line(-1,0){.036057692}}
\multiput(50.5,12.25)(-.036057692,-.033653846){104}{\line(-1,0){.036057692}}
\multiput(141.5,10.75)(-.036057692,-.033653846){104}{\line(-1,0){.036057692}}
\multiput(23.25,36)(-.036082474,-.033505155){97}{\line(-1,0){.036082474}}
\multiput(114.25,34.5)(-.036082474,-.033505155){97}{\line(-1,0){.036082474}}
\multiput(107.5,35.25)(-.036082474,-.033505155){97}{\line(-1,0){.036082474}}
\multiput(35.5,32.25)(-.036082474,-.033505155){97}{\line(-1,0){.036082474}}
\multiput(126.5,30.75)(-.036082474,-.033505155){97}{\line(-1,0){.036082474}}
\multiput(25,25.5)(-.036082474,-.033505155){97}{\line(-1,0){.036082474}}
\multiput(116,24)(-.036082474,-.033505155){97}{\line(-1,0){.036082474}}
\multiput(36.25,22.5)(-.036082474,-.033505155){97}{\line(-1,0){.036082474}}
\multiput(127.25,21)(-.036082474,-.033505155){97}{\line(-1,0){.036082474}}
\multiput(48.5,19)(-.036082474,-.033505155){97}{\line(-1,0){.036082474}}
\multiput(139.5,17.5)(-.036082474,-.033505155){97}{\line(-1,0){.036082474}}
\multiput(49,26.5)(-.036082474,-.033505155){97}{\line(-1,0){.036082474}}
\multiput(140,25)(-.036082474,-.033505155){97}{\line(-1,0){.036082474}}
\multiput(48,36)(-.036082474,-.033505155){97}{\line(-1,0){.036082474}}
\multiput(139,34.5)(-.036082474,-.033505155){97}{\line(-1,0){.036082474}}
\multiput(58,41.5)(-.036082474,-.033505155){97}{\line(-1,0){.036082474}}
\multiput(149,40)(-.036082474,-.033505155){97}{\line(-1,0){.036082474}}
\multiput(158.25,42.75)(-.036082474,-.033505155){97}{\line(-1,0){.036082474}}
\multiput(158.5,23.75)(-.036082474,-.033505155){97}{\line(-1,0){.036082474}}
\multiput(161,38)(-.036082474,-.033505155){97}{\line(-1,0){.036082474}}
\multiput(140.75,39.75)(-.036082474,-.033505155){97}{\line(-1,0){.036082474}}
\multiput(67.25,44.25)(-.036082474,-.033505155){97}{\line(-1,0){.036082474}}
\multiput(69.5,19.5)(-.036082474,-.033505155){97}{\line(-1,0){.036082474}}
\multiput(58.25,23)(-.036082474,-.033505155){97}{\line(-1,0){.036082474}}
\multiput(149.25,21.5)(-.036082474,-.033505155){97}{\line(-1,0){.036082474}}
\multiput(24.75,15.25)(-.036082474,-.033505155){97}{\line(-1,0){.036082474}}
\multiput(115.75,13.75)(-.036082474,-.033505155){97}{\line(-1,0){.036082474}}
\multiput(38.25,12)(-.036082474,-.033505155){97}{\line(-1,0){.036082474}}
\multiput(129.25,10.5)(-.036082474,-.033505155){97}{\line(-1,0){.036082474}}
\multiput(12,19.25)(-.036082474,-.033505155){97}{\line(-1,0){.036082474}}
\multiput(103,17.75)(-.036082474,-.033505155){97}{\line(-1,0){.036082474}}
\multiput(53.25,11.75)(-.036082474,-.033505155){97}{\line(-1,0){.036082474}}
\multiput(144.25,10.25)(-.036082474,-.033505155){97}{\line(-1,0){.036082474}}
\multiput(27,34.75)(-.038461538,-.033653846){104}{\line(-1,0){.038461538}}
\multiput(118,33.25)(-.038461538,-.033653846){104}{\line(-1,0){.038461538}}
\multiput(39.25,31)(-.038461538,-.033653846){104}{\line(-1,0){.038461538}}
\multiput(130.25,29.5)(-.038461538,-.033653846){104}{\line(-1,0){.038461538}}
\multiput(28.75,24.25)(-.038461538,-.033653846){104}{\line(-1,0){.038461538}}
\multiput(119.75,22.75)(-.038461538,-.033653846){104}{\line(-1,0){.038461538}}
\multiput(40,21.25)(-.038461538,-.033653846){104}{\line(-1,0){.038461538}}
\multiput(131,19.75)(-.038461538,-.033653846){104}{\line(-1,0){.038461538}}
\multiput(52.25,17.75)(-.038461538,-.033653846){104}{\line(-1,0){.038461538}}
\multiput(143.25,16.25)(-.038461538,-.033653846){104}{\line(-1,0){.038461538}}
\multiput(52.75,25.25)(-.038461538,-.033653846){104}{\line(-1,0){.038461538}}
\multiput(143.75,23.75)(-.038461538,-.033653846){104}{\line(-1,0){.038461538}}
\multiput(51.75,34.75)(-.038461538,-.033653846){104}{\line(-1,0){.038461538}}
\multiput(142.75,33.25)(-.038461538,-.033653846){104}{\line(-1,0){.038461538}}
\multiput(61.75,40.25)(-.038461538,-.033653846){104}{\line(-1,0){.038461538}}
\multiput(64.5,37)(-.038461538,-.033653846){104}{\line(-1,0){.038461538}}
\multiput(69.25,35.25)(-.038461538,-.033653846){104}{\line(-1,0){.038461538}}
\multiput(66.5,28.5)(-.038461538,-.033653846){104}{\line(-1,0){.038461538}}
\multiput(58,16.25)(-.038461538,-.033653846){104}{\line(-1,0){.038461538}}
\multiput(149,14.75)(-.038461538,-.033653846){104}{\line(-1,0){.038461538}}
\multiput(62.25,15)(-.038461538,-.033653846){104}{\line(-1,0){.038461538}}
\multiput(153.25,13.5)(-.038461538,-.033653846){104}{\line(-1,0){.038461538}}
\multiput(62,21.75)(-.038461538,-.033653846){104}{\line(-1,0){.038461538}}
\multiput(153,20.25)(-.038461538,-.033653846){104}{\line(-1,0){.038461538}}
\multiput(28.5,14)(-.038461538,-.033653846){104}{\line(-1,0){.038461538}}
\multiput(119.5,12.5)(-.038461538,-.033653846){104}{\line(-1,0){.038461538}}
\multiput(42,10.75)(-.038461538,-.033653846){104}{\line(-1,0){.038461538}}
\multiput(133,9.25)(-.038461538,-.033653846){104}{\line(-1,0){.038461538}}
\multiput(15.75,18)(-.038461538,-.033653846){104}{\line(-1,0){.038461538}}
\multiput(106.75,16.5)(-.038461538,-.033653846){104}{\line(-1,0){.038461538}}
\multiput(57,10.5)(-.038461538,-.033653846){104}{\line(-1,0){.038461538}}
\multiput(29.5,34.25)(-.033505155,-.041237113){97}{\line(0,-1){.041237113}}
\multiput(120.5,32.75)(-.033505155,-.041237113){97}{\line(0,-1){.041237113}}
\multiput(41.75,30.5)(-.033505155,-.041237113){97}{\line(0,-1){.041237113}}
\multiput(132.75,29)(-.033505155,-.041237113){97}{\line(0,-1){.041237113}}
\multiput(31.25,23.75)(-.033505155,-.041237113){97}{\line(0,-1){.041237113}}
\multiput(122.25,22.25)(-.033505155,-.041237113){97}{\line(0,-1){.041237113}}
\multiput(42.5,20.75)(-.033505155,-.041237113){97}{\line(0,-1){.041237113}}
\multiput(133.5,19.25)(-.033505155,-.041237113){97}{\line(0,-1){.041237113}}
\multiput(54.75,17.25)(-.033505155,-.041237113){97}{\line(0,-1){.041237113}}
\multiput(145.75,15.75)(-.033505155,-.041237113){97}{\line(0,-1){.041237113}}
\multiput(55.25,24.75)(-.033505155,-.041237113){97}{\line(0,-1){.041237113}}
\multiput(146.25,23.25)(-.033505155,-.041237113){97}{\line(0,-1){.041237113}}
\multiput(54.25,34.25)(-.033505155,-.041237113){97}{\line(0,-1){.041237113}}
\multiput(145.25,32.75)(-.033505155,-.041237113){97}{\line(0,-1){.041237113}}
\multiput(64.25,39.75)(-.033505155,-.041237113){97}{\line(0,-1){.041237113}}
\multiput(155.25,38.25)(-.033505155,-.041237113){97}{\line(0,-1){.041237113}}
\multiput(67,36.5)(-.033505155,-.041237113){97}{\line(0,-1){.041237113}}
\multiput(71.75,34.75)(-.033505155,-.041237113){97}{\line(0,-1){.041237113}}
\multiput(69,28)(-.033505155,-.041237113){97}{\line(0,-1){.041237113}}
\multiput(60.5,15.75)(-.033505155,-.041237113){97}{\line(0,-1){.041237113}}
\multiput(151.5,14.25)(-.033505155,-.041237113){97}{\line(0,-1){.041237113}}
\multiput(64.75,14.5)(-.033505155,-.041237113){97}{\line(0,-1){.041237113}}
\multiput(155.75,13)(-.033505155,-.041237113){97}{\line(0,-1){.041237113}}
\multiput(10.5,21.75)(-.033505155,-.041237113){97}{\line(0,-1){.041237113}}
\multiput(101.5,20.25)(-.033505155,-.041237113){97}{\line(0,-1){.041237113}}
\multiput(64.5,21.25)(-.033505155,-.041237113){97}{\line(0,-1){.041237113}}
\multiput(155.5,19.75)(-.033505155,-.041237113){97}{\line(0,-1){.041237113}}
\multiput(31,13.5)(-.033505155,-.041237113){97}{\line(0,-1){.041237113}}
\multiput(122,12)(-.033505155,-.041237113){97}{\line(0,-1){.041237113}}
\multiput(44.5,10.25)(-.033505155,-.041237113){97}{\line(0,-1){.041237113}}
\multiput(135.5,8.75)(-.033505155,-.041237113){97}{\line(0,-1){.041237113}}
\multiput(18.25,17.5)(-.033505155,-.041237113){97}{\line(0,-1){.041237113}}
\multiput(109.25,16)(-.033505155,-.041237113){97}{\line(0,-1){.041237113}}
\end{picture}

\caption{Cutting the map $M$ along the edges of $\Gamma$.}
\label{f:5}

\end{center}
\end{figure}
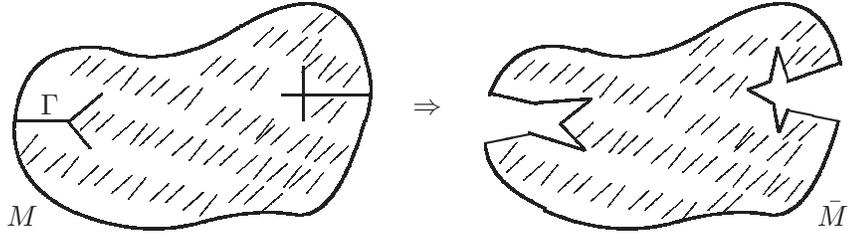

If $\bar M^0$ is empty, then the radius $\bar r$ of $\bar M$ is at most $r-1$ since by ($B$), the degree of every exterior face in $\bar M$ is less than $2p$. Hence by Lemma \ref{l:3} for $\bar M$, we have
$$\Area(M)=\Area(\bar M)\le \frac qp r(2p-1)n\le (\frac {3q}{2}+1)(r+p)n$$ since $\frac qp (2p-1) =\frac {3q}{2}+1$ by equality $\frac 1p+\frac 1q=\frac 12$, and the theorem is proved.

If $\bar M^0$ is not empty, then again by ($B$), it has a component $N$ of radius $\bar r^0\ge \bar r-p+1$
The map $N$ is a simple flat submap of $M$. Hence its radius $\bar r^0$ does not exceed $r$.
Therefore $\bar r -p+1\le r$ and $\bar r\le r+p-1$. By Lemma \ref{l:3},

$$\Area(M)=\Area(\bar M)\le \frac qp (2p-1)(r+p-1+1)n= (\frac{3q}{2} +1)(r+p)n$$
\endproof

\section{$(p,q)$-maps
that are
quasi-isometric to $\R^2$}\label{s:3}

Recall that a metric space $X$ with distance function $\dist_X$ is $(L,K)$-quasi-isometric to a metric space $Y$ with distance function $\dist_Y$, where $L>1$, $K>0$,  if there exists a mapping $\phi$ from $X$  to $Y$ such that $Y$ coincides with a tubular neighborhood of $\phi(X)$ and for every two vertices $o_1, o_2$ of $X$ we have
$$-K+\frac1L \dist_X(o_1,o_2)< \dist_{Y}(\phi(o_1),\phi(o_2)) < K+L\dist_X(o_1,o_2).$$
Two metric spaces $X$ and $Y$ are {\it quasi-isometric} if there is an $(L,K)$-quasi-isometry $X\to Y$ for some $L$ and $K$. This relation is reflexive, symmetric and transitive.

In this section we consider infinite planar maps. Here such a map $M$ is called {\it proper}
if the support of $M$ is the whole plane $\R^2$, and every disc on $\R^2$ intersects finitely
many faces, edges and vertices of $M$. The metric on $M$ is the combinatorial path metric on its
$1$-skeleton.

\begin{theorem}\label{t:qi} Let $M$ be a proper $(p,q)$-map, where positive integers $p$ and $q$ satisfy $\frac 1p + \frac 1q=\frac 12$.
Then  the $1$-skeleton of $M$ is quasi-isometric to  Euclidean plane if and only if $M$ has finite number of non-flat vertices and faces.
\end{theorem}

We will provide a proof for the case $p=q=4$, the other two cases are left for the reader (see Remark \ref{r:last}).

\subsection{The ``if" part of Theorem \ref{t:qi}}

A \emph{corridor} $\bbb$ of $M$ is a finite sequence of faces, where any two consecutive faces share a \emph{gluing} boundary edge, and two gluing edges of a face are not adjacent in the boundary path of the face.
In detail: a {\it corridor } is a sequence $$(e_0,\Pi_1, e_1,\Pi_2,\dots, e_{t-1},\Pi_t, e_t),$$ where $e_{i-1}$ and $e_i^{-1}$ are non-adjacent edges in the boundary of the face $\Pi_i$ for $i=1,\dots, t$.

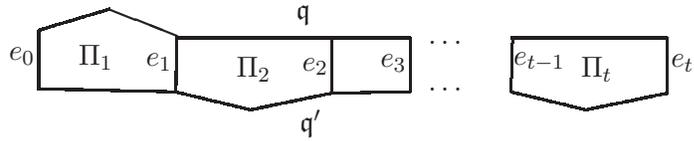
\begin{figure}[ht]
\begin{center}
\unitlength .8mm 
\linethickness{0.4pt}
\ifx\plotpoint\undefined\newsavebox{\plotpoint}\fi 
\begin{picture}(121.25,23.25)(0,0)
\thicklines
\put(14.25,9.5){\line(0,1){9.5}}
\multiput(14.25,19)(.112980769,.033653846){104}{\line(1,0){.112980769}}
\put(26,22.5){\line(5,-2){11.25}}
\multiput(37.25,18)(-.03125,-1.15625){8}{\line(0,-1){1.15625}}
\multiput(93,18)(-.03125,-1.15625){8}{\line(0,-1){1.15625}}
\multiput(37,8.75)(-1.5333333,.0333333){15}{\line(-1,0){1.5333333}}
\put(37,17.75){\line(1,0){26}}
\put(92.75,17.75){\line(1,0){26}}
\put(63,17.75){\line(0,-1){9.25}}
\put(118.75,17.75){\line(0,-1){9.25}}
\multiput(63,8.5)(-.16463415,-.03353659){82}{\line(-1,0){.16463415}}
\multiput(118.75,8.5)(-.16463415,-.03353659){82}{\line(-1,0){.16463415}}
\multiput(49.5,5.75)(-.14044944,.03370787){89}{\line(-1,0){.14044944}}
\multiput(105.25,5.75)(-.14044944,.03370787){89}{\line(-1,0){.14044944}}
\put(37,8.75){\line(0,1){0}}
\put(92.75,8.75){\line(0,1){0}}
\put(63,17.75){\line(1,0){13}}
\put(76,17.75){\line(0,-1){9}}
\multiput(76,8.75)(-1.625,-.03125){8}{\line(-1,0){1.625}}
\put(82.25,17){\makebox(0,0)[cc]{$\dots$}}
\put(82.25,9.25){\makebox(0,0)[cc]{$\dots$}}
\put(11.5,14.75){\makebox(0,0)[cc]{$e_0$}}
\put(34.25,13.75){\makebox(0,0)[cc]{$e_1$}}
\put(60.25,13){\makebox(0,0)[cc]{$e_2$}}
\put(121.25,13.5){\makebox(0,0)[cc]{$e_t$}}
\put(73.25,13.5){\makebox(0,0)[cc]{$e_3$}}
\put(97.50,13.75){\makebox(0,0)[cc]{$e_{t-1}$}}
\put(58.25,21.50){\makebox(0,0)[cc]{$\q$}}
\put(59.5,3.5){\makebox(0,0)[cc]{$\q'$}}
\put(23.75,14.5){\makebox(0,0)[cc]{$\Pi_1$}}
\put(50,12.5){\makebox(0,0)[cc]{$\Pi_2$}}
\put(107,12.5){\makebox(0,0)[cc]{$\Pi_t$}}
\end{picture}

\caption{A corridor}
\label{f:6}
\end{center}
\end{figure}

The boundary of $\bbb$ has the form $e_0\q e_t^{-1}(\q')^{-1}$, where the {\it sides} $\q$ and $\q'$ consist
of non-gluing edges of the faces $\Pi_1,\dots, \Pi_t$.

\begin{lemma}\label{no} In the above notation, no vertex of a ($4,4$)-map is passed by a side $\q$ or by $\q'$ twice.
\end{lemma}

\proof Assume that a corridor $\bbb$ is a counter-example.
Then we may assume that $q$ is a simple closed path bounding a submap $N$ of minimal possible area.
So $N$ contains no faces from $\bbb$. Every vertex of $\q$, except for the initial (= terminal) vertex $o$ has
degree at least $4$ in $M$, and so its degree in $N$ is at least $3$, as it follows from the definition of a corridor: two gluing edges of a face in a corridor are not adjacent in the boundary path of the face. Thus, the only vertex that can give a negative contribution to
the sum $I_f+I_v$ from Lemma \ref{l:0} for the map $N$ is $o$. But this contribution is at most
$1$, and so we have $I_f+I_v\le 1$, contrary to Lemma \ref{l:0}.\endproof

\begin{figure}[ht]
\begin{center}
\unitlength 1mm 
\linethickness{0.4pt}
\ifx\plotpoint\undefined\newsavebox{\plotpoint}\fi 
\begin{picture}(61.75,63.625)(0,0)
\thicklines
\multiput(12.75,8.5)(.069148936,.033687943){141}{\line(1,0){.069148936}}
\multiput(22.5,13.25)(.03333333,-.21111111){45}{\line(0,-1){.21111111}}
\multiput(13,8.75)(-.03333333,.10555556){45}{\line(0,1){.10555556}}
\multiput(11.25,13.75)(.069148936,.033687943){141}{\line(1,0){.069148936}}
\multiput(24,4.5)(.125,.03289474){38}{\line(1,0){.125}}
\multiput(28.75,5.75)(-.0333333,.3166667){30}{\line(0,1){.3166667}}
\qbezier(22.25,13.25)(23.375,38.5)(42,48.75)
\qbezier(42,48.75)(51.25,53.875)(51.5,44.5)
\qbezier(51.5,44.5)(51.5,34.375)(22.5,12.75)
\qbezier(20.5,18.75)(17.875,36)(40.75,54.25)
\qbezier(40.75,54.25)(52.25,63.625)(57.75,46.5)
\qbezier(57.75,46.5)(61.75,34.5)(27.75,14.5)
\multiput(16.25,15.75)(.03358209,-.06343284){67}{\line(0,-1){.06343284}}
\multiput(20.5,18.5)(.03365385,-.10576923){52}{\line(0,-1){.10576923}}
\multiput(20.5,26.25)(.08552632,-.03289474){38}{\line(1,0){.08552632}}
\multiput(24.5,36.75)(.07222222,-.03333333){45}{\line(1,0){.07222222}}
\multiput(33,47.25)(.04333333,-.03333333){75}{\line(1,0){.04333333}}
\multiput(53,55.75)(-.03353659,-.07926829){82}{\line(0,-1){.07926829}}
\multiput(58,44.25)(-.2826087,-.0326087){23}{\line(-1,0){.2826087}}
\put(51.5,43.5){\line(0,1){0}}
\multiput(55.75,36.5)(-.11538462,.03365385){52}{\line(-1,0){.11538462}}
\multiput(48.5,28.75)(-.041237113,.033505155){97}{\line(-1,0){.041237113}}
\multiput(42.25,23.75)(-.03353659,.03963415){82}{\line(0,1){.03963415}}
\multiput(36,19.25)(-.03666667,.03333333){75}{\line(-1,0){.03666667}}
\multiput(27.5,14.25)(-.13815789,-.03289474){38}{\line(-1,0){.13815789}}
\multiput(23.5,9.75)(.3166667,.0333333){15}{\line(1,0){.3166667}}
\put(40.75,54){\line(2,-5){2}}
\put(6.75,10){\makebox(0,0)[cc]{$\bbb$}}
\put(37.75,44){\makebox(0,0)[cc]{$\q$}}
\put(35.5,34.75){\makebox(0,0)[cc]{$N$}}
\thinlines
\multiput(24,25)(.0326087,-.1304348){23}{\line(0,-1){.1304348}}
\put(27.75,35.25){\line(1,0){4.25}}
\multiput(42.5,49.25)(.04,-.03333333){75}{\line(1,0){.04}}
\multiput(50.25,49.5)(-.0326087,-.1630435){23}{\line(0,-1){.1630435}}
\multiput(49.25,38.75)(-.03333333,.05416667){60}{\line(0,1){.05416667}}
\put(47.25,42){\line(0,1){0}}
\multiput(29.5,22.25)(.1630435,-.0326087){23}{\line(1,0){.1630435}}
\multiput(36.25,44.75)(-.03289474,-.11184211){38}{\line(0,-1){.11184211}}
\multiput(41,33.25)(.09210526,-.03289474){38}{\line(1,0){.09210526}}
\multiput(35.25,28.25)(.11184211,-.03289474){38}{\line(1,0){.11184211}}
\multiput(47.5,44.25)(.2666667,-.0333333){15}{\line(1,0){.2666667}}
\end{picture}

\caption{A corridor touching itself}
\label{f:7}
\end{center}

\end{figure}
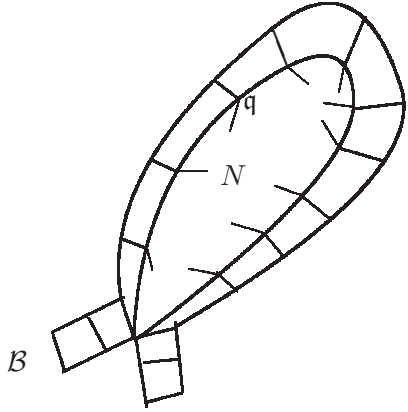

Lemma \ref{no} allows us to extend an arbitrary corridor infinitely in both directions: $$(\dots,e_{-1}, \Pi_0, e_0,\Pi_1, e_1,\Pi_2,\dots, e_{t-1},\Pi_t, e_t,\Pi_{t+1}\dots),$$ its sides are infinite simple paths subdividing the plane
in two parts (because the map $M$ is proper), these are \emph{infinite corridors}.
One can also consider \emph{semi-infinite} corridors
of the form $$(e_0,\Pi_1, e_1,\Pi_2,\dots, e_{t-1},\Pi_t, e_t,\dots).$$

\begin{lemma} \label{fv} Let $M$ be the map from Theorem \ref{t:qi} with a finite set of
non-flat faces and vertices. Then the $1$-skeleton of $M$ is quasi-isometric to the $1$-skeleton of a map with finitely many
non-flat vertices and without non-flat faces.
\end{lemma}
\proof Consider a infinite corridor $\bbb=(\dots,e_0,\Pi_1, e_1,\Pi_2,\dots, e_{t-1},\Pi_t, e_t,\dots)$ containing a non-flat face $\Pi_t$. Since the number of non-flat faces in $M$ is finite, we can assume that all the faces $\Pi_{t+1},\dots$ are flat.

\begin{figure}[ht]
\begin{center}
\unitlength .8mm 
\linethickness{0.4pt}
\ifx\plotpoint\undefined\newsavebox{\plotpoint}\fi 
\begin{picture}(171.25,24.75)(0,0)
\thicklines
\multiput(8.5,16.75)(.095652174,.033695652){230}{\line(1,0){.095652174}}
\multiput(102.5,17)(.095652174,.033695652){230}{\line(1,0){.095652174}}
\multiput(30.5,24.5)(.0608614232,-.0337078652){267}{\line(1,0){.0608614232}}
\multiput(124.5,24.75)(.0608614232,-.0337078652){267}{\line(1,0){.0608614232}}
\multiput(46.75,15.5)(.03125,-.90625){8}{\line(0,-1){.90625}}
\multiput(70.5,15.5)(-.03125,-.90625){8}{\line(0,-1){.90625}}
\put(47,8.25){\line(-1,0){38.25}}
\put(141,8.5){\line(-1,0){38.25}}
\put(8.75,8.25){\line(0,1){8.75}}
\put(102.75,8.5){\line(0,1){8.75}}
\put(46.5,15.25){\line(1,0){12.25}}
\put(140.5,15.5){\line(1,0){12.25}}
\put(70.75,15.25){\line(-1,0){12.25}}
\put(164.75,15.5){\line(-1,0){12.25}}
\put(47.25,8.25){\line(1,0){11.5}}
\put(141.25,8.5){\line(1,0){11.5}}
\put(70,8.25){\line(-1,0){11.5}}
\put(164,8.5){\line(-1,0){11.5}}
\put(58.75,15.25){\line(0,-1){6.75}}
\put(77.25,15.25){\makebox(0,0)[cc]{$\dots$}}
\put(171.25,15.25){\makebox(0,0)[cc]{$\dots$}}
\put(77.25,8.5){\makebox(0,0)[cc]{$\dots$}}
\put(3.25,8.5){\makebox(0,0)[cc]{$\dots$}}
\put(3.25,15.25){\makebox(0,0)[cc]{$\dots$}}
\put(98,8.5){\makebox(0,0)[cc]{$\dots$}}
\put(98,15.25){\makebox(0,0)[cc]{$\dots$}}
\put(171.25,8.5){\makebox(0,0)[cc]{$\dots$}}
\put(85.5,10.75){\makebox(0,0)[cc]{$\Rightarrow$}}
\multiput(141,8.75)(-.0352631579,.0336842105){475}{\line(-1,0){.0352631579}}
\multiput(152,8.25)(-.05326087,.033695652){230}{\line(-1,0){.05326087}}
\multiput(163.75,8.75)(-.052238806,.03358209){201}{\line(-1,0){.052238806}}
\multiput(58.75,13)(-.03125,-.03125){8}{\line(0,-1){.03125}}
\put(18.25,24){\makebox(0,0)[cc]{$\vv$}}
\put(39.5,23){\makebox(0,0)[cc]{$u$}}
\put(113,23.5){\makebox(0,0)[cc]{$\vv$}}
\put(135,22.75){\makebox(0,0)[cc]{$u$}}
\put(13.25,13){\makebox(0,0)[cc]{$e_{t-1}$}}
\put(48.25,10.5){\makebox(0,0)[cc]{$e_t$}}
\put(63.25,11.5){\makebox(0,0)[cc]{$e_{t+1}$}}
\put(107,13){\makebox(0,0)[cc]{$e_{t-1}$}}
\put(133.25,13.75){\makebox(0,0)[cc]{$f_t$}}
\put(143.75,12){\makebox(0,0)[cc]{$f_{t+1}$}}
\put(29.75,4.5){\makebox(0,0)[cc]{$\bbb$}}
\put(29.75,15){\makebox(0,0)[cc]{$\Pi_t$}}
\put(123.75,15){\makebox(0,0)[cc]{$\Pi_t$}}
\end{picture}

\caption{Modifying non-flat faces in a corridor}
\label{f:8}

\end{center}
\end{figure}
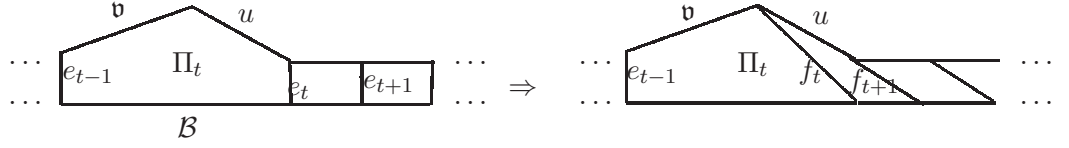

Let $\p, \p'$ be the two sides of $\bbb$ so that each $e_i$ connects a vertex $o_i$ on $\p$ with a vertex $o_i'$ on $\p'$.

Since $d=d(\Pi_t)\ge 5$, without loss of generality we can assume that the subpath $\w$ of $\p$ connecting $o_{t-1}$ and $o_t$ has length at least $2$, so it can be decomposed as $\w=\vv u$, where $u$ is one edge connecting $o''$ with $o_t$ and $|\vv|>0$. Then we modify the faces in $\bbb$ as follows:
replace the gluing edge $e_t$ by a new gluing edge $f_t$ connecting  $o''$ and $o'_t$, and replace every gluing edge
$e_s$ with $s>t$ by the new gluing edge  $f_s$ connecting $o_{s-1}$ with $o'_s$ (see Figure \ref{f:8}).
Then the degree of $\Pi$ decreases by $1$ and the degrees of all other faces are preserved. To complete the proof by induction, it suffices to notice that the $1$-skeleton of new map $M'$ is quasi-isometric to
$1$-skeleton of $M$, since the distances between the
vertices cannot increase/decrease more than two times when we passing from $M$ to $M'$.
\endproof

\begin{lemma} \label{off} Let $\bbb=(\dots,e_0,\Pi_1, e_1,\Pi_2,\dots, e_{t-1},\Pi_t, e_t,\dots)$ be a infinite corridor in $M$,
where every face $\Pi_i$ has degree $4$, and so $\Pi$ has the boundary of the form $e_{i-1}f_i(e_i)^{-1}g_i$,
where $f_i$ and $g_i$ are edges. Excising the faces of $B$ from $M$ and identifying the
edges $f_i$ and $g_i$ (for every $-\infty<i<\infty$), one obtains a new map $M'$. We claim that $M'$ is a
($4,4$)-map whose $1$-skeleton is quasi-isometric to the $1$-skeleton of $M$.
\end{lemma}
\proof Every end vertex of $e_i$ has degree $\ge 4$ in $M$. So the same must be true
in $M'$. Hence $M'$ is a ($4,4$)-map.

\begin{figure}[ht]

\begin{center}
\unitlength .7mm 
\linethickness{0.4pt}
\ifx\plotpoint\undefined\newsavebox{\plotpoint}\fi 
\begin{picture}(215.25,23)(0,0)
\put(13,23){\line(0,-1){22.75}}
\put(128.5,22.5){\line(0,-1){22.75}}
\put(33.75,23){\line(0,-1){22.75}}
\put(149.25,22.5){\line(0,-1){22.75}}
\put(54.5,23){\line(0,-1){22.75}}
\put(170,22.5){\line(0,-1){22.75}}
\put(75.25,23){\line(0,-1){22.75}}
\put(190.75,22.5){\line(0,-1){22.75}}
\put(96,23){\line(0,-1){22.75}}
\put(211.5,22.5){\line(0,-1){22.75}}
\multiput(5.5,13.5)(12.34375,-.03125){8}{\line(1,0){12.34375}}
\multiput(5.5,4.25)(12.34375,-.03125){8}{\line(1,0){12.34375}}
\put(42.5,16.5){\makebox(0,0)[cc]{$f_i$}}
\put(43,2){\makebox(0,0)[cc]{$g_i$}}
\put(109.25,9){\makebox(0,0)[cc]{$\Rightarrow$}}
\multiput(124,9.75)(11.40625,.03125){8}{\line(1,0){11.40625}}
\put(156,6.5){\makebox(0,0)[cc]{$f_i=g_i$}}
\put(37.25,9.25){\makebox(0,0)[cc]{$e_{i-1}$}}
\put(4.5,8.75){\makebox(0,0)[cc]{$\bbb$}}
\put(56.75,9.5){\makebox(0,0)[cc]{$e_i$}}
\end{picture}

\caption{Collapsing gluing edges of a corridor}
\label{f:10}
\end{center}

\end{figure}
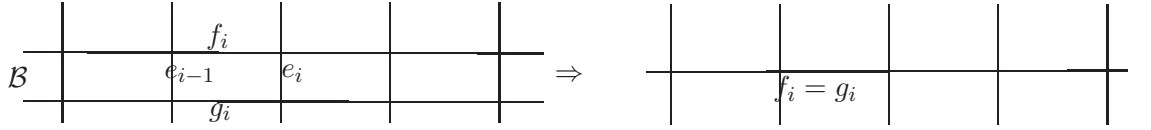

If two vertices  can be connected in $M'$ by a path $\p$ of length $m$, then
their preimages in $M$ can be connected in $M$ by a path $\q$ of length at most $2m+1$ since
$\q$ can be constructed from the edges of $\p$ and the gluing edges of $\bbb$. Conversely,
the distance between two vertices in $M'$ does not exceeds the distance between them in $M$.
The quasi-isometry of the $1$-skeletons follows.\endproof

\proof[Proof of the "if part" of Theorem \ref{t:qi}.] Suppose that $M$ has finitely many non-flat vertices and faces.
By lemma \ref{fv}, we may assume that $M$ has no non-flat faces,
and so every face has degree $4$. Let two distinct non-flat vertices $o$ and $o'$ be connected by a
shortest path $\p$. Then we chose an edge $e$ on $\p$ and consider a infinite corridor $\bbb$, where $e$ is one of the gluing edges of $\bbb$. Then the transformation $M\to M'$ from Lemma \ref{off} decreases the sum of
distances between non-flat vertices and replaces $M$ by a quasi-isometric map $M'$. So after a number
of such transformations, we shall have a ($4,4$)-map without non-flat faces and with at most one
non-flat vertex $o$. We will use the same notation $M$ for it.

Let us enumerate the edges $e_1,\dots, e_k$ with initial vertex $o$ in clockwise order; so $o$ lies
on the boundaries $e_if_ig_ie_{i+1}^{-1}$ (indexes are taken modulo $k$) of $k$ quadrangles
$\Pi_1,\dots,\Pi_k$.  Consider the semi-infinite corridors $\bbb_1 =(e_1,\Pi_1, g_1, \Pi'_1 \dots)$
and ${\mathcal C}_1 =(e_2,\Pi_1, f_1,\dots)$ starting with the face $\Pi_1$. They define semi-infinite
sides $\q_1$ and $\q'_1$ starting at $o$. Since $M$ is proper, the paths $\q_1, \q_2$ bound a submap ${\mathcal Q}_1$  of the plane.

There is a semi-infinite corridor ${\mathcal C}'_1$ starting
with the second edge of $\q'_1$ and the face $\Pi'_1$. This corridor has to share the whole side with
${\mathcal C}_1$ since it is made of quadrangles and all the vertices, except for $o$, have degree $4$.
Similarly, the semi-infinite corridor ${\mathcal C}''_1$ starting with the third edge of $\q'_1$ is glued up
to ${\mathcal C}'_1$ along the whole side, and so on.

Therefore $\mathcal Q_1$ with its path metric is isometric to a standard quadrant of the square grid $\Z^2$.
\begin{figure}[ht]

\begin{center}

\unitlength .8mm 
\linethickness{0.4pt}
\ifx\plotpoint\undefined\newsavebox{\plotpoint}\fi 
\begin{picture}(115.5,77.5)(0,0)
\multiput(39.5,77.5)(.03125,-8.875){8}{\line(0,-1){8.875}}
\multiput(7.75,24.5)(.06466262976,.03373702422){1156}{\line(1,0){.06466262976}}
\multiput(39.5,41.25)(.05296896086,-.03373819163){741}{\line(1,0){.05296896086}}
\multiput(39.5,41)(-.05407407407,.0337037037){675}{\line(-1,0){.05407407407}}
\multiput(39.5,51.25)(.065789474,.033625731){171}{\line(1,0){.065789474}}
\multiput(39.5,51.25)(-.065789474,.033625731){171}{\line(-1,0){.065789474}}
\multiput(39.5,30.75)(-.065789474,-.033625731){171}{\line(-1,0){.065789474}}
\multiput(39,30.75)(.065789474,-.033625731){171}{\line(1,0){.065789474}}
\put(50.75,57){\line(0,-1){10}}
\put(28.25,25){\line(0,1){10}}
\put(50.25,25){\line(0,1){10}}
\put(28.25,48){\line(-2,-1){11.5}}
\multiput(16.75,42.25)(.052325581,-.03372093){215}{\line(1,0){.052325581}}
\multiput(50.25,34.25)(.0658499234,.0336906585){653}{\line(1,0){.0658499234}}
\multiput(61.25,27.25)(.0658499234,.0336906585){653}{\line(1,0){.0658499234}}
\multiput(72.25,20.25)(.0658499234,.0336906585){653}{\line(1,0){.0658499234}}
\multiput(39.5,41.25)(.05295735901,-.03370013755){727}{\line(1,0){.05295735901}}
\multiput(50.25,47)(.05295735901,-.03370013755){727}{\line(1,0){.05295735901}}
\multiput(61,52.75)(.05295735901,-.03370013755){727}{\line(1,0){.05295735901}}
\multiput(71.75,58.5)(.05295735901,-.03370013755){727}{\line(1,0){.05295735901}}
\put(44.75,35.5){\makebox(0,0)[cc]{$e_1$}}
\put(46.75,43){\makebox(0,0)[cc]{$e_2$}}
\put(41,48.5){\makebox(0,0)[cc]{$e_3$}}
\put(51,40.5){\makebox(0,0)[cc]{$\Pi_1$}}
\put(59,25.5){\makebox(0,0)[cc]{$\q_1$}}
\put(65.25,58.5){\makebox(0,0)[cc]{$\q_2$}}
\put(115.5,50){\makebox(0,0)[cc]{$Q_1$}}
\put(84,58.5){\makebox(0,0)[cc]{$\bbb_1$}}
\put(83.5,19){\makebox(0,0)[cc]{$\ccc_1$}}
\put(93.5,26.25){\makebox(0,0)[cc]{$\ccc_1'$}}
\put(104,30.75){\makebox(0,0)[cc]{$\ccc_1''$}}
\multiput(28.5,56.75)(-.0333333,-.5833333){15}{\line(0,-1){.5833333}}
\end{picture}

\caption{Representing a map as a union of several quadrants of $\Z^2.$}
\label{f:11}
\end{center}
\end{figure}
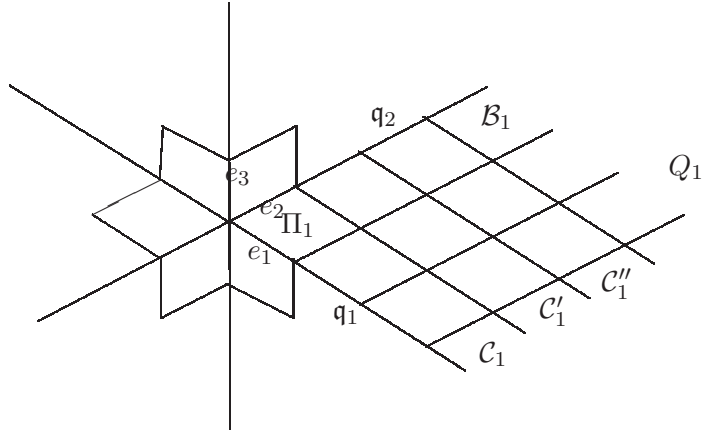

Similarly we have quadrants ${\mathcal Q}_2,\dots, {\mathcal Q}_k$, where each ${\mathcal Q}_i$ is bounded by
semi-infinite paths $\q_i$ and $\q'_i$, and as above, we have $\q'_i=\q_{i+1}$ (indices modulo $k$).
The $1$-skeleton of every submap ${\mathcal Q}_i$ is quasi-isometric to a quadrant on the Euclidean plane $\R^2$ with the Euclidean metric which, in turn, is quasi-isometric to a part $S_i$ of the plane bounded by two rays with common origin and angle $2\pi/k$ so that the union of all $S_i$ is the whole plane $\R^2$ (use polar coordinates). Combining all these quasi-isometries and using the fact that a quadrant of $\R^2$ is convex, we get a quasi-isometry between the $1$-skeleton of $M$ and $\R^2$.\endproof

\subsection{The ``only if" part of Theorem \ref{t:qi}}

Let $M$ be a proper $(4,4)$-map having infinitely many non-flat vertices and faces. By contradiction, suppose that $M$ with its path metric $\dist_M$ is $(L,K)$-quasi-isometric ($L>1$, $K>0$) to $\R^2$
If $\dist_M(o_1,o_2)>c$ for some $c\ge KL$, then $\dist_{\R^2}(\phi(o_1),\phi(o_2))>c_0=(c-KL)/L$
and so the growth of every $c$-separated set\footnote{A set of points $S$ in a metric space $X$ is called $c$-separated if $\dist_X(o_1,o_2)> c$ for every two distinct points $o_1,o_2\in S$.} $S$ of vertices of $M$ is at most quadratic, that is the function $\gamma_{S,o}(r)=|\{o'\in V' \mid \dist_M(o,o')\le r\}|$ is at most quadratic in $r$.

The number of vertices in a submap $M'$ of $M$ will be denoted by $\area(M')$ (recall that $\Area(M')$ denotes the number of faces in $M'$).

We start with the following well known

\begin{lemma}\label{qi:10}(See Theorem 6.2 in \cite[Chapter V]{LS}. Also it immediately follows from Lemmas \ref{l:4} and \ref{l:1} by induction on $n$.) If $N$ is a $(4,4)$-map with perimeter $n$, then $\area(N)\le kn^2$ for some constant $k> 0$.
\end{lemma}

Now we are going to modify faces of high degree.

\begin{lemma} \label{l:qi3} There exists a map $M'$ on the plane
\begin{enumerate}
\item with the same set $V$ of vertices as $M$,
\item with infinite set of non-flat vertices and faces,
\item for a marked vertex $o\in V$ and every $o'\in V$, $\dist_M(o,o')=\dist_{M'}(o,o')$,
\item for arbitrary vertices $o', o''$, we have $\dist_{M'}(o',o'')\le \dist_M(o',o'')$,
\item the degrees of all faces are at most $6$.
\end{enumerate}
\end{lemma}

\proof Let $\Pi$ be a face with $d=d(\Pi)\ge 7$, and vertices $o_1,...,o_d$ in the clockwise order.
Consider the difference  $f(i)=\dist(o_i,o)-\dist(o_{i+j},o)$ (indices modulo $d$). It is non-negative if $o_i$ is the farthest vertex from $o$ among $o_1,...,o_d$, and it is non-positive
if $o_i$ is the closest one. Since $|\dist(o_m,o)-\dist(o_{m+1},o)|\le 1$, we have $|f(m)-f(m+1)|\le 2$ for any $m$, and so there is $i$ such that $f(i)=|\dist(o_i,o)-\dist(o_{i+j},o)|\le 1.$

Let us connect vertices $o_i, o_{i+j}$ by a new, diagonal edge $e$ inside the cell $\Pi$, so that $\Pi$ is divided into two new cells $\Pi'$, $\Pi''$ both of degrees at least
$4$ and at least one of degree at least
$5$. This operation does not introduce any new vertices.
Let $M'$ be the new map on the plane. It is clear that Properties 1 and 4 of the lemma hold.

Let us show that distance from every vertex $o'$ to $o$ in $M'$ is the same as in $M$ (Property 3).
Let $g$ be a geodesic in $M'$ connecting $o'$ and $o$. If $g$ does not contain the new edge $e$, then the distance between $o'$ and $o$ did not change. So suppose that $g$ contains $e$. Without loss of generality we can assume that the vertices of $g$ in the natural order are $o',...,o_i, o_{i+j},...,o$. Since $g$ is a geodesic, $e$ appears in $g$ only once, and $\dist_{M'}(o_{i+j},o)=\dist_M(o_{i+j},o)$, and $\dist_{M'}(o_i,o)=\dist_{M}(o_{i+j},o)+1$. By the choice of the pair $(o_i,o_{i+j})$, $f(i)\in \{0,1,-1\}$. Since $\dist_M(o_i,o)\ge \dist_{M'}(o_i,o)$ we can deduce that $\dist_M(o_i,o)= \dist_{M'}(o_i,o)$, so there exists a geodesic $g'$ in $M'$ connecting $o'$ and $o$ and avoiding $e$. Hence $\dist_{M'}(o',o)=\dist_M(o',o)$.

This implies that we can cut all faces of degree $\ge 7$ by diagonals into several parts so that the resulting map on $\R^2$ satisfies all five properties of the lemma.
\endproof

Lemma \ref{l:qi3} implies that for the map $M'$ the growth function $\gamma_{S,o}$ of every $c$-separated set $S$ of vertices with respect to vertex $o$ is at most quadratic if $c$ is large enough (because every $c$-separated set $S$ of vertices in $M'$ is $c$-separated in $M$ by Property 4, and the functions $\gamma_{S,o}(r)$ for $M$ and $M'$ are the same by Property 3).
To obtain a contradiction with this quadratic growth, Lemma \ref{l:qi3} allows us to assume from  now on that the degrees of all faces in $M$ are at most $6$.

\begin{lemma} \label{qi:7} For every $r>0$ there exists a simple submap $N=N(r)$ of $M$ containing the vertex $o$, such that $\dist_M(o,\partial(N))\ge r$ and the maximal distance (in $M$) from $o$ to an exterior vertex of $N$ is at most $r+2$.
\end{lemma}

\proof Let $N$ be the smallest (with respect to the length of the boundary)  submap of $M$ containing all faces $\Pi$ of $M$
with $\dist_M(o,\Pi)\le r-1$. Such submap $N$ exists since $M$ is locally finite.  We claim that the boundary path of $N$ has no cut points. Indeed, let $o'$ be a cut point on $\partial(N)$ and $o'$ subdivides $N$ into two submaps $N_1, N_2$ containing faces, $N_1\cap N_2=\{o'\}$, $o\in N_1$. Suppose that $N_2$ contains a face $\Pi$ at distance (in $M$) at most $r-1$ from $o$.  Let $\g$ be a geodesic connecting $o$ and $\Pi'$ in $M$. Then every vertex on $\g$ is at distance (in $M$) at most $r-1$ from $o$. Hence every face of $M$ having a common vertex with  $\g$ is in $N$. Thus $\g$ is a path in the interior of $N$. Since $\g$ connects a vertex in $N_1$ with a vertex in $N_2$, $\g$ must contain $o'$.
Hence $o'$ is an interior vertex of $N$, a contradiction.

Since the boundary path $\partial N$ contains no cut points and has minimal length, it is simple.

There are no vertices $o'\in \partial(N)$ at distance (in $M$) at most $r-1$ from $o$. Indeed, otherwise every face of $M$ containing $o'$ would be at distance $\le r-1$ from $o$, and would be contained in $N$, hence $o$ would be an interior vertex in $N$, a contradiction.

Suppose that $N$ contains an exterior vertex $o'$ at distance (in $M$) at least $r+3$ from $o$. Then $o'$ belongs to an exterior face $\Pi$ of $N$. Since $d(\Pi)\le 6$, we have $\dist_M(o,\Pi)\ge r$. Therefore if we remove $\Pi$ from $M$ together with the longest subpath of $\partial(\Pi)$ containing $o'$ and contained in $\partial(N)$, we get a smaller submap $N'$ of $N$ containing all faces of $M$ at distance $\le r-1$ from $o$, a contradiction. Hence $\dist(o,o')\le r+2$ for every $o'\in \partial(M)$.
\endproof

Let $\Phi(r)$ be the number of vertices $o'$ of $M$ with $\dist_M(o',o)\le r$.

\begin{lemma}\label{Phi} The function $\Phi$ is super-quadratic, i.e $\lim_{r\to\infty} \Phi(r)/r^2 =\infty$.
\end{lemma}
\proof Let us denote by $\phi(r)$ the minimum of the numbers of vertices on the boundaries of the finite
submaps $Q$ with the property
that $Q$ is simple and $\dist_M(o,\partial Q)\ge r$. For any  $Q$ with this property, let $N$ be the component of the interior $Q^0$ containing $o$. Since every (exterior) face of $Q$ has degree at most $6$,
the boundary $\y_N$ of $N$ satisfies inequality $|\y_N|\ge \phi(r-3)$.
If $\x$ is the boundary path
of $Q$, then by Lemma \ref{l:1}, $|\x|\ge |\y_N|+J+4$, where $J= - I_f(Q)- 2I_v^i(Q)$.
Since $J$ is not less that the number $K=K(Q)$ of non-flat faces in $Q$ plus the number of interior in $Q$ non-flat vertices, we have $|\x| > |\y_N|+K\ge \phi(r-3)+K$, and so $\phi(r)>\phi(r-3)+K$.

If a non-flat face $\Pi$ (or vertex $o'$) lies in $M$ at a distance $\le r-1$ from $o$, then it
belongs to $Q$ (resp., it is interior in $Q$). Indeed, if $\Pi$ is not in $Q$ or $o'$ is not an interior vertex of $Q$, then any path connecting $\Pi$ or $o'$ with $o$ has to insersect $\partial Q$,
which contradicts the property $\dist_M(o,\partial Q)\ge r$. Hence $K\ge \phi(r-1)$,
where $\psi(r-1)$ is the number of non-flat vertices and faces of $M$ at the distance $\le r-1$ from $o$.
Therefore $\phi(r)>\phi(r-3)+\psi(r-1)$.

 Since $\psi(r)\to\infty$ as $r\to \infty$, we have $\phi(r)/r\ge \frac{1}{r} \sum_{0\le i\le r/6} \psi(r-1-3i)\ge \psi(\lfloor r/2\rfloor-1)/6\to \infty$.

Since the boundaries of the maps $N(r)$ and $N(r')$ from Lemma \ref{qi:7} do not intersect if $|r-r'|\ge 3$, we have $$\Phi(r)/r^2 \ge \frac{1}{r^2} \sum_{0\le i\le r/6} \phi(r-2-3i)
\ge \frac{1}{6r} \phi([r/2]-2)\to\infty $$ \endproof

\begin{lemma}\label{al} For an arbitrary integer $c\ge 1$, there is a super-linear function $\alpha(r)=\alpha_c(r)$ such that the boundary of every submap $N(r)$ satisfying the condition of Lemma \ref{qi:7} contains a $c$-separated set of vertices $S$ with $|S|\ge \alpha(r)$.
\end{lemma}
\proof We may assume that $r>c$. The simple boundary path $\x$ of $N(r)$ has winding number $\pm 1$ around the vertex $o$, and it is a sequence of edges, i.e., of subpaths of length $1\le c$. Therefore there is the smallest $t\ge 1$ and the vertices $o_1,\dots, o_t$ of $\x$ such that $\dist_M (o_i,o_{i+1})\le c$ (indices modulo $t$) and the geodesic paths $\z_i=o_i - o_{i+1}$ (in $M$) form
 a product $\z=\z_1\dots \z_t$ with non-zero winding number around $o$. (Self-intersections are allowed for $\z$, but no $\z_i$ goes through $o$ since $|\z_i|\le c<r \le \dist_M(o, o_i)$, and so the winding number is well defined.)

Suppose there is a pair of distinct vertices $o_i, o_j$, where $i<j$ and $i-j\ne \pm 1$ (mod $t$),
with $\dist(o_i,o_j)\le c$. Then a geodesic  path $\bar \z =o_i - o_j$ defines two new
closed paths $\z'=\bar \z \z_j\dots \z_t\z_1\dots \z_{i-1}$ and $\z''=\z_i^{-1}\dots \z_{j-1}^{-1}\bar \z$ with  numbers of factors less than $t$. Since one of the paths $\z'$ and $\z''$ has nonzero
winding number with respect to $o$, this contradicts the minimality in the choice of $t$. Hence
$\dist(o_i,o_j)>c$, and so one can chose the set $S$ of cardinality $\ge(t-1)/2$.

Assume now that $t< \sqrt[4]{r^2\Phi(r-C/2)}$, where $\Phi (r)$ is defined before Lemma \ref{Phi}. Note that $o$
belongs to a submap $L$ bounded by a simple closed path $\w$, which is a product of some
pieces of the paths $\z_i$-s. Therefore $|\w|\le |\z|\le tc$, and by Lemma \ref{qi:10},
$\area(L)\le k(ct)^2\le Dr\sqrt{\Phi(r-c/2)}$, where $D=kc^2$.

On the other hand, $\dist_M(o',o)\ge r-c/2$ for every vertex $o'$ of $\w$ since $\z_i\le c$ for every
$i$. Therefore $\area(L)\ge \Phi(r-c/2)$ by Lemma \ref{Phi}. We obtain inequality
$\Phi(r-c/2)\le Dr\sqrt{\Phi(r-c/2)}$, which can hold only for finitely many values of $r$
since the function $\Phi$ is super-quadratic by Lemma \ref{Phi}.
Thus $t=t(r)\ge \sqrt[4]{r^2\Phi(r-c/2)}$ for every $r\ge r_0$, and since $|S|\ge \frac12 (t-1)$, one can define
$\alpha(r)$ to be equal to $\frac12(\sqrt[4]{r^2\Phi(r-c/2)}-1)$ if $r\ge r_0$  and $\alpha(r)=1$
if $r< r_0$. \endproof

Now we can prove that for any given integer $c\ge 1$, the map $M$ contains an infinite $c$-separated set $S$ of vertices which grows super-quadratically with respect to the vertex $o$, which would give the desired contradiction. The boundary $\partial N(r)$ has a $c$-separated
subset $S_r$ with at least $\alpha(r)$ vertices. Since the distance between $\partial N(r)$
and $\partial N(r')$ is greater than $c$ for $r-r'\ge c+3$, the union $S(r)=S_r\cup S_{r-(c+3)}\cup S_{r-2(c+3)}\cup\dots$ is $c$-separated and
$$|S(r)|/r^2\ge \frac{1}{r^2}\sum_{0\le i\le \frac{r}{2(c+3)}}\alpha(r-i(c+3))
\le \frac{1}{2(c+3)r}\min_{0\le i\le \frac{r}{2(c+3)}}\alpha(r-i(c+3))\to\infty$$
as $r\to \infty$ since $r-i(c+3)\ge r/2$ and the function $\alpha$ is super-linear by Lemma \ref{al}.

\begin{rk}\label{r:last} The proof of Theorem \ref{t:qi} for $(4,4)$-maps can be easily adapted for $(6,3)$- and $(3,6)$-maps. The ``only if" part does not need virtually any modification.

To prove the "if" part for $(3,6)$-maps by contradiction, again one can use Lemma  \ref{l:qi3} to uniformly bound  the degrees of all faces.
Then one can subdivide non-flat faces by diagonals and
 obtain a quasi-isometric $(3,6)$ map $M'$, where all the faces have degree $3$. If two distinct triangles of $M'$ share an edge $e$, we say that they form a \emph{diamond} with the \emph{hidden} edge $e$.
We can view diamonds as new faces, and build analogs of  corridors made of diamonds, where the gluing edges of a
 diamond are not adjacent. The additional requirement is that the hidden edges of neighbor
 diamonds in a corridor have no common vertices (see Fig. 11). The vertices on sides $\q$ and $\q'$ of a corridor $\cal B$ have degrees at most $4$ in $\cal B$. Then the statement of Lemma \ref{no} holds since
 every exterior vertex (except for one) of the submap $N$ should have degree $\ge 4$. Hence
 one obtains the notion of an infinite and semi-infinite corridors. The Lemma \ref{off} reduces the
 task to a map $M$ with single non-flat vertex, and the rest of the proof of Theorem \ref{t:qi}
 is as above: the quadrangles $e_if_ig_ie_{i+1}$ should be replaced by diamonds and
 the corridors $\cal B$, ${\ccc}_1$, ${\ccc}'_1$, ${\ccc}''_1$,... are now built from
 diamonds. If now one erases the hidden edges of all these diamonds in the quadrant ${\cal Q}_1,\dots$, then the obtained quadrants ${\cal Q}'_i$-s  are $(4,4)$-maps quasi-isometric
 to $Q_i$-s. So our task is reduces to the case of $(4,4)$-maps.

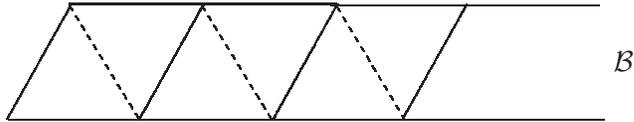
\begin{figure}[ht]
\begin{center}
\unitlength .8mm 
\linethickness{0.4pt}
\ifx\plotpoint\undefined\newsavebox{\plotpoint}\fi 
\begin{picture}(106,23.5)(0,0)
\put(3.75,4.25){\line(1,0){22.25}}
\put(25.75,4.25){\line(1,0){22.25}}
\put(47.75,4.25){\line(1,0){22.25}}
\put(36,23.5){\line(-1,0){22.25}}
\put(58.25,23.5){\line(-1,0){22.25}}
\multiput(14.25,23.5)(-.0336990596,-.0603448276){319}{\line(0,-1){.0603448276}}
\multiput(36.25,23.5)(-.0336990596,-.0603448276){319}{\line(0,-1){.0603448276}}
\multiput(25.5,4.25)(.0336990596,.0603448276){319}{\line(0,1){.0603448276}}
\multiput(47.75,4.25)(.0336990596,.0603448276){319}{\line(0,1){.0603448276}}
\multiput(69.25,4.25)(.0336990596,.0603448276){319}{\line(0,1){.0603448276}}
\put(57.25,23.25){\line(1,0){44.75}}
\put(69.5,4.25){\line(1,0){33.25}}
\multiput(13.93,23.18)(.0333333,-.0543478){15}{\line(0,-1){.0543478}}
\multiput(14.93,21.549)(.0333333,-.0543478){15}{\line(0,-1){.0543478}}
\multiput(15.93,19.919)(.0333333,-.0543478){15}{\line(0,-1){.0543478}}
\multiput(16.93,18.288)(.0333333,-.0543478){15}{\line(0,-1){.0543478}}
\multiput(17.93,16.658)(.0333333,-.0543478){15}{\line(0,-1){.0543478}}
\multiput(18.93,15.028)(.0333333,-.0543478){15}{\line(0,-1){.0543478}}
\multiput(19.93,13.397)(.0333333,-.0543478){15}{\line(0,-1){.0543478}}
\multiput(20.93,11.767)(.0333333,-.0543478){15}{\line(0,-1){.0543478}}
\multiput(21.93,10.136)(.0333333,-.0543478){15}{\line(0,-1){.0543478}}
\multiput(22.93,8.506)(.0333333,-.0543478){15}{\line(0,-1){.0543478}}
\multiput(23.93,6.875)(.0333333,-.0543478){15}{\line(0,-1){.0543478}}
\multiput(24.93,5.245)(.0333333,-.0543478){15}{\line(0,-1){.0543478}}
\multiput(35.93,23.18)(.0333333,-.0543478){15}{\line(0,-1){.0543478}}
\multiput(36.93,21.549)(.0333333,-.0543478){15}{\line(0,-1){.0543478}}
\multiput(37.93,19.919)(.0333333,-.0543478){15}{\line(0,-1){.0543478}}
\multiput(38.93,18.288)(.0333333,-.0543478){15}{\line(0,-1){.0543478}}
\multiput(39.93,16.658)(.0333333,-.0543478){15}{\line(0,-1){.0543478}}
\multiput(40.93,15.028)(.0333333,-.0543478){15}{\line(0,-1){.0543478}}
\multiput(41.93,13.397)(.0333333,-.0543478){15}{\line(0,-1){.0543478}}
\multiput(42.93,11.767)(.0333333,-.0543478){15}{\line(0,-1){.0543478}}
\multiput(43.93,10.136)(.0333333,-.0543478){15}{\line(0,-1){.0543478}}
\multiput(44.93,8.506)(.0333333,-.0543478){15}{\line(0,-1){.0543478}}
\multiput(45.93,6.875)(.0333333,-.0543478){15}{\line(0,-1){.0543478}}
\multiput(46.93,5.245)(.0333333,-.0543478){15}{\line(0,-1){.0543478}}
\multiput(58.43,22.93)(.0326087,-.0543478){15}{\line(0,-1){.0543478}}
\multiput(59.408,21.299)(.0326087,-.0543478){15}{\line(0,-1){.0543478}}
\multiput(60.386,19.669)(.0326087,-.0543478){15}{\line(0,-1){.0543478}}
\multiput(61.364,18.038)(.0326087,-.0543478){15}{\line(0,-1){.0543478}}
\multiput(62.343,16.408)(.0326087,-.0543478){15}{\line(0,-1){.0543478}}
\multiput(63.321,14.778)(.0326087,-.0543478){15}{\line(0,-1){.0543478}}
\multiput(64.299,13.147)(.0326087,-.0543478){15}{\line(0,-1){.0543478}}
\multiput(65.278,11.517)(.0326087,-.0543478){15}{\line(0,-1){.0543478}}
\multiput(66.256,9.886)(.0326087,-.0543478){15}{\line(0,-1){.0543478}}
\multiput(67.234,8.256)(.0326087,-.0543478){15}{\line(0,-1){.0543478}}
\multiput(68.212,6.625)(.0326087,-.0543478){15}{\line(0,-1){.0543478}}
\multiput(69.191,4.995)(.0326087,-.0543478){15}{\line(0,-1){.0543478}}
\put(106,14){\makebox(0,0)[cc]{$\bbb$}}
\end{picture}

\caption{Corridor and hidden edges in $(3,6)$-maps}
\end{center}

\end{figure}

 The case of $(6,3)$-map $M$  can be easily reduces to $(3,6)$. For this goal, one bounds the degrees of faces as above, then chooses a new vertex inside of every face $\Pi$ and connects it with the vertices of $\partial\Pi$. The resulting map is a $(3,6)$-map which is quasi-isometric to $M$, it has finitely many non-flat vertices and faces.

\end{rk}

\section{Maps with angle functions}\label{s:4}

Let $M$ be a map with an angle function (for the definition, see Section \ref{s:i}).

For every face $\Pi$ (vertex $o$) we denote by $\Sigma_\Pi$ (resp. $\Sigma_o$) the sum of the angles of the corners of $\Pi$ (resp. corners at $o$).  Note that if there are no corners at a vertex $o$, then $\Sigma_o=0$.
We define the {\it curvature} $\irr(\Pi)$ of a face $\Pi$ with degree $d=d(\Pi)$ as $\Sigma_\Pi- \pi(d-2))$.
The {\it curvature} $\irr(o)$ of a vertex $o$ is defined as $(2-\mu(o))\pi-\Sigma_o$, where, as before,
$\mu(o)$ is the multiplicity of $o$ in the boundary path of $M$ .

We denote by $I_f$ (by $I_v$) the sum of curvatures of the faces (vertices) of a finite map $M$. The following discrete analog of Gauss - Bonnet formula is well known but we include its proof here  anyway.

\begin{lemma}\label{GB} Let a map $M$ with angle function have at least one edge. Then $I_f+I_v=2\pi$.
\end{lemma}
\proof Let $V, E$ and $F$  be the numbers of vertices, non-oriented edges and faces in $M$
and $n$ be the perimeter of $M$.
It was observed in the proof of Lemma \ref{l:0} that $n=\sum_o\mu(o)$ (the sum over all vertices in $M$).
Since  $\sum_{\Pi} d(\Pi)$ (the sum over all faces in $M$) is equal to the number of exterior edges of the faces in $M$ plus twice
the number of the interior edges in $M$, we have $2E=\sum_{\Pi} d(\Pi)+n = \sum_{\Pi} d(\Pi)+\sum_o \mu(o)$.
Hence $$I_f+I_v = \sum_\Pi (2 -d(\Pi))\pi +\sum_o(2-\mu(o))\pi +
 \left(\sum_{\Pi}\Sigma_\Pi -
\sum_o\Sigma_o\right)$$ $$= \sum_{\Pi} 2\pi+\sum_o 2\pi-\pi\left(\sum_{\Pi}d(\Pi)+ \sum_o \mu(o)\right) +0=2\pi F+2\pi V-2\pi E=2\pi$$
\endproof

In the next lemma, $I_v^i$ (resp., $I_v^b$) is the sum of the curvatures
of interior (exterior) vertices in $M$.

\begin{lemma} \label{A} Let $M$ be a map of perimeter $n\ge 1$ with angle function. Assume that the curvatures of the faces and of the interior vertices
of a map $M$ are non-positive. Then $n\pi\ge -I_f-I_v^i+2\pi$.
\end{lemma}

\proof On the one hand, it follows from the definition that
$$I_v^b=\sum_{o\in\partial M} ((2-\mu(o))\pi-\Sigma_o)\le 2n\pi-n\pi-\sum_{o\in\partial M}\Sigma_o\le n\pi$$
 On the other hand, Lemma \ref{GB} gives us $I_v^b+I_v^i+I_f=2\pi$.
Therefore we have $n\pi+I_f+I_v^i\ge 2\pi$, as required.\endproof

%

Recall that a map $M$ is called $(\delta, b)$-{\em map}
for some $\delta>0$ and a natural number $b > 0$ if

 \begin{itemize}
 \item[(1)]the curvature of every non-flat vertex or face does not exceed $-\delta$ and

 \item[(2)] the degree of every face and of every vertex in $M$ is at most $b$.
 \end{itemize}

We denote by $B(d,o)$ the {\it ball} of radius $d$ centered at $o$ in a graph $G$, i.e.,
the set of vertices $o'$ of $G$ such that $\dist(o',o)\le d$.

The following lemma is well known and obvious.

 \begin{lemma} \label{B} The inequality $|B(d,o)|\le b^d+1$ holds for any graph where degrees of all vertices are at most $b$ (hence for $(\delta, b)$-maps).
 \end{lemma}



 \proof[Proof of Theorem \ref{t:dense}]
Let $\cal V$ be the set of vertices of $M$ which are either exterior or non-flat or belong to a non-flat face.
From the $(\delta,b)$-condition and Lemma \ref{A}, one deduces that

 \begin{equation}\label{O}
|{\cal V}|\le n +\delta^{-1}(-I_v^i) +\delta^{-1}b(-I_f)\le n\pi+\delta^{-1}bn\pi=(\delta^{-1}b+1)n\pi
 \end{equation}
 For arbitrary vertex $o\in \cal V$, we consider the ball $B(o, r).$ By the assumption of the theorem, every vertex $o'$ of $M$ belongs in one of these balls.
 Therefore by Lemma \ref{B}, $\area(M)\le |{\cal V}|\times (b^r+1)\le
 (1+\delta^{-1}b)n\pi(b^r+1)$.  Since every vertex belongs to the boundaries
 of at most $b$ faces, the inequality $\Area(M)\le L n$ follows, provided
 $L\ge  \pi b(1+\delta^{-1}b)(b^r+1)$. \endproof

\begin{rk}\label{last} Theorem \ref{t:dense} generalizes Theorem \ref{t:is}. Indeed, it is enough to establish Theorem \ref{t:is} for simple maps. As it was explained in Subsection \ref{ad}, every simple $(p,q)$-map $M$ can be modified so that the new $(p,q)$-map $M'$ satisfies Condition (B) from Section \ref{ad}. Condition (B) implies  Condition (2) above with $b\ge 11$. Moreover the area of $M'$ is not smaller than the area of $M$, the perimeter of $M'$ is the same as the perimeter of $M$ and the maximal distance from a vertex to an exterior vertex or non-flat vertex or face in $M'$ does not exceed that for $M$.
The $(p,q)$-map $M'$ can be naturally viewed as a map with angle function which assigns to every corner of a $d$-gon face the angle $\frac{\pi (d-2)}{d}$. Again, Condition (B) implies Condition (1) above with $\delta>\frac \pi{21}$. It remains to note that the function $L(r)$ in Theorem \ref{t:dense} is exponential as in Theorem \ref{t:is}.
\end{rk}

\begin{minipage}[t]{3 in}
\noindent Alexander Yu. Ol'shanskii\\ Department of Mathematics\\
Vanderbilt University \\alexander.olshanskiy@vanderbilt.edu\\ and\\ Department of
Higher Algebra, MEHMAT\\
 Moscow State University\\
\end{minipage}
\begin{minipage}[t]{3 in}
\noindent Mark V. Sapir\\ Department of Mathematics\\
Vanderbilt University\\
m.sapir@vanderbilt.edu\\
\end{minipage}

\end{document}